\documentclass[11pt]{article}
\usepackage{latexsym} 
\usepackage{enumitem}
\usepackage{tikz}
\usepackage{rotating}
\input{amssym.def} 
\input{amssym} 
\date{} 
\newfont{\fra}{eufm10 scaled 1095} 
\newfont{\Bb}{msbm10 scaled 1095} 
\newfont{\Bbg}{msbm10 scaled 1280} 
\newcommand\CC{{{\mbox{\Bb C}}}} 
\newcommand\RR{{\mbox{\Bb R}}} 
\newcommand\NN{{\mbox{\Bb N}}} 
\newcommand\ZZ{{\mbox{\Bb Z}}} 
\newcommand\QQ{{\mbox{\Bb Q}}}

\newcommand\fg{{\frak{g}}} 
\newcommand\fh{{\frak h}}

\newcommand\fl{{\frak l}} 
 
\newcommand\fn{{\frak n}} 
\newcommand\fa{{\frak a}} 
\newcommand\fb{{\frak b}}

\newcommand\fk{{\frak k}}

\newcommand\ft{{\frak t}} 
\newcommand\fu{{\frak u}} 
\newcommand\fq{{\frak q}} 
\newcommand\fz{{\frak z}} 
\newcommand\fS{{\frak S}} 
\newcommand\uk{{\underline k}}

\newcommand\cZ{{\cal Z}} 
\newcommand\cH{{\cal H}} 
\newcommand\cM{{\cal M}} 
\newcommand\cN{{\cal N}} 
\newcommand\cA{{\cal A}} 

\newcommand\cP{{\cal P}}  
\newcommand\cQ{{\cal Q}}
\newcommand\cL{{\cal L}}
\newcommand\cF{{\cal F}}
\newcommand\ph{\varphi} 
\newcommand\eps{\varepsilon}

\newcommand\osc{\frak o \frak s \frak c} 
 
\newcommand{\so}{\mathop{{\frak s \frak o}}}
\newcommand{\fo}{\mathop{{\frak o}}}
\newcommand{\fsp}{\mathop{{\frak s \frak p}}}
\newcommand{\End}{\mathop{{\rm End}}} 
\newcommand{\Iso}{\mathop{{\rm Iso}}}
\newcommand{\Aut}{\mathop{{\rm Aut}}} 
\newcommand{\Hol}{{\rm Hol}}
\newcommand{\Der}{\mathop{{\rm Der}}} 
\newcommand{\GL}{\mathop{{\rm GL}}} 
\newcommand{\Sp}{\mathop{{\rm Sp}}} 
\newcommand{\SL}{\mathop{{\rm SL}}} 
\newcommand{\SO}{\mathop{{\rm SO}}} 
\newcommand{\PO}{\mathop{{\rm PO}}} 
\newcommand{\PGL}{\mathop{{\rm PGL}}} 
\newcommand{\grO}{{\rm O}}
\newcommand{\Hom}{\mathop{{\rm Hom}}} 
\newcommand{\U}{{{\rm U}}}

\newcommand{\id}{\mathop{{\rm id}}} 
\newcommand{\ad}{\mathop{{\rm ad}}} 
\newcommand{\tr}{\mathop{{\rm tr}}} 
 
\newcommand{\Ad}{\mathop{{\rm Ad}}}

\newcommand{\Ker}{\mathop{{\rm ker}}} 
\newcommand{\im}{\mathop{{\rm im}}} 

\newcommand{\diag}{\mathop{{\rm diag}}} 
 
\newcommand{\Span}{\mathop{{\rm span}}} 
 
\newcommand{\ind}{{{\rm ind}}} 
\newcommand{\reg}{{\mathop{{\rm reg}}}}
\newcommand{\spec}{{\mathop{{\rm spec}}}}
\newcommand{\ms}{{\}_{_!}}}
\newcommand{\pro}{{\rm pr}} 
\newcommand{\pv}{\mbox{{\small $\oplus$}}}
\newcommand{\mv}{\mbox{{\small $\ominus$}}}
\newcommand\ip{{\langle\cdot \,,\cdot \rangle}}

\newcommand\proof{{\sl Proof. }} 
\newcommand{\qed}{\hspace*{\fill}\hbox{$\Box$}\vspace{2ex}} 
\newcommand{\qedohne}{\hspace*{\fill}\hbox{$\Box$}}
\newcommand{\benur}{\begin{enumerate}[label=(\roman*)]}
\newtheorem{theo}{Theorem}[section] 
\newtheorem{pr}[theo]{Proposition} 
\newtheorem{de}[theo]{Definition} 
\newtheorem{ex}[theo]{Example} 
\newtheorem{re}[theo]{Remark} 
\newtheorem{co}[theo]{Corollary} 
\newtheorem{lm}[theo]{Lemma} 
\begin{document} 
\title{Compact quotients of Cahen-Wallach spaces}
\author{Ines Kath and Martin Olbrich} 
\maketitle 
\begin{abstract}
\noindent Indecomposable symmetric Lorentzian manifolds of non-constant curvature are called Cahen-Wallach spaces. Their isometry classes are described by continuous families of real parameters.  
We derive necessary and sufficient conditions for the existence of compact quotients of Cahen-Wallach spaces in terms of these parameters.
\end{abstract}
\tableofcontents
\section{Introduction}
A Clifford-Klein form of a homogeneous space $X= G/H$ of a Lie group $G$ is the quotient
manifold $\Gamma\setminus X$, where $\Gamma \subset G$ is a discrete subgroup of $G$ acting properly and freely on $X$.  This paper will deal with compact Clifford-Klein forms, which will also be called compact quotients. The existence or non-existence of  compact Clifford-Klein forms of a given homogenuous space has been the subject of intense study for many years.  If the isotropy subgroup $H$ is compact, then the action of any discrete subgroup of $G$ on $X$ is proper. Hence the study of compact quotients of $X$ is essentially equivalent to the study of cocompact lattices of $G$. 
On homogeneous spaces with non-compact stabiliser, the action of a discrete group $\Gamma\subset G$ on $X$ is not automatically proper and the existence of compact quotients is a far more involved problem. Properly discontinuous actions, especially those on reductive spaces, were intensely studied for instance by Benoist and Kobayashi, see, e.g., \cite{Be, Ko1,Ko2, KY}.  The homogeneous spaces that will be considered in the present paper have a non-compact stabiliser and are, moreover, non-reductive.

Of particular interest is the situation where $G$ preserves some kind of geometry (affine, conformal, pseudo-Riemannian, etc.) on $X$. The study of Clifford-Klein forms fits into the more general concept of geometric structures on manifolds locally modelled on a homogeneous space $X=G/H$ of a Lie group $G$. These are called $(G,X)$-structures. For a review on what is known on such structures see \cite{Gicm}. In particular, the existence of compact manifolds carrying a $(G,X)$-structure for certain groups $G$ is discussed. 

Let us now assume that $X$ is a symmetric space. Riemannian symmetric spaces $X$ have a compact stabiliser and admit compact Clifford-Klein forms \cite{Bo}. For non-Riemannian symmetric spaces, in general, the stabiliser is non-compact and it becomes rather difficult to prove existence or non-existence of compact quotients, see \cite{KY} for a review. Even for Lorentzian symmetric spaces no complete answer is known. It is the aim of the paper to shed more light to the existence of compact quotients in this special situation, i.e., we want to discuss the 

{\bf Problem:} Which Lorentzian symmetric spaces admit compact Clifford-Klein forms?   

More exactly, we want to consider Lorentzian symmetric spaces $X=G/G_+$, where $G$ is the isometry group of $X$, which can essentially differ from the transvection group of $X$.   
Let us first review what is known for Lorentzian symmetric spaces $X$ of constant sectional curvature.

For positive sectional curvature the Calabi-Markus phenomenon occurs: Every subgroup of the isometry group of the de Sitter spacetime $S^{1,n}$, $n\ge2$, that acts properly discontinous on $S^{1,n}$ is finite \cite{CM}. Hence compact quotients $\Gamma\setminus S^{1,n}$ of the de Sitter spacetime $S^{1,n}$ do not exist \cite{CM}. 

Kulkarni \cite{Ku} proved that compact quotients $\Gamma\setminus \tilde H^{1,n}$ of the universal anti de Sitter spacetime $\tilde H^{1,n}$ are odd-dimensional. Moreover, he showed that for each odd dimension such a quotient exist. He used that the group $\U(1,m) \subset\SO(2,2m)$ acts isometrically on $H^{1,2m}$. This action is transitive and proper. Hence every torsion-free lattice $\Gamma\subset \U(1,m)$ defines a compact quotient.  Any Lorentzian manifold that is obtained, up to finite coverings, by this construction is called standard. It is conjectured that in dimension fuer $2m+1>5$ all quotients of $\tilde H^{1,2m}$ are standard \cite{Zads}. In dimension three there exist non-standard quotients \cite{Gh, Gads}.
 
If $Y$ is a compact quotient of the Minkowski space $\RR^{1,n}$, then there is a connected solvable group $U$ acting isometrically and simply transitively on $\RR^{1,n}$ and a lattice $\Gamma\subset U$ such that $Y=\Gamma\setminus \RR^{1,n}$ up to finite coverings \cite{GK}.

Besides these spaces, which have a reductive transvection group, there exist many non-reductive Lorentzian symmetric spaces. These spaces were classified by Cahen and Wallach \cite{CW}.
Each indecomposable Lorentzian symmetric space is either semisimple or solvable.  A non-flat simply-connected indecomposable solvable Lorentzian symmetric space $X$ is called Cahen-Wallach space. Any such space is isometric to some $X_{p,q}(\lambda,\mu):=(\RR^n, g_{\lambda,\mu})$, where $$g_{\lambda,\mu}=2dzdz'+ \sum_{i=1}^{p+q} dx_i^2 + \Big(  \sum_{i=1}^p\lambda_i^2x_i^2 -  \sum_{j=1}^q\mu_j^2x_{p+j}^2 \Big)dz'^2$$ for parameters $(\lambda,\mu)\in(\RR^*)^p\times(\RR^*)^q$, $n=2+p+q$. We will say that $X$ is of real type if $q=0$, of imaginary type if $p=0$ and of mixed type if $p,q\not =0$. For spaces of real or imaginary type, we write just $X_{p,0}(\lambda)$ and $X_{0,q}(\mu)$ instead of $X_{p,q}(\lambda,\mu)$.

Our aim is to find conditions for the parameters $(\lambda,\mu)$ that are equivalent to the existence of a compact quotient of $X$. To our knowledge this question has not been much studied before. 
For certain choices of $(\lambda,\mu)$, Cahen and Wallach \cite{CW} claim to construct examples of compact quotients.  However, the action of the discrete subgroups they consider is not proper, see Remark~\ref{CWfalsch}. 

Before we will start to explain our results we want to recall a proven approach to this kind of problems. To find compact quotients of a homogeneous space $X=G/H$ one can try to find a (virtually) connected subgroup $U\subset G$ acting properly and cocompactly (or even transitively) on $X$ and a cocompact lattice $\Gamma$ in $U$. For instance, in the already mentioned paper \cite{Ku}, Kulkarni did this for pseudo-Riemannian space forms. 
For some homogeneous spaces $X=G/H$ all compact quotients are of this kind. This observation was one of the key points in the classification of three-dimensional affine crystallographic groups by Fried and Goldman \cite{FG}. They proved that every subgroup $\Gamma\subset {\rm Aff}(E)$, $\dim E=3$, acting properly discontinuously and cocompactly on $E$ is virtually solvable. Then, in order to find a suitable  group $U$,  they proved the existence of a syndetic hull  for virtually solvable subgroups of $\GL(n, \RR)$. A syndetic hull of a closed subgroup $\Gamma\subset G$ is defined to be a connected subgroup
$S$ of $G$ containing $\Gamma$ such that  $\Gamma\setminus S$ is compact.

For Cahen-Wallach spaces, we will proceed in a similar way. That is, we will show that $\Gamma$ is, essentially, a lattice in a certain closed subgroup of $G$, which, however, now can have an infinite cyclic component group. Let us state this in a slightly more precise way. The transvection group $\hat G$ of a Cahen-Wallach space $X$ is 
isomorphic to a semi-direct product $H_n\rtimes \RR$  for some Heisenberg group $H_n$, and the isometry group of $X$ is a semi-direct product $\hat G\rtimes (K\times\ZZ_2)$, where $K$ is compact. We may consider $G=\hat G\rtimes K$ instead of the whole isometry group. In Section~\ref{S3} we study discrete subgroups of $G$. We slightly generalise the setup to groups  $G=N\rtimes(\RR \times K)$, where $N$ is simply-connected nilpotent, $K$ is compact, and $\RR\times K$ acts by semisimple automorphisms on $N$. For a discrete subgroup $\Gamma\subset G$, let $\Delta\subset\RR$ be the closure of the projection of $\Gamma$ to the $\RR$-factor of $G$. We show that $\Gamma$ is, essentially,  
a lattice in a subgroup $(U\cdot\psi(\Delta))\times C_K\subset G$, where $\psi:\Delta\rightarrow G$ is a section of the projection to $\RR$, $U\subset N$ is connected and $C_K\subset  K$ is connected and abelian. This generalises a classical result of Auslander \cite{Aus} on discrete subgroups of semidirect products of nilpotent with compact Lie groups.

Let us return to discrete subgroups $\Gamma\subset G=H_n\rtimes (\RR\times K)$ of isometries of a Cahen-Wallach space $X$. If $\Gamma$ acts properly and cocompactly on $X$, then obviously $\Delta\not=\{0\}$, thus $\Delta$ is infinite cyclic or equal to $\RR$. Both cases are possible. However, the case $\Delta=\RR$ is very special and can occur only if $X$ is a Lie group with biinvariant Lorentian metric, see Thm.~\ref{AA}. All these groups also admit a compact quotient $\Gamma\setminus X$ with infinite cyclic $\Delta=\langle t_0\rangle$. Hence, if we are only interested in conditions characterising the existence of a compact quotient we may concentrate on the case $\Delta=\langle t_0\rangle$. Prop.~\ref{wulf} gives a first criterion for the existence of such compact quotients. Its proof relies on the fact that a discrete subgroup $\Gamma\subset G$ acts properly and cocompactly on $X$ if and only if $U\cdot\psi(\Delta)$ does so. In Sections~\ref{neukirch}--\ref{hans-werner} we try to make this criterion as explicit as possible in terms of the parameters $(\lambda,\mu)$ of $X$.

Let us first consider a Cahen-Wallach space $X$ of real type. Then Theorem~\ref{BB}  gives a purely arithmetic criterion for the existence of compact quotients:  $X$ admits a compact quotient if and only if there exists a polynomial
$f\in\ZZ[x]$ of the form
$f(x)= x^n+a_{n-1}x^{n-1}+\dots+a_1x\pm 1$
with no roots on the unit circle such that 
$X\cong X_{n,0}(\log |\nu_1|,\log |\nu_2|,\dots,\log |\nu_n|) ,$
where $\nu_1,\nu_2,\dots,\nu_n$ are the roots of $f$. Obviously, Theorem~\ref{BB} yields a `recipe' to find all spaces of real type admitting compact quotients. However, for a {\it given} Cahen-Wallach space $X_{n,0}(\lambda)$ it might be rather difficult to decide whether this condition is satisfied.  One of the difficulties is caused by the fact that $X_{n,0}(\lambda)\cong X_{n,0}(t\lambda)$ for all $t\in\RR^*$ whereas the condition for the parameters $\log |\nu_1|,\log |\nu_2|,\dots,\log |\nu_n|$ in Theorem~\ref{BB} is not scaling invariant. There are, however, some necessary conditions, which are easy to check, see Propositions~\ref{transc}, \ref{werwolf}, and the remark on the trace condition below.  The investigation of the moduli space of isometry classes of Cahen-Wallach spaces  admitting compact quotients using Theorem~\ref{BB} leads to classical problems in number theory. For instance, a complete description of the moduli space for $n=3$ depends on whether the four exponentials conjecture is true. 

In order to formulate the result for spaces of imaginary type, let us introduce the notion of $\RR$-admissibility. Given a $d$-tuple $\uk=(k_1,\dots,k_d)\in\ZZ^d$ we define a linear map $L_\uk:\CC^d\rightarrow \CC^d$ by $L(z_1,\dots,z_m)=i(k_1 z_1,\dots, k_d z_d)$. We will say that $\uk$ is $\RR$-admissible if there exists a real $d$-dimensional subspace $V\subset \CC^d$ such that $\exp(tL)(V)\cap\RR^d=0$ for all $t\in\RR$. We will prove that a Cahen-Wallach space  $X$ of imaginary type admits a compact quotient if and only if there exists an $\RR$-admissible $d$-tuple $(k_1,\dots,k_d)\in(\ZZ_{\not=0})^d$ such
that 
$X\cong X_{0,n}(k_1,\dots, k_d, \mu_{d+1},\dots,\mu_n)\,,$
where the parameters $\mu_i\in\RR^*$, $i=d+1,\dots,n$, all appear with even multiplicity, see Theorem~\ref{CC}. This reduces the problem to an elementary geometric one on rotations in $\CC^d$,  which is closely related to problems on positive trigonometric polynomials,  on the matrix Riccati equation and on the topology of Grassmannians. Although its formulation is simple, it seems to be not easy to find a general solution. We can prove that if $\uk$ is $\RR$-admissible, then $\sum_{i=1}^{d} c_ik_i = 0$ for suitable choice of $c_i\in\{1,-1\}$.  For $n\le 4$, we show that this condition is also sufficient. It is an open question whether it is sufficient for $n\ge 5$, too.

Theorem~\ref{DD} gives a criterion for the existence of compact quotients for spaces of general type. We do not want to give an exact formulation here, the conditions are a nontrivial combination of those for spaces of purely real and of purely imaginary type. However, we want to mention the following necessary trace condition.
If $X_{p,q}(\lambda,\mu)$ admits a compact quotient, then 
$ \sum_{i=1}^{p} c_i\lambda_i = 0$ and $\sum_{j=1}^{q} \hat c_j\mu_j = 0$ for suitable $c_i,\hat c_j\in\{1,-1\}$.

Although the paper concentrates on criteria for the existence of compact quotients rather than on a systematic study of all these quotients we want to remark that our proofs are constructive and yield explicit examples of compact quotients.  Moreover, they 
contain  information on the shape of a discrete group $\Gamma$ defining a compact quotient $\Gamma\setminus X$, especially on the structure of $\Gamma$ when considered as an abstract group.
Roughly speaking, $\Gamma$ contains a subgroup of finite index that is a semi-direct product of a (possibly degenerate) discrete Heisenberg group and  $\ZZ$. The type of the Heisenberg group and the action of $\ZZ$ on it is described in Prop.~\ref{fundi}. In particular, we see that $\Gamma$ is never abelian.

There are compact quotients $\Gamma\setminus X$ of Cahen-Wallach spaces for which $\Gamma$ is not only contained in the isometry group of $X$ but even in the transvection group. These are exactly the quotients whose holonomy group is abelian. We will decide for which parameters $(\lambda,\mu)$ the space $X_{p,q}(\lambda,\mu)$ admits such a quotient in Subsection~\ref{trans}.
        
Subsection~\ref{solv} deals with compact manifolds of the form $\Gamma\setminus S$  called Lorentzian solvmanifolds. Here $S$ is a 1-connected solvable Lie group equipped with a left-invariant Lorentzian metric and $\Gamma \subset S$ is a lattice.  We decide which Cahen-Wallach spaces have compact quotients that are solvmanifolds. Essentially, this leads to the question for which compact quotients $\Gamma\setminus X$ there 
is a connected subgroup $S \subset G$ containing $\Gamma$ as a lattice (i.e., a syndetic hull).  Note that this problem is not yet solved by the above mentioned construction of the group $(U\cdot\psi(\Delta))\times C_K\subset G$ containing  $\Gamma$ as a lattice since in general this group is not connected.

In Subsection~\ref{low}, we will give a rather explicit description of moduli spaces of low-dimensional Cahen-Wallach spaces admitting compact quotients.

The results of this paper provide a basis for future investigations concerning problems as: (1) the existence of compact quotients of decomposable Lorentzian symmetric spaces, that is, of products of Cahen-Wallach spaces by flat or semisimple Riemannian ones, (2) the classification of all compact quotients for a given Cahen-Wallach space, (3) the determination of deformation spaces of compact quotients.

Finally, note that compact quotients of Cahen-Wallach spaces are just the same as compact manifolds that are locally isometric to a Cahen-Wallach space and geodesically complete. So, a natural question is whether there exist incomplete compact manifolds locally isometric to a Cahen-Wallach space. New results by Leistner and Schliebner \cite{LS} say that this is not the case. Actually, they proved that any compact $pp$-wave is complete. This generalises and complements results on the completeness of compact Lorentzian manifolds of constant curvature by Carri\`ere and Klingler \cite{Ca, K}.

\subsection*{Some conventions} 
\begin{tabbing}
\hspace*{2cm}\=\kill
$\NN$ \> set of positive integers \\ 
$\NN_{0}$\> $\NN\cup \{0\}$ \\
$\ZZ_{\not=0}$\>$\ZZ\setminus\{0\}$\\
$\QQ^+$\> set of positive rational numbers\\
$\RR^+$\> set of positive real numbers\\
$\RR^*$ \> $\RR\setminus\{0\}$\\
$G_0$ \> identity component of a group $G$\\
$Z(G)$\> center of a group $G$\\
$N_G(U)$ \> normaliser of a subgroup $U$ of a group $G$\\
$M^a$\> set of fixed points of a map $a:M\rightarrow M$\\
$M^G$ \> set of fixed points for an action of a group $G$ on a space $M$
\end{tabbing}
\section{Cahen-Wallach spaces}
\subsection{Classification}\label{class}
Let us recall the construction and classification of Cahen-Wallach spaces. We will use the description of simply connected symmetric spaces by their associated infinitesimal objects called symmetric triples. A symmetric triple $(\hat \fg,\theta,\ip)$ consists of a Lie algebra $\hat \fg$, an indefinite scalar product $\ip$ on $\hat \fg$ and an involutive automorphism $\theta$ of $\hat \fg$, which satisfy the following conditions.  The scalar product $\ip$ is invariant under the adjoint representation of $\hat \fg$, $\theta$ is an isometry with respect to $\ip$ and the eigenspaces $\hat \fg_+$ and $\hat \fg_-$ of $\theta$ with eigenvalues $1$ and $-1$ satisfy $[\hat \fg_-,\hat \fg_-]=\hat \fg_+$. The correspondence between a simply connected symmetric space $X$ and the associated symmetric triple $(\hat \fg,\theta,\ip)$ is given in such a way that $\hat \fg$ is the Lie algebra of the transvection group $\hat G$ of $X$. The connected Lie subgroup $\hat G_+\subset \hat G$ with Lie algebra $\hat \fg_+$ is the stabiliser of a fixed base point $x_0\in X$. Moreover, $\hat \fg_-$ can be identified with the tangent space of $X$ at $x_0$  and $\ip|_{\hat \fg_-\times \hat \fg_-}$ is the metric of $X$ in $x_0$. We have used the somewhat unusual notation $\hat G$ for the transvection group since we want to reserve the notation $G$ for (a subgroup of index 2 of) the isometry group.

Cahen-Wallach spaces as introduced in the introduction are associated with non-abelian indecomposable solvable Lorentzian symmetric triples, which we want to call Cahen-Wallach triples. Such triples can be constructed in the following way. 
Let $\omega$ be a non-degenerate 2-form on $\RR^{2n}$. The $(2n+1)$-dimensional Heisenberg group $H_n(\omega)$ is a central extension of the abelian Lie group $\fa:=\RR^{2n}$ by $\fz:=\RR$ defined by 
\begin{equation}
(z,a)\cdot (z', a')=(z+z'+\textstyle{\frac12} \omega(a,a'), a+ a') \label{E*Hn}
\end{equation}
for $z,z' \in\frak z$ and $a,a'\in\fa$. The isomorphism class of $H_n(\omega)$ does not depend on  $\omega$ and we just write $H_n$ or $H$ instead of $H_n(\omega)$ if we are not interested in the explicit realisation of this group.
The Lie algebra of $H_n(\omega)$ equals $\fh_n(\omega):=\fz\oplus\fa$ (as a vector space) with Lie bracket
$$[(z,a),(\hat z, \hat a)] =(\omega(a,\hat a),0).$$

Now let $\theta_\fa:\fa\rightarrow\fa$ be a linear map such that $\theta_\fa^2=\id$ and $\theta_\fa^*\omega=-\omega$. Then the eigenspaces $\fa_+$, $\fa_-$ of $\theta_\fa$ are Lagrange spaces of $\omega$. Moreover, let $L:\fa\rightarrow \fa$ be an invertible linear map such that 
\benur
\item $\omega(La,a')+\omega(a,La')=0$ for all $a,a'\in\fa$,
\item $L\circ\theta_\fa =-\theta_\fa\circ L$,
\item $ \ip_\fa:=\omega(L^{-1} \cdot\,,\cdot)$ restricted to the $(-1)$-eigenspace $\fa_-$ of $\theta_\fa$ is positive definite.
\end{enumerate}
Then $L$ defines actions $l$ and $l_*$ of $\RR$ on $H=H_n(\omega)$ and $\fh:=\fh_n(\omega)$, respectively, by  
\begin{eqnarray} l: \RR\longrightarrow\Aut(H),&& l(t)(z,a) = (z, e^{tL} a) \label{El}\\
l_*: \RR\longrightarrow\Der(\fh),&& l_*(t)(z,a) = (0, tL a). \nonumber
\end{eqnarray}
Let us consider the semidirect product $\hat G:=H\rtimes_l\RR$ with Lie algebra $\hat \fg:=\fh\rtimes_{l_*}\RR$. We will write also $h$ for $(h,1)\in \hat G$,  $t$ for $(0,t)\in \hat G$ and $h\cdot t$ instead of $(h,t)$ for $h\in H$ and $t\in \RR$. If $\ip_\fa$ is positive definite on the whole vector space $\fa$, then $\hat G$ is called oscillator group. Otherwise we will call it generalised oscillator group.

We define an involution $\theta$ on $\fz\times \fa\times \RR$ by 
$$\theta(z,a,t)=(-z, \theta_\fa(a),-t)$$ 
for $(z,a,t)\in\fz\times \fa\times \RR$. Then $\theta$ is an automorphism of the Lie group $\hat G$ as well as an automorphism of the Lie algebra $\hat \fg$. Obviously, $\theta\in \Aut(\fg)$ is the differential of $\theta\in \Aut(\hat G)$.

Finally, we define an indefinite scalar product $\ip$ on $\hat \fg$ by
$$ \fz\perp\fz\oplus\fa, \ \fa\oplus\RR\perp\RR,\ \ip|_{\fa\times\fa} =\ip_\fa,\ \langle z,t\rangle =zt$$
for $z\in\fz$ and $t\in\RR$. Then $\ip$ is invariant under the adjoint representation of $\hat \fg$.
Moreover, $\theta\in \Aut(\fg)$ is an isometry with respect to $\ip$. 
In this way we obtain an indecomposable solvable Lorentzian symmetric triple
$\osc(\omega,\theta_\fa,L):= (\hat \fg,\theta,\ip)$.  

Let $X$  be the symmetric space that is associated with $\osc(\omega,\theta_\fa,L)$.  The transvection group of $X$ is isomorphic to $\hat G= H\rtimes_{l}\RR$, where $H=H_n(\omega)$ and $l$ is defined as in (\ref{El}). Moreover, $X$ is the homogeneous space $\hat G/\hat G_+$, where $\hat G_+$ is the connected subgroup of $\hat G$ whose Lie algebra equals $\hat \fg_+=\fa_+\subset \hat \fg$. Note that the exponential map $\exp:\fh \rightarrow H$ is equal to the identity on $\fz\oplus\fa$ under the identifications $\fh=\fz\oplus\fa$ (as vector spaces) and $H=\fz\oplus\fa$ (as manifolds). Thus $\hat G_+=\fa_+\subset \hat G$.
The symmetric space $X$ can be identified with $\hat \fg_-$ by 
\begin{eqnarray*} X=\hat G/\hat G_+ & \longrightarrow & \fz\times\fa_-\times\RR\\
t\cdot(z,a)\cdot \hat G_+& \longmapsto & (z+\textstyle{\frac12}\omega(a_+,a_-),a_-,t)
\end{eqnarray*}
for $t\in\RR$, $(z,a)\in H$, where $a=a_++a_-$ for $a_+\in\fa_+$ and $a_-\in\fa_-$. Using that the embedding $X=\fz\times\fa_-\times\RR\hookrightarrow \hat G$ is a section of the projection $\hat G\rightarrow X$ it is easy to see that with this identification the metric of $X$ at $(z,a,t)\in\fz\times\fa_-\times\RR$ equals
$$2dzdt+\ip_\fa|_{\fa_-\times\fa_-}-\langle La,La\rangle_\fa \cdot dt^2.$$

Cahen and Wallach proved that every non-abelian indecomposable solvable Lorentzian symmetric triple is isomorphic to some $\osc(\omega,\theta_\fa,L)$ for suitable data $\fa$, $\omega$, $\theta_\fa$ and $L$. Moreover, two symmetric triples $\osc(\omega,\theta_\fa,L)$ and $\osc(\hat\omega,\theta_{\hat\fa},\hat L)$ are isomorphic if and only if there is an orthogonal map $A:(\fa_-,\ip_\fa|_{\fa_-\times\fa_-})\rightarrow(\hat\fa_-,\ip_{\hat\fa}|_{{\hat\fa_-}\times{\hat\fa_-}})$ and a real number $c>0$ such that $cAL^2A^{-1}=\hat L^2$ on $\hat \fa_-$ \cite{CW}. This is the case if and only if the eigenvalues of the symmetric maps $L^2|_{\fa_-}:\fa_-\rightarrow \fa_-$ and $\hat L^2|_{\hat\fa_-}:\hat\fa_-\rightarrow \hat\fa_-$ coincide up to a common positive factor. 

In particular, the number $p$ of positive eigenvalues and the number $q$ of negative eigenvalues of $L^2$ on $\fa_-$ are invariants of the isomorphism class of $\osc(\omega,\theta_\fa,L)$. We will call $(p,q)$ the type of $\osc(\omega,\theta_\fa,L)$. If $X$ is a Cahen-Wallach space associated with $\osc(\omega,\theta_\fa,L)$ we will also say that $X$ is of type $(p,q)$. 
If $\lambda_1^2,\dots,\lambda_p^2$, $-\mu_1^2,\dots,-\mu_q^2$ are the  eigenvalues of $L^2$ on $\fa_-$, then $ \pm\lambda_1,\dots,\pm \lambda_p$ and  $\pm i \mu_1,\dots,\pm i\mu_q$ are the eigenvalues of $L$ considered as a complex linear map on the complexification $\fa_{\Bbb C}$ of $\fa$. This motivates the following definition. 

\begin{de}
Let $X$ be a Cahen-Wallach space of type $(p,q)$. We will say that $X$ is of real type if $q=0$, of imaginary type if $p=0$ and of mixed type if $p>0$ and $q>0$. 
\end{de}

Let $\cM_{p,q}$ denote the set of isomorphism classes of Cahen-Wallach triples of type $(p,q)$. We will denote the set of isometry classes of the associated Cahen-Wallach spaces by the same symbol. 
The classification explained above gives us a surjection 
$$
\Phi_{p,q}:(\RR^*)^{p+q}\longrightarrow \cM_{p,q},
$$
where a symmetric triple $\osc(\omega,\theta_\fa,L)$ for which $\lambda_1^2,\dots,\lambda_p^2$, $-\mu_1^2,\dots,-\mu_q^2$ are the  eigenvalues of $L^2$ on $\fa_-$ belongs to the isomorphism class $\Phi_{p,q}(\lambda,\mu)$ for $\lambda=(\lambda_1,\dots,\lambda_p)$, $\mu=(\mu_1,\dots,\mu_q)$. We want to describe the fibres of $\Phi_{p,q}$.
Let ${\fS}_m$ denote the symmetric group of degree $m$. The group $\bar{\fS}_m:={\fS}_m\ltimes (\ZZ_2)^m$ acts on $(\RR^*)^m$ by 
$$(\sigma,\kappa)\cdot x= (\kappa_1 x_{\sigma(1)},\dots,\kappa_m x_{\sigma(m)})$$ for $\sigma\in\fS_m,\ \kappa=(\kappa_1,\dots,\kappa_m)\in(\ZZ_2)^m$ and $x=(x_1,\dots,x_m)\in\RR^m$. 
We define an action of $\RR^*\times \bar{\fS}_p\times \bar {\fS}_q$ on $(\RR^*)^{p+q}$ by 
$$(r,s_p,s_q)\cdot (\lambda,\mu)=(r\cdot s_p\cdot\lambda, r\cdot s_q\cdot\mu) $$
for $r\in\RR^*$, $s_p\in\bar{\fS}_p$,  $s_q\in\bar{\fS}_q$, $\lambda\in(\RR^*)^p$ and $\mu\in(\RR^*)^q$.
Then the fibres of $\Phi_{p,q}$ are exactly the orbits of this action. 

We endow $\cM_{p,q}$ with the quotient topology with respect to $\Phi_{p,q}$.

One can obtain an alternative description of $\cM_{p,q}$ using the  bijection
\begin{equation}\label{humus} 
\left\{ (\lambda,\mu)\left| 
\begin {array}{l}\ \lambda=(\lambda_1,\dots,\lambda_p)\in\RR^p,\ \mu=(\mu_1,\dots,\mu_q)\in\RR^q,\\
 0< \lambda_1\le\lambda_2\le\dots\le\lambda_p,\  1= \mu_1\le\mu_2\le \dots\le \mu_q \end{array} \right.\right\}\ \longrightarrow \cM_{p,q}
\end{equation}
sending $(\lambda,\mu)$ to $\Phi_{p,q}(\lambda,\mu)$ if $p,q>0$. In fact, if the left hand side carries the topology inherited by $\RR^p\times\RR^q$,
then it is a homeomorphism. Similarly, we have homeomorphisms
\begin{eqnarray*} 
\left\{ \lambda\in(\RR^*)^p\mid   
1=\lambda_1\le\lambda_2\le\dots\le\lambda_p\right\}& \longrightarrow & \cM_{p,0},\\[1ex]
\left\{ \mu\in(\RR^*)^q\mid 
1=\mu_1\le\mu_2\le\dots\le\mu_q\right\}& \longrightarrow & \cM_{0,q}.
\end{eqnarray*}
For each element of $\cM_{p,q}$, we are going to provide an explicit description of some representatives of the isomorphism/isometry class. Before we start, let us define certain endomorphisms of $\RR^{2m}\cong\CC^m$ and $\RR^{4m}\cong\CC^{2m}$, which we will use frequently in this paper.
\begin{de} For $\mu=(\mu_1,\dots,\mu_m)\in\RR^m$, we set
\begin{eqnarray}
&L_\mu:\, \CC^m\longrightarrow \CC^m,\ \ &\ L(z_1,\dots,z_m)=i(\mu_1 z_1,\dots, \mu_m z_m) \label{L}\\
&\phi_\mu:\, \CC^{2m}\longrightarrow \CC^{2m},&\ \phi(z_1,\dots,z_m)=(-\mu_1 z_2,\mu_1z_1,\dots,-\mu_mz_{2m}, \mu_m z_{2m-1})\,.\label{phi}
\end{eqnarray} 
\end{de}

\begin{ex} \label{Exosc}{\rm This example will show that, for fixed $n\in\NN$, we can obtain representatives of all elements of $\cM_{p,q}$ with $p+q=n$ starting from the same data $\fa$, $\theta_\fa$ and $\omega$ by varying $L$.
We consider $\fa=\RR^{2n}\cong \CC^n$ and the standard symplectic form on $\fa$ defined by $\omega(a,a')=\Im(\overline a^\top\cdot a')$.
Let $\theta_\fa$ be the complex conjugation on $\fa$.  We fix elements $\lambda\in\RR^p$ and $\mu\in\RR^q$ with {\it positive} coordinates and define an endomorphism $L$ on $\fa=\CC^p\oplus\CC^q$ by
\begin{equation}\label{EdefL}
L=(L_\lambda\circ \theta_p) \oplus L_\mu\,,
\end{equation}
where $\theta_p$ is the complex conjugation on $\CC^p$. 
Then $\omega,\ \theta_\fa$ and $L$ define an indecomposable solvable Lorentzian symmetric triple $\osc(\omega,\theta_\fa,L)$, which is a representative of the isomorphism class $\Phi_{p,q}(\lambda,\mu)\in \cM_{p,q}$.
}\end{ex}
\begin{ex} \label{Exosc2}{\rm In this example we will give another description of representatives of the isomorphism classes of Cahen-Wallach triples. This description is adapted to the notation in \cite{KO1} and \cite{KOesi}, where symmetric triples were constructed by quadratic extensions of a Lie algebra with involution $(\fl,\theta_\fl)$ by an orthogonal $(\fl,\theta_\fl)$-module $(\fa,\theta_\fa,\ip_\fa)$. This construction is closely related to double extensions introduced by Medina and Revoy \cite{MR1}. For $p,q\in\NN$, we choose $\fa=\RR^{2p}\oplus\RR^{2q}\cong\CC^p\oplus\CC^q$ and consider again the complex conjugation $\theta_\fa$ on $\fa$. On $\RR^{2p}\cong\CC^p$ we define a scalar product $\ip_{p,p}$ of signature $(p,p)$ by $\langle u,v\rangle_{p,p}= -\Re(u^\top v)$ and on $\RR^{2q}$ we consider the Euclidean standard scalar product $\ip_{2q}$ given by $\langle u,v\rangle_{2q}=\Re (u^\top \bar v)$. Now we fix the scalar product $\ip_\fa:=\ip_{p,p}\oplus\ip_{2q}$ on $\fa$.  For arbitrary $\lambda\in (\RR^*)^p$ and $\mu\in(\RR^*)^q$, we define $L$ as in (\ref{EdefL}). The map $L$ is antisymmetric with respect to $\ip_\fa$ and we put $\omega:=\langle L\,\cdot\,,\cdot\rangle_\fa$. Hence, $\omega,\ \theta_\fa$ and $L$ define an indecomposable solvable Lorentzian symmetric triple $\osc(\omega,\theta_\fa,L)$, which is also denoted by $\osc_{p,q}(\lambda,\mu)$. Analogously we define $\osc_{p,0}(\lambda)$ and $\osc_{0,q}(\mu)$. Then $\osc_{p,q}(\lambda,\mu)$ is in the isomorphism class $\Phi_{p,q}(\lambda,\mu)\in\cM_{p,q}$.
Analogous statements hold for $\osc_{p,0}(\lambda)$ and $\osc_{0,q}(\mu)$.

}\end{ex}
\begin{de}
For $\lambda\in (\RR^*)^p$ and $\mu\in(\RR^*)^q$, where $p,q\in\NN$, let $X_{p,q}(\lambda,\mu)$, $X_{p,0}(\lambda)$ and $X_{0,q}(\mu)$ denote the Cahen-Wallach spaces associated with $\osc_{p,q}(\lambda,\mu)$, $\osc_{p,0}(\lambda)$ and $\osc_{0,q}(\mu)$, respectively. We will also use the notation $X_{p,q}(\lambda,\mu)$ if either $p$ or $q$ equals zero. In this case, $\lambda$ is empty if $p=0$ and $\mu$ is empty if $q=0$. 

We will call the coordinates of $(\lambda,\mu)\in(\RR^*)^{p+q}$ parameters of the symmetric space $X_{p,q}(\lambda,\mu)$. For a coordinate $\lambda_i$ of $\lambda=(\lambda_1,\dots,\lambda_p)$, the {\it multiplicity} of $\lambda_i$ as a parameter of $X_{p,q}(\lambda,\mu)$ is the number of coordinates of $\lambda$ that are equal to $\lambda_i$ or $-\lambda_i$. In the same way we define the multiplicity of a coordinate of $\mu$.
\end{de}

\subsection{The isometry group}\label{S22}
Next we are going to determine the isometry group $\Iso(X)$ of a Cahen-Wallach space $X$. Suppose that $X$ is associated with the symmetric triple $\osc(\omega,\theta_\fa,L)$.
\begin{pr}\label{Piso}
The isometry group of $X$ is isomorphic to $\hat G\rtimes (K\times\ZZ_2)$, where 
$$K:=\{\ph\in\grO(\fa)\mid \theta_\fa\ph=\ph\theta_\fa,\ L\ph=\ph L \}
$$ acts on $\hat G$ by its standard representation on $\fa\subset \hat G$ and  $-1\in\ZZ_2$ acts on $\hat G$ by~$\theta$.
\end{pr}
\proof 
Let us consider the stabiliser $P:=\{f\in \Iso(X)\mid f(x_0)=x_0\}$ of $x_0=e\hat G_+\in \hat G/\hat G_+$. This group acts on $\hat G$ by conjugation and the homomorphism $\hat G\rtimes P \rightarrow \Iso(X)$, $(g,p)\mapsto gp$ is surjective. The kernel of this map equals $\{(g,p)\in \hat G\rtimes P\mid p=g^{-1}\}=\{(g,g^{-1})\mid g\in \hat G_+\}\cong \hat G_+$, where the latter isomorphism is given by $(g,g^{-1})\mapsto g^{-1}$. Hence $\Iso(X)\cong \hat G_+\setminus(\hat G\rtimes P)$, where $g\in \hat G_+$ acts by $(\hat g,p)\mapsto(\hat gg^{-1}, gp)$. 

Note that $\hat G_+$ is a normal subgroup of $P$. The quotient $\hat G_+\backslash P$ is isomorphic to the subgroup $$P_0:=\{f\in P\mid \pro_{\fa_-}(df_{x_0}(t_0))=0\}$$ 
of $P$, where $t_0:=(0,0,1)\in \fz\oplus\fa\oplus\RR$.
Indeed, for all $a\in\fa_+=\hat G_+$, the differential of the isometry $a:X\rightarrow X$ at $x_0$ equals $\Ad (a)|_{\hat \fg_-}$ on $T_{x_0}X\cong\hat \fg_-$. We calculate $\Ad(a)$ from (\ref{E*Hn}) and (\ref{El}) and obtain
$$\pro_{\fa_-}(\Ad(a)df_{x_0}(t_0))=\pro_{\fa_-}(df_{x_0}(t_0))-\pro_{{\Bbb R}}(df_{x_0}(t_0))\cdot La.$$  Hence, for all $f\in P$, there is exactly one element $\hat a\in \hat G_+$ such that $\hat a\circ f$ is in $P_0$. 

We obtain $\Iso(X)\cong  \hat G_+\setminus (\hat G\rtimes P)\cong \hat G\rtimes P_0$.
It remains to show that $P_0\cong K\times \ZZ_2$ and to determine the induced action of $K\times \ZZ_2$ on $\hat G$ by conjugation.  Conjugation by $f\in P\subset \Iso(X)$ is an automorphism $F$ of $\hat G$. The differential $F_*:=dF_e$ of $F$ at the identity belongs to $\Aut(\hat \fg,\theta,\ip)$ and equals $df_{x_0}$ on $\hat \fg_-=T_{x_0}X$. Hence we obtain a homomorphism
$$P_0\longrightarrow\cA:=\{\Phi \in \Aut(\hat \fg,\theta,\ip) \mid \pro_{\fa_-}(\Phi(t_0))=0\},\quad f\longmapsto F_*.$$ 
This homomorphism is injective since $f\in P_0$ is determined by $df_{x_0}=F_*|_{\hat \fg_-}$. We fix an element $z_0\in\fz$, $z_0\not=0$. Using that $\fz=\RR z_0$ is the centre of $\hat \fg$ and must be preserved by any $\Phi\in\Aut(\fg)$ it is not hard to prove that 
\begin{equation}\label{Prenzlau}
\cA\cong\left\{\Phi\in\Aut(\hat\fg)\ \left|\  \begin{array}{l} \Phi(z_0)=\kappa z_0,\  \Phi(t_0)=\kappa t_0,\ \kappa\in\{1,-1\},\\
\Phi(\fa)\subset\fa,\ \ph:=\Phi |_\fa\in \grO(\fa),\
\theta_\fa\ph=\ph\theta_\fa,\ \ph L=\kappa L\ph
\end{array}\right.\right\}.
\end{equation}
The homomorphism $P_0\longrightarrow\cA$ is also surjective. Indeed, suppose that $\Phi$ is an element of the set on the right hand side of (\ref{Prenzlau}).
Then there is an automorphism $F$ of $\hat G$ such that $\Phi$ is the differential $F_*$ of $F$ at the identity, namely, 
\begin{equation} \label{EF}
F:\ \hat G\longrightarrow \hat G,\quad (z,a,t)\longmapsto (\kappa z, \ph(a), \kappa t).
\end{equation} Since $F_*$ commutes with $\theta$ and preserves $\ip$, $F$ defines an isometry $f$ of $X$ by $f(g\hat G_+)= F(g)\hat G_+$.   We have $df_{x_0}=F_*|_{\hat\fg_-}$, which implies that $f$ is in $P_0$ and that $F_*$ coincides with the automorphism of $\hat \fg$ induced by the conjugation by $f$. 

We proved that $P_0$ is isomorphic to $\cA$. On the other hand, identifying $\cA$ with the right hand side of (\ref{Prenzlau}) we get an isomorphism $$\cA\longrightarrow K\times \ZZ_2,\quad \Phi \longmapsto (\ph\sigma_\kappa,\kappa),$$
 where $\sigma_{-1}=\theta_\fa\in \grO(\fa)$ and $\sigma_1$ is the identity on $\fa$.
Combining these two isomorphisms we get an isomorphism $P_0\cong K\times \ZZ_2$. Take $f\in P_0$ and let $(\ph \sigma_\kappa,\kappa)$ be its image under this isomorphism. By construction, the conjugation by $f$ equals the map $F:\hat G\rightarrow \hat G$ given by (\ref{EF}). This proves the assertion on the action of $K\times \ZZ_2$ on $\hat G$.
\qed

Choose $p,q,\lambda,\mu$ such that $X\cong X_{p,q}(\lambda;\mu)$. If $\dim K>0$, then some parameter $\lambda_i$ or $\mu_j$ has multiplicity greater than one. The subgroup
\begin{equation}\label{defG} G:= \hat G\rtimes K\subset \Iso(X)
\end{equation}
acts transitively on $X$. The subgroup $K\subset G$ normalises $\hat G_+$ since all elements of $K$ commute with $\theta_\fa$. Hence, the stabiliser of $e\hat G_+\in X=\hat G/\hat G_+$ equals $G_+:=\hat G_+\rtimes K$. Thus we obtain $X=G/G_+$.

\subsection{Lie groups with a biinvariant Lorentzian metric}\label{groupcase}
An interesting subclass of Cahen-Wallach spaces is constituted by solvable Lie groups endowed with a biinvariant Lorentzian metric. The infinitesimal object that is associated with such a group is a solvable metric Lie algebra of index one, i.e., a Lie algebra endowed with an $\ad$-invariant non-degenerate scalar product of signature $(1,n+1)$. Take, for example, the symmetric triple $\osc_{0,m}(\mu)$ and forget about the involution. Then you get a solvable metric Lie algebra of signature $(1, 2m+1)$, which we will also denote by $\osc_{0,m}(\mu)$. Medina \cite{M} proved that each indecomposable solvable metric Lie algebra of signature $(1,n+1)$ is isomorphic to $\osc_{0,m}(\mu)$ for exactly one $\mu\in\RR^m$ with $1=\mu_1\le\mu_2\le\dots\le\mu_m$, where $n=2m$. 
\begin{pr} A Cahen-Wallach space $Q$ is a Lie group endowed with a biinvariant Lorentzian metric if and only if it is isometric to  some $X_{0,2m}(\tilde\mu)$, where $\tilde \mu=(\mu_1,\mu_1,\mu_2,\mu_2,\dots, \mu_m,\mu_m)\in\RR^{2m}$.
\end{pr}
\proof
Let $Q$ be a Lie group with biinvariant Lorentzian metric. We consider $Q$ as a symmetric space
and we wish to determine the associated symmetric triple. Note first that the action of $Q\times Q$ on $Q$ defined by $(q_1,q_2)\cdot q=q_1qq_2^{-1}$  is isometric since the metric on $Q$ is biinvariant. The kernel of this action is isomorphic to the centre $Z(Q)$ of $Q$. More exactly, it equals $\{(z,z)\in Q\times Q\mid z\in Z(Q)\}$. Hence $I:=(Q\times Q)/Z(Q)$ is a subgroup of the isometry group of $Q$. This subgroup contains the transvection group of $Q$ since the reflection at a point $q\in Q$ is given by $Q\ni p\mapsto qp^{-1}q\in Q$.  Moreover, it is invariant under the conjugation by the reflection of $Q$ at the identity, which we denote by $\theta$. If $\fq$ denotes the Lie algebra of $Q$, then the Lie algebra of $I$ equals $(\fq\oplus\fq)/\fz(\fq)$. Now we consider the eigenspace decomposition of $(\fq\oplus\fq)/\fz(\fq)$ with respect to the differential of $\theta$. The (-1)-eigenspace $\hat \fg_-$ equals the anti-diagonal $\{(X,-X)\mid X\in\fq\}\cong \fq$. The subspace $\hat \fg_+:=[ \hat \fg_-,\hat \fg_-]$ of the (+1)-eigenspace equals
$\{(X,X)\mid X\in[\fq,\fq]\}/(\fz(\fq)\cap [\fq,\fq]).$
 Hence the Lie algebra $\hat\fg$ of the transvection group of $Q$ is isomorphic to 
$$
\hat\fg=\hat\fg_+\oplus\hat \fg_-=([\fq,\fq]/(\fz(\fq)\cap[\fq,\fq]))\oplus \fq
$$
with Lie bracket 
\begin{equation}\label{E7} 
[(X_1,Y_1), (X_2,Y_2)]=([X_1,X_2]+[Y_1,Y_2],\ [X_1,Y_2] + [Y_1,X_2]).
\end{equation}
By the discussion above, the metric Lie algebra $\fq$ is isomorphic to $\osc_{0,m}(\mu)$ for some $\mu=(\mu_1,\dots,\mu_m)\in(\RR^*)^m$. In the notation of Example \ref{Exosc2} we have $[\fq,\fq]/(\fz(\fq)\cap[\fq,\fq])\cong\fa$, hence 
$$\hat \fg=\hat\fg _+\oplus \hat\fg_-=\fa \oplus (\fz\oplus\fa\oplus\RR).$$
We denote by $t_0$ the element $1\in\RR$ in the last summand of this direct sum. Since we already know that $\hat\fg$ is isomorphic to some $\osc_{\tilde p,\tilde q}(\tilde \lambda,\tilde \mu)$ it suffices to determine the adjoint action of $t_0$ on $[\hat \fg,\hat \fg]/\fz(\hat\fg)=\fa\oplus\fa$. More exactly, it suffices to determine the eigenvalues of $\ad(t_0)^2$ on $(\fa\oplus\fa)\cap \hat\fg_-=\fa$.  By (\ref{E7}) these are exactly the eigenvalues of $L^2$ on $\fa$, which are $-\mu_1^2,\dots,-\mu_m^2$, each with multiplicity two. Consequently, the symmetric triple associated with $Q$ is isomorphic to $\osc_{0,2m}(\tilde\mu)$, where $\tilde \mu=(\mu_1,\mu_1,\mu_2,\mu_2,\dots, \mu_m,\mu_m)\in\RR^{2m}$. These considerations also show that each $X_{0,2m}(\tilde\mu)$ is a group with a biinvariant Lorentzian metric.
\qed

If $\Gamma\subset Q:=X_{0,2m}(\tilde\mu)$ is a lattice, then $\Gamma \setminus Q$ is a compact quotient of the symmetric space~$Q$. The investigation of lattices in oscillator groups was started in \cite{MR2}. However, the results in \cite{MR2} are not correct.  It turns out that the structure of a general lattice is more complicated than claimed in that paper. A description of these lattices including a complete classification for four-dimensional oscillator groups, i.e., for $m=1$, can be found in~\cite{F}. 
In Section \ref{klumpfuss}, we will see that quotients by lattices of $Q$ only give very special examples of compact quotients of the symmetric space $Q$. In particular, we will see that every symmetric space $Q=X_{0,2m}(\tilde\mu)$ admits compact quotients that do not come from a lattice in the Lie group $Q$. 

\subsection{The canonical fibration}\label{canfib}
Recall from Subsection \ref{S22} that $X=G/G_+$ and that $G_+\subset H\rtimes K$, where $H$ denotes the Heisenberg group. Hence there is a natural projection
$$\pi:\ X=G/G_+\longrightarrow G/(H\rtimes K)\cong \RR.$$
This projection defines a locally trivial fibration. The fibres are flat, coisotropic and connected. The radical of the restriction of the Lorentzian metric to the fibres is a one-dimensional subbundle of the tangent bundle. Hence it defines a foliation with one-dimensional leaves called null-leaves. 

Since also the action of $\ZZ_2\subset\Iso(X)$ on $G$ leaves invariant both subgroups $G_+$ and $H\rtimes K$  of $G$, the fibration is equivariant with respect to the action of $\Iso(X)$.


\section{Discrete subgroups of the isometry group}\label{S3}

We are interested in subgroups $\Gamma\subset \Iso(X)$ of the group of isometries of a Cahen-Wallach space $X$ that act properly discontinuously and cocompactly on $X$. 
As a first step we prove a structural result for arbitrary discrete subgroups of $\Gamma\subset \Iso(X)$. Using this result properness and cocompactness of the action of $\Gamma$
will be studied in Section~\ref{klumpfuss}.

In fact, we will investigate discrete subgroups of slightly more general Lie groups. Namely, we will study discrete subgroups $\Gamma$ of arbitrary Lie groups $G$ of the form
\begin{equation}\label{setting}
G=N\rtimes_\rho (\RR\times K)\,,
\end{equation}
where 
\begin{itemize}
\item $N$ is 1-connected nilpotent,
\item $K$ is compact, and
\item $\rho:\RR\times K\rightarrow \Aut(N)\cong \Aut(\fn)$ is an action by semisimple automorphisms.
\end{itemize}
The results will be applicable to the isometry group since $\Iso(X)$ admits a subgroup $G$ of index 2 that has exactly this form, see Prop.~\ref{Piso} and Equ.~(\ref{defG}). 

 Let $r:G\rightarrow \RR\times K$ and $p:G\rightarrow \RR$ be the natural projections.
\begin{de}\label{unfug}
Let $G$ be as in {\rm (\ref{setting})}. A discrete subgroup $\Gamma\subset G$ is called tame, if the closure $\overline{\rho(r(\Gamma))}$ in $\Aut(N)$ has only finitely many connected components. 
\end{de}
\begin{lm} \label{Ltame}
The subgroup $\Gamma\subset G$ is tame if and only if $\overline{\rho(\RR)}$ is compact or $\overline{p(\Gamma)}$ is connected.
\end{lm}
\proof
If $\overline{\rho(\RR)}$ is compact, then also $\overline{\rho(\RR\times K)}$ is compact, hence $\overline{\rho(r(\Gamma))}\subset\overline{\rho(\RR\times K)}$ is compact. In particular, $\overline{\rho(r(\Gamma))}$ has only finitely many connected components, thus $\Gamma$ is tame.

Now let $\overline{p(\Gamma)}$ be connected. The restriction of the projection $\RR\times K\rightarrow \RR$ to $\overline{r(\Gamma)}$  gives a Lie group homomorphism $\overline{r(\Gamma)}\rightarrow   \overline{p(\Gamma)}$. This morphism is surjective since $K$ is compact. Hence it defines a fibre bundle $C\rightarrow \overline{r(\Gamma)}\rightarrow   \overline{p(\Gamma)}$, where the fibre $C$ is contained in $\{0\}\times K$  and thus is compact. Now we see from the long exact homotopy sequence of this fibre bundle that $\overline{r(\Gamma)}$ has only finitely many connected components and we conclude that the same is true for $\overline{\rho(r(\Gamma))}$, hence $\Gamma$ is tame.

Now suppose that $\overline{\rho(\RR)}$ is not compact and that $\overline{p(\Gamma)}$ is not connected. In this case, we have $p(\Gamma)=\ZZ\cdot t_0$. All elements of the one-parameter subgroup $\rho(\RR)$ of $\Aut(N)\cong \Aut(\fn)$ act semisimply on $\fn$ and commute with $\rho(K)$. Hence there is a $K$-invariant decomposition of the complexification $\fn_{\Bbb C}$ of $\fn$ into common eigenspaces for all elements of $\overline{\rho(\RR)}$. We choose a $K$-invariant norm $\|\cdot\|$ on $\fn_{\Bbb C}$.  We find a common eigenvector $v$ of $\overline{\rho(\RR)}$ such that $\rho(t_0)(v)=\lambda v$, where $|\lambda|\not=1$ since $\overline{\rho(\RR)}$ is not compact. We may assume $\|v\|=1$. Now we consider the map $\overline{\rho(r(\Gamma))}\rightarrow \RR^*$, $a\mapsto \|\rho(a)(v)\|$. This map is continuous and its image equals $\{ |\lambda|^k\mid k\in\ZZ\}$. Hence $\overline{\rho(r(\Gamma))}$ has infinitly many connected components.
\qed

The following two propositions generalise a classical result of Auslander \cite{Aus} on discrete subgroups of semidirect products of nilpotent with compact Lie groups.
They say that $\Gamma$ is, essentially, a lattice in a certain closed subgroup of $G$, which is either connected or has an infinite cyclic component group.

\begin{pr}\label{nottame}
Let $G$ be as in {\rm (\ref{setting})}, and let $\Gamma\subset G$ be a discrete subgroup. Then there exist
\begin{enumerate}
\item[(a)] an element $n\in N$ and a subgroup of finite index $\Gamma_0\subset n\Gamma n^{-1}$,
\item[(b)] a closed abelian subgroup $C\subset \RR\times K$ with $C_K:=C\cap K$ connected,
\item[(c)] a connected $C$-invariant subgroup $U\subset N^{C_K}$, and
\item[(d)] a group homomorphism $\psi: \Delta:=\overline{p(\Gamma_0)}\rightarrow N_N(U)^C\times C$ with $p\circ\psi=\id_\Delta$ and $C=r(\psi(\Delta))\times C_K$
\end{enumerate}
such that $\Gamma_0\subset (U\cdot\psi(\Delta))\times C_K$ is cocompact.

If $\Delta\not=\RR$, then $\psi$ may be chosen such that $\psi(\Delta)\subset\Gamma_0$.
\end{pr}

\begin{pr}\label{tame}
Let $G$ be as in {\rm(\ref{setting})}, and let $\Gamma\subset G$ be a tame discrete subgroup. Then the conclusion of Proposition \ref{nottame} holds with some
$U\subset N^C$.
\end{pr}

Note that $\Delta$ is either trivial, infinite cyclic, or equal to $\RR$. If $\Gamma$ is not tame, then $\Delta$ is infinite cyclic.

The proof of the propositions will occupy the remainder of the section. We follow quite closely the (very sketchy) arguments in \cite{Aus}.
Let us remark that it would be also possible to base the proof on Witte's result \cite{WM1}, \cite{WM2} on existence of syndetic hulls in solvable Lie groups.

We start with a couple of certainly well-known lemmas of preparatory character.

\begin{lm}\label{eins}
Let $N$ be a 1-connected nilpotent Lie group, $n\in N$, and let $a\in\Aut(N)$ be a semisimple element. Then there exists $n_1\in N$ such that
\begin{equation}\label{putin} n_1n\, a(n_1)^{-1}\in N^a\ .
\end{equation}
In other words: The element $(n,a)\in N\rtimes\langle a\rangle$ is conjugate via $n_1$ to $(n',a)$ with $n'\in N^a$.
\end{lm}

\proof 
We prove (\ref{putin}) by induction on the nildegree of $N$. First, let $N$ be abelian.
Then $N$ is the additive group of a vector space, $a$ is a semisimple linear map. This yields
the decompositions
$$ N=\Ker(\id-a)\oplus\im(\id-a)=N^a\oplus\im(\id-a),\quad n=n_0+(a-\id)n_1 $$
for some $n_0\in N^a$, $n_1\in N$. Then $n_1n\, a(n_1)^{-1}$ is equal to (written additively)
$$ n_1+n-a(n_1)=n-(a-\id)n_1=n_0\in N^a\ .$$
Now let $N$ be arbitrary. We assume that (\ref{putin}) holds for $\bar N:=N/Z(N)$.
This implies the existence of $n_2\in N$ with $\bar n_2\bar n a(\bar n_2)^{-1}\in \bar N^a$.
Here $\bar n:=nZ(N)\in\bar N$. Since $a$ is semisimple the natural map $N^a\rightarrow \bar N^a$ is surjective. Therefore we have 
\begin{equation}\label{kohl} n_2n\, a(n_2)^{-1}z^{-1}\in N^a
\end{equation}
for some $z\in Z(N)$. Since $Z(N)$ is abelian we find $z_1\in Z(N)$ with
\begin{equation}\label{rasputin} zz_1\, a(z_1)^{-1}=z_1z\, a(z_1)^{-1}\in N^a\ .
\end{equation}
Multiplying (\ref{kohl}) with (\ref{rasputin}) we obtain
$ N^a\ni n_2n\, a(n_2)^{-1}z_1\, a(z_1)^{-1}=(n_2z_1)n\, a(n_2z_1)^{-1}\ .$
Thus $n_1:=n_2z_1$ does the job.
\qed

\begin{lm}\label{zwei}
Let $N$ be a 1-connected nilpotent Lie group, and let $a\in\Aut(N)$ be a semisimple element. Let $n_1\in N^a$, $n_2\in N$ be such that
$$ n_1\, a(n_2)n_1^{-1}n_2^{-1}\in N^a\ .$$
Then $n_2\in N^a$.
\end{lm}

\proof
We consider the descending central series  $N^1=N$, $N^{k+1}=[N,N^{k}]$,  $k\ge 1$. Suppose $n_2\in N^k$. The assertion is obvious if $k$ is sufficiently large. Let us assume that the assertion is true for $k+1$ and prove it for $k$.  As in the proof of Lemma~\ref{eins} we consider the abelian group $\bar N:=N^k/N^{k+1}$ as the additive group of a vector space. By assumption, $(a-\id)\bar n_2\in \bar N^a$. This implies $\bar n_2\in \bar N^a$. Using the surjectivity of the map $(N^k)^a\rightarrow \bar N^a$ we find elements $n_2^0\in (N^k)^a$ and
$n_2'\in N^{k+1}$ such that $n_2=n_2^0n_2'$. We obtain 
$$ N^a\ni n_1\, a(n_2^0n_2')n_1^{-1}(n_2^0n_2')^{-1}
=n_1n_2^0\, a(n_2')n_1^{-1}n_2'^{-1}(n_2^0)^{-1}\ ,$$
hence $a(n'_2)n_1^{-1}(n'_2)^{-1}\in N^a$, which gives $n_1\, a(n_2')n_1^{-1}n_2'^{-1}\in N^a$. By induction hypothesis, $n'_2\in N^{a}$, hence $n_2=n_2^0n_2'\in N^a$.
\qed

\begin{lm}\label{zusatz}
Let $N$ be a 1-connected nilpotent Lie group, let $C\subset\Aut(N)$ be a subgroup, and let $U\subset N$ be a connected subgroup. Then $N^C$ and $N_N(U)$ are connected, hence 1-connected.
\end{lm}

\proof
The exponential map $\exp: \fn\rightarrow N$ is an $\Aut(N)$-equivariant diffeomorphism.
Hence $N^C=\exp(\fn^C)$, and $\fn^C$ is a vector space. 
As for $N_N(U)$, assume that $\exp(X)\in N_N(U)$ for some $X\in\fn$. We claim that $\exp(tX)\in N_N(U)$
for all $t\in \RR$. Indeed, let $Y$ be an element of the Lie algebra $\fu$ of $U$ and let $\ph\in\fn^*$ be a functional vanishing on $\fu$. Then the polynomial $t\mapsto \ph(\Ad(\exp(tX))Y)$ vanishes at integral $t$, hence it vanishes identically. Since $\ph$ was arbitrary, $\Ad(\exp(tX))Y\in\fu$. 
\qed

\begin{lm}\label{drei}
For all $d\in\NN$, there exists a neighbourhood $U(d)\subset\CC$ of $1$ such that all
$d\times d$-matrices with integer entries and all eigenvalues in $U(d)$ are unipotent.
\end{lm}

\proof
We consider the map $s: \CC^d\rightarrow\CC^d$ given by the elementary symmetric polynomials $s(\lambda):=(s_1(\lambda), s_2(\lambda),\dots,s_d(\lambda))$. If $\lambda=(\lambda_1,\lambda_2,\dots,\lambda_d)$ is the $d$-tuple of eigenvalues of an integral matrix, then $s(\lambda)\in\ZZ^d$. Choose a neighbourhood $U(d)$ of $1\in\CC$ such that
$s(U(d)\times U(d)\times\dots\times U(d))\cap \ZZ^d= \{s(1,1,\dots,1)\}$.
\qed

We will also use repeatedly the following classical result, see e.g. \cite{Ragh}, Prop.~2.5.

\begin{lm}\label{drei.neun} 
Let $N$ be 1-connected nilpotent, and let $\Gamma\subset N$ be a closed subgroup. Then there is a unique connected subgroup $U\subset N$ such that $\Gamma\backslash U$ is compact. If $\Gamma$ is abelian, then so is $U$. \qedohne
\end{lm}
{\sl Proof of Proposition \ref{tame}. } Let $\Gamma\subset G$ be a tame discrete subgroup.
We first construct $C\subset \RR\times K$. We distinguish between 2 cases.

{\sl Case 1: }$\overline{p(\Gamma)}\subset\RR$  connected.  As already observed in the proof of Lemma \ref{Ltame}, $\overline{r(\Gamma)}\subset K\times\RR$ has only finitely many connected components in this case.
We set $C:=\overline{r(\Gamma)}_0$. The subgroup $\tilde \Gamma_0:=r^{-1}(C)\cap \Gamma$ has finite index in $\Gamma$. By a theorem of Auslander (\cite{Ragh}, Thm.~8.24)
the group $C\subset \RR\times K$ is solvable, hence abelian.

{\sl Case 2: }${p(\Gamma)}=\langle p(\gamma_0)\rangle$ for some $\gamma_0\in\Gamma$. We set $C_K:=\overline{r(\Gamma\cap\Ker p)}_0\subset K$. As in Case~1, Auslander's theorem implies that $C_K$ is abelian. Conjugation by $r(\gamma_0)$ induces an automorphism of the torus $C_K$. This automorphism is of finite order since only the $K$-component of $r(\gamma_0)$ matters. Thus there is some $k>0$ such that $r(\gamma_0^k)$
commutes with $C_K$. The group $\tilde C:=\langle C_K, r(\gamma_0^k)\rangle$ is abelian and closed. Since $\Gamma$ is tame and $\tilde C\subset\overline{r(\Gamma)}$ has finite index, the closure $\overline{\rho(\tilde C)}\subset \Aut(N)$ has only finitely many connected components. We set $C:=\{c\in\tilde C\mid \rho(c)\in \overline{\rho(\tilde C)}_0\}$. Then $\tilde \Gamma_0:=r^{-1}(C)\cap \Gamma$ has finite index in $\Gamma$.

In both cases we have: $\overline{\rho(C)}\subset \Aut(N)$ is connected, abelian, and
${\rho(r(\tilde\Gamma_0))}\subset\overline{\rho(C)}$ is dense. This ensures that for every neighbourhood $V_1$ of $\id\in \overline{\rho(C)}$ there exists an element $\gamma_1\in \tilde\Gamma_0$ such that $\rho(r(\gamma_1))\in V_1$ and $N^C=N^{r(\gamma_1)}$. Indeed, $\overline{\rho(C)}$ acts on $\fn/\fn^C$. This representation of $\overline{\rho(C)}$ decomposes into non-trivial irreducible subrepresentations, which have kernels $D_j\subset \overline{\rho(C)}$, $j=1,\dots, J$, of codimension one.
Now we choose $\gamma_1$ such that $\rho(r(\gamma_1))\in V_1\setminus \bigcup_{j}
D_j$. 

We consider the discrete subgroup $\tilde \Gamma_0\cap N$ of $N$. It is a lattice in a connected subgroup $\tilde U_0\subset N$. Let $d=\dim\tilde U_0\in\NN_0$. If $d>0$ we choose $V_1$ sufficiently small such that all eigenvalues of elements of $V_1$ on $\fn$ belong to the neighbourhood $U(d)$ provided by Lemma \ref{drei}. Let us fix a corresponding element $\gamma_1\in \tilde\Gamma_0$ as in the previous paragraph.

We write $\gamma_1=(n_1,c_1)\in N\rtimes C$. By Lemma \ref{eins} there exists an element
$n\in N$ such that $\gamma':=n\gamma_1n^{-1}=(n',c_1)$ for some $n'\in N^{c_1}=N^C$.
We set $\Gamma_0:=n\tilde \Gamma_0 n^{-1}$, $U_0:=n\tilde U_0 n^{-1}$.
Then $\Gamma_0$ has finite index in $n\Gamma n^{-1}$, and $\Gamma_0\cap N$ is a lattice in $U_0$.

\begin{lm}\label{vier} We have $U_0\subset N^C$.
\end{lm}

\proof 
The element $\gamma'$ normalises $\Gamma_0\cap N$ and therefore also $U_0$.
This implies that the linear map $\Ad(\gamma')$ on $\fn$ leaves invariant the Lie algebra $\fu_0$
of $U_0$ as well as a lattice in $\fu_0$. Thus, in a suitable basis,  $\Ad(\gamma')|_{\fu_0}$
is given by a matrix in $\GL(d,\ZZ)$.

On the other hand, $\gamma'=(n',c_1)$, $\Ad(n')$ is unipotent, $\Ad(c_1)$ is semisimple, and these maps commute. Hence
$\Ad(\gamma')=\Ad(n')\Ad(c_1)$ is the multiplicative Jordan decomposition of $\Ad(\gamma')$.
This implies that $\fu_0$ is invariant under $\Ad(c_1)$ and that the eigenvalues of
$\Ad(\gamma')|_{\fu_0}$ and $\Ad(c_1)|_{\fu_0}=\rho(r(\gamma_1))|_{\fu_0}$ coincide.
Since $\rho(r(\gamma_1))\in V_1$ we can apply Lemma \ref{drei} to $\Ad(\gamma')|_{\fu_0}$ and conclude that all these eigenvalues are equal to $1$.
It follows that $\Ad(c_1)|_{\fu_0}=\id_{\fu_0}$, thus $U_0\subset N^{c_1}=N^C$.
\qed

\begin{lm}\label{fuenf} If $\gamma=(n,c)\in\Gamma_0$, then $n\in N^C$.
\end{lm}

\proof
Let $\gamma'=(n',c_1)$ be as above. Then $[\gamma',\gamma]\in \Gamma_0\cap N\subset U_0\subset N^C$. The last inclusion comes from Lemma \ref{vier}. Using that $C$ is abelian and that $n'\in N^{c_1}=N^C$ we compute
$$
[\gamma',\gamma]=n'c_1ncc_1^{-1}{n'}^{\, -1}c^{-1}n^{-1} 
= n'\rho(c_1)(n){n'}^{\, -1} n^{-1}\ .$$
Thus $n'\rho(c_1)(n){n'}^{\, -1} n^{-1}\in N^{c_1}$. Now Lemma \ref{zwei} implies that
$n\in N^{c_1}=N^C$.
\qed

Now we consider the normal subgroup $\Gamma_1:=\Ker p\,\cap\Gamma_0\subset\Gamma_0$.
It is contained in $N^C\times C_K$. Let $q$ be the projection on the first component.
Since $C_K$ is compact, the subgroup $q(\Gamma_1)\subset N^C$ is discrete.
$N^C$ is 1-connected nilpotent (see Lemma \ref{zusatz}). Thus $q(\Gamma_1)$ is a lattice in a certain connected subgroup $U_1\subset N^C$, $\Gamma_1$ is a lattice in $U_1\times C_K$.

\begin{lm}\label{sex} If $\gamma=(n,c)\in\Gamma_0$, then $n\in N_N(U_1)$.
\end{lm}

\proof Let $\gamma=(n,c)\in\Gamma_0$. It suffices to show that $n$ normalises $q(\Gamma_1)$. Let $(n_2,c_2)\in\Gamma_1$. Using Lemma \ref{fuenf} we compute
$$\Gamma_1\ni \gamma(n_2,c_2)\gamma^{-1}=(n\rho(c)(n_2)n^{-1}, c_2)
=(nn_2n^{-1}, c_2)\ .
$$ 
We conclude that $nn_2n^{-1}\in q(\Gamma_1)$ as desired.
\qed

Eventually we construct $U\subset N^C$ and $\psi: \Delta\rightarrow N_N(U)^C\times C$.
If $\Delta=\{0\}$, then we set $U:=U_1$, and we are done. If $\Delta\cong\ZZ$, then we also set  
$U:=U_1$ and choose $\gamma_2\in \Gamma_0$ such that $\delta:=p(\gamma_2)$ generates $\Delta$. We set $\psi(\delta^k):=\gamma_2^k$ for $k\in\ZZ$.  According to Lemma \ref{fuenf} and Lemma \ref{sex} we have $\psi(\Delta)\subset N_N(U)^C\times C$. Since $\Gamma_1$ is a lattice in $U\times C_K$ the group $\Gamma_0$ is a lattice in $(U\cdot \psi(\Delta))\times C_K$.

It remains to discuss the case $\Delta=\RR$.  By 
Lemma \ref{fuenf} and Lemma \ref{sex} we have $\Gamma_0\subset N_N(U_1)^C$. We consider
the abelian group
$$ \bar\Gamma_0:= \Gamma_0/\Gamma_1\subset (N_N(U_1)^C\times C)/(U_1\times C_K)
\cong V\times \Delta\,,
$$
where $V:=N_N(U_1)^C/U_1$. 
The image of a discrete group $\Gamma\subset G$ in $G/N$, $N\subset G$ being a normal subgroup of a Lie group $G$, is discrete provided
that $\Gamma\cap N$ is a lattice in $N$. 
Since $\Gamma_1$ is a lattice in $U_1\times C_K$, we conclude that $\bar\Gamma_0\subset V\times \Delta$ is discrete.
By Lemma \ref{zusatz} the group $V\times\Delta$ is 1-connected nilpotent. It follows that there is an abelian connected
subgroup $W\subset V\times \Delta$ such that  $\bar\Gamma_0\subset W$ is a lattice. The group $W$ projects surjectively to $\Delta$. Let $\psi_0: \Delta\rightarrow W$ be a lift along this projection. Let $\pi: N_N(U_1)^C\times C\rightarrow V\times \Delta$ be the natural projection. We lift $\psi_0$ further along $\pi$ and obtain a homomorphism $\psi: \Delta\rightarrow N_N(U_1)^C\times C$. Set $U:=\pi^{-1}(W)\cap N_N(U_1)^C$. Since $W$ is abelian, the group $\pi(\psi(\Delta))=\psi_0(\Delta)\subset W$ normalises $W\cap V$. Hence $\psi(\Delta)$ normalises $\pi^{-1}(W\cap V)=U\times C_K$.  Thus, $\psi(\Delta)\subset N_N(U)^C\times C$. It remains to prove that $\Gamma_0$ is a lattice in $(U\cdot\psi(\Delta))\times C_K$.
Since, by construction, $\Gamma_0$ is a lattice in $\pi^{-1}(W)$, it suffices to show that $\pi^{-1}(W)=(U\cdot\psi(\Delta))\times C_K$. Obviously, $\pi^{-1}(W)\supset(U\cdot\psi(\Delta))\times C_K$. To show equality, take $x\in\pi^{-1}(W)$. Since $W=(W\cap V)\cdot \psi_0(\Delta)$, we get $\pi(x)\psi_0(t)^{-1}\in W\cap V$ for some $t\in\Delta$, thus  $x\psi(t)^{-1}\in \pi^{-1}(W\cap V)=U\times C_K$.
\qed

{\sl Proof of Proposition \ref{nottame}. } 
Since Proposition \ref{tame} is already proved it
suffices to discuss discrete subgroups $\Gamma\subset G$ that are not tame. Then
$\Gamma=\langle \Gamma',\gamma_0\rangle$, where $\Gamma':=\Gamma\cap\Ker p$ is tame, and
$\gamma_0\in\Gamma$ is an element such that $p(\gamma_0)$ generates $p(\Gamma)\subset\RR$. Proposition \ref{tame} applied to $\Gamma'$ yields
\begin{enumerate}
\item[(a)] an element $n'\in N$ and a subgroup of finite index $\Gamma_0'\subset n'\Gamma' {n'}^{\, -1}$,
\item[(b)] a torus $C'=:C_K\subset K$,
\item[(c)] a connected subgroup $U'\subset N^{C_K}$
\end{enumerate}
such that $\Gamma_0'\subset U'\times C_K$ is a lattice. As in the proof of Proposition \ref{tame} (discussion of {\sl Case 2}) we see that $C_K$ commutes with $r(\gamma_0)^k$  for some $k\in\NN$. We set $\gamma_2=(n_2,c_2):=n'\gamma_0^k {n'}^{\, -1}$,
$C:=\langle c_2\rangle\times C_K$. 

Arguing as in the proof of Lemma \ref{fuenf} (with a suitable element $\gamma'\in\Gamma_0'\subset N^{C_K}\times C_K$) we see that $n_2\in N^{C_K}$. We now apply Lemma \ref{eins} to $N^{C_K}\rtimes\langle c_2\rangle$ and find
an element $n''\in N^{C_K}$ such that $\gamma_2'=(n_2',c_2):=n''\gamma_2 {n''}^{\, -1}$
satisfies $n_2'\in N^{c_2}\cap N^{C_K}=N^C$.

We now define
$\Gamma_0:=\langle n''\Gamma_0'{n''}^{\, -1},\gamma_2'\rangle$, $U:=n''U'{n''}^{\, -1}$
and $\psi: \Delta\rightarrow N^C\times C$ by $\psi(p(\gamma_2')^k):={\gamma_2'}^k$.
Then $\Gamma_0$ has finite index in $n\Gamma n^{-1}$ ($n=n''n'$) and is a lattice in $(U\cdot\psi(\Delta))\times C_K$. Moreover, $U\subset N^{C_K}$.

The element $\gamma_2'$ normalises $U$, and the Jordan decomposition of $\Ad(\gamma_2')$ on $\fn$ is given by $\Ad(n_2')\rho(c_2)$ (compare the proof of Lemma \ref{vier}). We conclude that both components $n_2'$ and $c_2$ normalise $U$. Therefore
$U\subset N^{C_K}$ is $C$-invariant and $\psi(\Delta)\subset N_N(U)^C\times C$. This finishes the proof of the proposition.
\qed

A discrete group is called virtually nilpotent if it contains a nilpotent subgroup of finite index.
\begin{co}\label{nuff}
Let $\Gamma\subset G$ be a tame discrete subgroup. Then $\Gamma$ is virtually nilpotent.
\end{co}
\proof
Let $\Gamma_0\subset (U\cdot\psi(\Delta))\times C_K$ be as in Proposition \ref{nottame}. By Proposition \ref{tame} we have $U\subset N^C$. This implies that $(U\cdot\psi(\Delta))\times C_K$ is nilpotent. We conclude that $\Gamma_0$ is nilpotent, too.
\qed
 
     
\section{Proper and cocompact actions on Cahen-Wallach spaces} \label{klumpfuss}

In this section we will derive criteria for a Cahen-Wallach space $X$ to admit a compact quotient, see Proposition~\ref{wulf} and Theorem~\ref{AA} below. For this we have to decide whether $\Iso(X)$ contains a subgroup $\Gamma$ acting properly discontinuously, cocompactly and freely on $X$. 
By replacing $\Gamma$ by a subgroup of index 2, if necessary, we may assume that $\Gamma\subset G$, where
$$ G= H\rtimes (\RR\times K) \subset \Iso(X) $$
is as defined in (\ref{defG}).
The freeness condition will turn out to be harmless.
Thus we are looking for discrete $\Gamma\subset G$ acting properly and cocompactly on $X$. The notions of properness and cocompactness make sense for non-discrete groups, too. 
Recall that the action of a locally compact Hausdorff topological group $T$ on a locally compact Hausdorff space $M$ is called proper if the map $T\times M\rightarrow M\times M$, $(t,m)\mapsto (m,tm)$, is proper. 
In other words, for every compact $D\subset M$ the subset $T_D:=\{t\in T\mid tD\cap D\ne\emptyset\}\subset T$ is compact. Orbit
spaces $T\backslash M$ of proper actions are Hausdorff with respect to the quotient topology. A proper action is called cocompact, if $T\backslash M$ is compact. Let $T_1\subset T$ be a cocompact subgroup. Then the action $T_1\times M\rightarrow M$ is proper (and cocompact) if and only if $T\times M\rightarrow M$ is proper (and cocompact).

\begin{de}\label{vier.eins}
Let $\Gamma\subset G$ be a discrete subgroup. We consider the subgroups 
$\psi(\Delta)$ and $U$ of $G$ provided by Proposition \ref{nottame}. 
We define $S_\Gamma:=U\cdot\psi(\Delta)$.
\end{de}
Note that $S_\Gamma$ is only determined up to certain conjugations and, if it is not connected, up to replacements by subgroups (or overgroups) of finite index. But this is sufficient for our purposes. 

\begin{lm}\label{null}
Let $\Gamma\subset G$ be a discrete subgroup. Then $\Gamma$  acts properly discontinuously and cocompactly on $X$ if and only if $S_\Gamma$
acts properly and cocompactly on $X$. 
\end{lm}

\proof Recall from Prop.~\ref{nottame} that there is a finite index  subgroup $\Gamma_0$ of a conjugate of $\Gamma$ that is cocompact in $S_\Gamma\times C_K$ for some closed subgroup $C_K\subset K$. 
The group $\Gamma$  acts properly discontinuously and cocompactly on $X$ iff $\Gamma_0$ does so. Since the inclusions $\Gamma_0\subset S_\Gamma\times C_K\supset S_\Gamma$ are cocompact, the lemma follows. 
\qed

We now study the actions of groups like $S_\Gamma$ on $X$, regardless whether they come from a discrete
subgroup or not.
First we look at the action on the fibres of the canonical fibration $\pi: X\rightarrow \RR$ (see Subsection~\ref{canfib}).

\begin{lm}\label{ulf}
Let $U$ be a connected subgroup of the Heisenberg group $H$ and $\fu$ be its Lie algebra.
Let $\fb\subset (\fa,\omega)$ be a Lagrange subspace, and let $B\subset H$ be the connected subgroup with 
Lie algebra $\fb$. Then the
following conditions are equivalent:
\begin{enumerate} 
\item[(i)] $U$ acts properly and cocompactly on $H/B$;
\item[(ii)]$U$ acts simply transitively on $H/B$;
\item[(iii)] $\fh=\fu\oplus\fb$;
\item[(iv)] $\fz\subset\fu$ and $\fa=(\fa\cap\fu)\oplus\fb$.
\end{enumerate}
\end{lm}

\proof
The implications $(iv)\Rightarrow (iii)$ and $(ii)\Rightarrow (i)$ are obvious.
We show $(iii)\Rightarrow (iv)$, $(iv)\Rightarrow (ii)$ and $(i)\Rightarrow (iv)$. 

Assume $(iii)$ and that $\fz\cap\fu=\{0\}$. Let $0\ne Z\in\fz$. Then $Z=X+Y$ for some
$X\in\fu$, $0\ne Y\in\fb$. Moreover, since $\fb$ is Lagrange, there is some $X_1\in\fu$ with $Z=[Y,X_1]$. We obtain
$$ Z=[Y-Z,X_1]=[X_1,X]\in\fu\ .$$
This is a contradiction. Condition $(iv)$ follows. 

Assume $(iv)$. Then $B$ acts simply transitively on the affine space $\fa/(\fa\cap\fu)\cong H/U$. Hence $U$ acts simply transitively on $H/B$. 

It remains to show the implication $(i)\Rightarrow (iv)$. Let $(Z,Y)\in\fz\oplus\fb$ with $Y\ne 0$. Then $(Z,Y)=\Ad(h)Y$ for some $h\in H$. Thus
$(Z,Y)$ belongs to the Lie algebra $\Ad(h)\fb$ of the stabiliser in $H$ of $hB\in H/B$.  Assume that $U$ acts properly on $H/B$. Then all stabilisers are compact, hence trivial. It follows that $(Z,Y)\not\in\fu$ and $\fu\cap\fb=\{0\}$. We obtain $\fu\cap(\fz\oplus\fb)=\fu\cap\fz$. 
If $\fu\cap(\fz\oplus\fb)=\{0\}$, then $U\backslash H/B$ is not compact. If $\fu\cap(\fz\oplus\fb)\ne\{0\}$, then 
$\fz\subset\fu$ and $U\backslash H/B\cong \fa/((\fa\cap\fu)\oplus\fb)$. The latter space is compact only if 
$(\fa\cap\fu)\oplus\fb=\fa$.
\qed

Let $\Delta\subset\RR$ be a closed subgroup, and let $\psi:\Delta\rightarrow G$ be a homomorphism such that $p\circ\psi=\id_\Delta$. Let $U\subset H$ be a connected subgroup normalised by $\psi(\Delta)$.

Let $n=\frac{1}{2}\dim\fa$ as usual. 
\begin{lm}\label{kirsten}
The group $U\cdot\psi(\Delta)$ acts properly and cocompactly on $X$ if and only if
$$\Delta\ne\{0\}, \quad\fz\subset \fu, \quad\dim(\fa\cap\fu)=n, $$
and $e^{tL}(\fa\cap\fu)\cap\fa_+=\{0\}$ for all $t\in\RR$. If these conditions are satisfied, then $U\cdot\psi(\Delta)$ acts freely on $X$. If in addition $\Delta=\RR$, then the action is simply transitive.
\end{lm}

\proof Let $t\in\RR$, and let $B_t\subset H$ be the connected subgroup with Lie algebra $e^{tL}(\fa_+)$. The latter is a Lagrange
subspace of $(\fa,\omega)$.  There exists an $H$-equivariant diffeomorphism from the fibre $\pi^{-1}(t)$ of the canonical fibration  to $H/B_t$.  
Since $\Delta\subset\RR$ is closed the action of $\Delta$ on $\RR$ is proper.

Assume now that $U\cdot\psi(\Delta)$ acts properly and cocompactly on $X$. Then $U$ acts properly on the fibres of $\pi$. 
Moreover, $\pi$ induces a continuous fibration
$\bar \pi:  (U\cdot\psi(\Delta))\backslash X\rightarrow \Delta \backslash \RR$ between Hausdorff spaces with compact total space. Hence base and fibres of the fibration $\bar\pi$ are compact, too. We conclude that $\Delta\ne 0$ and that $U$ acts properly and cocompactly on $H/B_t$ for all $t\in\RR$. Now the claimed properties of $\fu$ are implied by Condition $(iv)$ in Lemma~\ref{ulf}.

Vice versa, let us assume that the claimed properties of $\fu$ are satisfied. Then, according to Lemma~\ref{ulf}, the group $U$ acts simply transitively on all fibres of $\pi$. Since $\Delta$ acts freely (and transitively, if $\Delta=\RR$) on the base the last two assertions of the lemma follow.
We claim that the composition of natural maps $U\times\RR\hookrightarrow H\rtimes \RR\hookrightarrow G\rightarrow  G/G_+=X$ is a diffeomorphism. 
First of all, this map (let us denote it by $\Phi$) is bijective. The differential $d\Phi$ maps vectors tangent to $U$ bijectively to vectors tangent to the fibres of $\pi$. Moreover, $d\pi\circ d\Phi$ is equal to the projection $\fu\times\RR\rightarrow \RR$. Thus $d\Phi$ is bijective at all points of $U\times\RR$. The claim follows.
The induced action of $U\cdot\psi(\Delta)\ni n\psi(\delta)$ on $U\times\RR$ is of the form 
$n\psi(\delta)(u,t)=(n\Psi(\delta,u,t),\delta+t)$ for some smooth map $\Psi: \Delta\times U\times\RR\rightarrow U$. 
Let now $D\subset U$ and $E\subset\RR$ be compact subsets. Suppose that $n\psi(\delta)(D\times E)\cap (D\times E)\ne\emptyset$. Then $\delta\in E - E$ and 
$$n\in D\cdot \left(\Psi((E-E)\cap \Delta,D,E)\right) ^{-1}\ .$$
Hence $n\psi(\delta)$ belongs to a compact subset of $U\cdot\psi(\Delta)$ depending only on $D\times E$. This proves properness of the $U\cdot\psi(\Delta)$-action. If $\Delta\ne\{0\}$, then the orbit space $(U\cdot\psi(\Delta))\backslash X\cong \Delta\backslash \RR$ is compact.
\qed

\begin{co}\label{adam}
Let $\Gamma\subset G$ be a discrete subgroup acting properly and cocompactly on $X$. Then the group $C_K\subset K$ provided by Proposition \ref{nottame} is trivial.
In other words: A conjugate of a finite index subgroup of $\Gamma$ is contained in $S_\Gamma$.
\end{co}

\proof We consider $S_\Gamma=U\cdot\psi(\Delta)$. Then $C_K$ centralises $U$, see Prop.~\ref{nottame}. By Lemma~\ref{null} and Lemma~\ref{kirsten} we have $\fa=(\fu\cap\fa)\oplus\fa_+$. Hence the projection of $\fu\cap\fa$ to $\fa_-$ is a $C_K$-equivariant bijection. We conclude that $C_K$ acts trivially on $\fa_-$, i.e. $C_K$ is trivial.
\qed

We now characterise the Cahen-Wallach spaces $X$ admitting a compact quotient $Y=\Gamma\backslash X$ with 
$\overline{p(\Gamma)}=\RR$ completely. Let us first introduce some terminology.

\begin{de}\label{vier.sex} Let $X$ be a Cahen-Wallach space, and let $Y=\Gamma\backslash X$ be a compact quotient of $X$.
Then $Y$ (and $\Gamma$) is called straight, if ${p(\Gamma)}\subset\RR$ is discrete. 
\end{de}

Note that for straight compact quotients the canonical fibration $\pi:X\rightarrow \RR$ induces a fibration $\bar\pi: Y\rightarrow S^1$. Non-straight quotients  $Y$ inherit a foliation with dense leaves instead.

We formulate the result in terms of the classification of Cahen-Wallach spaces explained in Subsection~\ref{class}.

\begin{theo}\label{AA}
A Cahen-Wallach space $X$ admits a non-straight  compact quotient if and only if $X$ is of imaginary type and all the parameters $\mu_i$ of $X$ appear with
even multiplicity, i.e., $X$ is a group manifold with biinvariant Lorentzian metric (cf. Subsection~\ref{groupcase}).
All these spaces $X$ also admit straight compact quotients.
\end{theo}

\proof Let $X\cong X_{p,q}(\lambda_1,\dots,\lambda_p,\mu_1,\dots,\mu_q)$, and let $Y=\Gamma\backslash X$ be a compact quotient of $X$ that is not straight. Then  $\overline{p(\Gamma)}=\RR$.
Therefore $\Gamma\subset G$ is tame by Lemma~\ref{Ltame}. We consider $S_\Gamma=U\cdot\psi(\Delta)$. Then $\Delta=\RR$. By Proposition \ref{tame} we have
\begin{equation}\label{arno}
\fu\cap\fa \subset \fa^{r(\psi(\Delta))}\ .
\end{equation}
There is a generator of the Lie algebra of $r(\psi(\Delta))$ acting on 
$\fa$ by $L+\phi$ for some $\phi\in \fk$. Thus $\fa^{r(\psi(\Delta))}=\fa^{L+\phi}$. There exist real numbers $\beta_1,\dots,\beta_l$, $0\le 2l\le p$, and $\gamma_1,\dots,\gamma_m$, $0\le 2m\le q$,  and a reordering of the parameters 
$\lambda_i,\mu_j$ of $X$ with 
\begin{equation}\label{schmidt}
|\lambda_{2k-1}|=|\lambda_{2k}|,\quad k=1,\dots, l\ , \qquad\  |\mu_{2k-1}|=|\mu_{2k}|,\quad k=1,\dots, m
\end{equation}
such that the eigenvalues of $L+\phi$ are given by 
$$\lambda_1+i\beta_1,\ \lambda_1-i\beta_1,\dots, \lambda_{2l-1}+i\beta_l,\ \lambda_{2l-1}-i\beta_l,\ \lambda_{2l+1},\dots,\lambda_p$$
multiplied by $\pm 1$ and 
$$\mu_1+\gamma_1,\ \mu_1-\gamma_1,\dots, \mu_{2m-1}+\gamma_m,\ \mu_{2m-1}-\gamma_m,\ \mu_{2m+1},\dots,\mu_q$$
multiplied by $\pm i$. It follows that $\dim \fa^{r(\psi(\Delta))}\le 2m$. Thus $\dim \fu\cap\fa\le 2m$ by (\ref{arno}). On the other hand, Lemma \ref{kirsten} tells us that $\dim \fu\cap\fa=p+q$. Hence $2m=q$ and $p=0$. By (\ref{schmidt}) all $\mu_i$ have even multiplicity.

Vice versa, let $X=X_{0,2m}(\mu_1,\mu_1,\cdots,\mu_m,\mu_m)$, $2m=n$, be a Cahen-Wallach space of imaginary type with parameters of even multiplicity.  Put $\phi:=\phi_\mu$ as defined in (\ref{phi}) for $\mu=(\mu_1,\dots,\mu_m)$. Then $\phi\in\fk$. Moreover, $\fa^{L+\phi}\subset(\fa,\omega)$ is a symplectic subspace of dimension $n$ complementary to $\fa_+$.
In particular,  $\fu:=\fz\oplus\fa^{L+\phi}$ is isomorphic to a Heisenberg algebra. We consider the corresponding connected group $U\subset G$.
Define $\psi: \RR\rightarrow G$ by $\psi(t):=(0,\exp(t\phi),t)$ and put $S:=U\cdot\psi(\RR)$. By Lemma~\ref{kirsten} 
the group $S$ acts
simply transitively on $X$. Every lattice $\Gamma\subset S$ acts freely, properly discontinuously and cocompactly on $X$ and thus defines a compact quotient $Y_\Gamma$ of $X$. Now it becomes important that $\psi(\RR)$ centralises $U$, i.e. that $S\cong U\times\RR$. 
The quotient $Y_\Gamma$
is straight if and only if $\Gamma$ has a discrete projection on the $\RR$-factor of $S$.
Lattices $\Gamma\subset U\times\RR$ of both types do really exist. For example, let $\Gamma_0$ be a lattice of the Heisenberg group $U$, and let 
$\chi: U\rightarrow \RR$ be a group homomorphism such that $\chi(\Gamma_0)\not\subset\QQ$.  Then 
$$\Gamma_1=:\Gamma_0\times \ZZ\quad \mbox{ and }\quad 
\Gamma_2:=\{(\gamma_0, \chi(\gamma_0)+k)\mid \gamma_0\in\Gamma_0, k\in\ZZ\}$$ 
are lattices in $S$. The group $\Gamma_1\subset G$ is straight, $\Gamma_2$ is not straight.
\qed 

The theorem says in particular that all group manifolds $X=Q$ among the Cahen-Wallach spaces admit compact quotients.
This is in contrast to the fact that not all oscillator groups $Q$ admit lattices  (see \cite{MR2}). 
The group $S\cong U\times\RR\cong H_m\times\RR$ that was crucial in the above proof 
has the following nice alternative description.
Namely, $Q=H_m\rtimes \RR$, and $H_m\times \RR$ acts isometrically and simply transitively
on $Q$ via $(h,t)q:= hqt^{-1}$.

Now we derive a criterion for the existence of general
compact quotients. We will say that a subset $M_1$ of a set $M_2$ is 
{\em stable} under a map $f:M_2\to M_2$ if $f(M_1)=M_1$.

\begin{pr}\label{wulf}
Let $X=X(\omega,\theta_\fa,L)$ be a Cahen-Wallach space of dimension $n+2$. 
Then $X$ has a compact quotient if and only if there exist
\begin{enumerate}
\item[(a)] an $n$-dimensional subspace $V\subset \fa$ such that $e^{tL}V\cap\fa_+=\{0\}$ for all $t\in\RR$;
\item[(b)] elements $t_0\in\RR\setminus\{0\}$, $\ph_0\in K$, $h_0\in H^{t_0\ph_0}$ and a lattice $\Lambda$ of the subgroup $\fz\oplus V\subset H$ that is stable under conjugation by $h_0t_0\ph_0$.
\end{enumerate}
If $X$ is of imaginary type, then (b) can be replaced by
\begin{enumerate}
\item[(b')] elements $t_0\in\RR\setminus\{0\}$, $\ph_0\in K$ such that $(e^{t_0L}\ph_0)|_{ V}=\id_V$.
\end{enumerate}
\end{pr}
\proof The group $\fz\oplus V$ always contains a lattice. Therefore $(b')$ implies $(b)$ with $h_0=0$.

Assume now that $(a)$ and $(b)$ are satisfied. Let $\gamma_0=h_0 t_0 \ph_0\in H\rtimes (\RR\times K)=G$. We consider the discrete subgroup $\Gamma:=\langle\Lambda,\gamma_0\rangle\subset G$. Then $S_\Gamma=\langle \fz\oplus V,\gamma_0\rangle$.
By Lemma~\ref{kirsten} the group $S_\Gamma$ acts freely, properly and cocompactly on $X$. Hence, so does $\Gamma$, and
$Y=\Gamma\backslash X$ is a compact quotient.

Vice versa, assume that $X$ admits a compact quotient. Then, by Theorem \ref{AA} and Corollary~\ref{adam}
the space $X$ admits a straight compact quotient $Y=\Gamma\backslash X$ such that
$$\Gamma\subset S_\Gamma =: U\cdot\psi(\langle t_0\rangle)$$ 
for some $t_0\in\RR$. We may assume that $\psi(t_0)\in\Gamma$. Define $\Lambda:=\Gamma\cap U$. 
The group $\Lambda$ is a lattice in $U$ that is invariant under conjugation by $\psi(\langle t_0\rangle)$. Let $h_0, \ph_0$ defined by $H\rtimes(\RR\times K)\ni h_0t_0 \ph_0=\psi(t_0)$. By Prop.~\ref{nottame} we have $h_0\in H^{t_0\ph_0}$. 
The group $S_\Gamma$ acts properly and cocompactly on $X$. Now Lemma~\ref{kirsten} tells us that $U=\fz\oplus V$ for some
$V\subset \fa$ and that the remaining conditions in $(a)$ and $(b)$ are satisfied. If $X$ is of imaginary type, then $\Gamma$ is tame, see Lemma~\ref{Ltame}. By Prop.~\ref{tame} we can assume that $\psi(t_0)$ centralises $U$. Condition $(b')$ follows. 
\qed

In the following sections we try
to make the criteria provided by Proposition~\ref{wulf} as explicit as possible in terms 
of the parameters $(\lambda,\mu)$ of $X$. We do this separately for spaces of real, imaginary and mixed type. 
The most
complete understanding will be achieved for spaces of real type, because in this case the intricate condition
\begin{equation}\label{harald}
e^{tL}V\cap\fa_+=\{0\}\qquad \mbox{for all }t\in\RR
\end{equation}
turns out to be a consequence of the remaining criteria and $V\cap\fa_+=\{0\}$.

\begin{re}\label{CWfalsch} {\rm 
In \cite{CW}, Section 4, Cahen and Wallach claim to construct examples of compact quotients of Lorentzian symmetric spaces. They consider the symmetric spaces $X_{0,n}(2\pi p_1/q_1,\dots, 2\pi p_n/q_n)$, where $p_i, q_i\in \ZZ_{\not= 0}$. Of course we may assume $q_1=\dots=q_n=1$. For each such space $X$ and an arbitrary lattice $\Gamma_0$ of the abelian group $\fa_-$ they consider the discrete subgroup $\Gamma=\ZZ\times \Gamma_0\times \ZZ\subset \fz\times\fa_-\times\RR\subset\hat G$  of the transvection group and they state that $\Gamma$ acts properly discontinuously and cocompactly on $X$. However, this is not correct. The proof of Proposition~\ref{wulf} shows that, if the action were proper, the space $V:=\fa_-$ would satisfy $e^{tL}V\cap \fa_+=0$ for all $t\in\RR$, which is obviously not true.  This can also be seen without using the preceding investigations. Indeed, if the action of $\Gamma$ were proper then also the action of the subgroup $\fa_-\subset \hat G$ on $X$ would be proper since $\Gamma_0$ is a lattice in $\fa_-$. Now take $v\in\fa_-$ and $t\in\RR$ such that $e^{-tL}v\in\fa_+$. Abbreviating as usual $(0,0,t)\in \hat G$ to $t$ we get 
$$(0,sv)\cdot t \cdot \hat G_+=t\cdot(0,e^{-tL}(sv))\cdot \hat G_+=t\cdot \hat G_+\in \hat G /\hat G_+=X$$
for all $s\in \RR$. Hence the line $\RR\cdot v\subset\fa_-$ is contained in the stabiliser of $t\cdot \hat G_+\in X$, which contradicts properness.
}
\end{re}

\section{The real case}\label{neukirch}

We will need the following well-known fact. For convenience of the reader, we include its short proof here.

\begin{lm}\label{gigi}
Let $V$ be a real vector space, and let $A\in \GL(V)$ have characteristic polynomial $f_A$. We consider the following two assertions:
\begin{enumerate}
\item[(a)] There exists a lattice $\Lambda\subset V$ stable under $A$.
\item[(b)] $f_A\in\ZZ[x]$ and the constant term of $f_A$ has absolute value $1$.
\end{enumerate} 
Then $(a)$ implies $(b)$. If $A$ is semisimple, then $(a)$ and $(b)$ are equivalent.
\end{lm}

\proof If $A$ stabilises a lattice, then it can be represented by a matrix in $\GL(n,\ZZ)$, hence $f_A\in\ZZ[x]$ and $|\det(A)|=1$.

We now assume that $A$ is semisimple and satisfies $(b)$. We first show that for each $f\in\ZZ[x]$ whose constant term is $\pm1$ there is a semisimple matrix $C_f\in \GL(n,\ZZ)$ having characteristic polynomial $f$. We factorise $f$ into monic polynomials irreducible in $\QQ[x]$:  $f=f_1\cdot\dots\cdot f_r$. The roots of $f_i$, $i=1,\dots,r$, are simple. Let 
$C_{f_i}$ be the companion matrix of $f_i$.
Then $C_{f_i}$ has characteristic polynomial $f_i$. In particular, all eigenvalues of $C_{f_i}$ have algebraic multiplicity one, thus $C_{f_i}$ is semisimple. Now the block diagonal matrix $C_f=\diag(C_{f_1},\dots, C_{f_r}) \in \GL(n,\ZZ)$ has the required properties.

Choose a basis $e_1,\dots,e_n$ of $V$. Then the matrix $C_{f_A}$ defines an element $B\in \GL(V)$ that stabilises the lattice $\Lambda_0$ spanned by $e_1,\dots,e_n$. Now $A$ and $B$ are semisimple and their eigenvalues over $\CC$ (with multiplicities) coincide. Thus there  exists an $M\in \GL(V)$ such that $A=MBM^{-1}$. Thus $A$ satisfies $(a)$ with $\Lambda=M\Lambda_0$.
\qed

\begin{theo}\label{BB}
Let $X$ be an $(n+2)$-dimensional Cahen-Wallach space of real type. Then $X$ admits a compact quotient if and only if there exists a polynomial
$f\in\ZZ[x]$, 
\begin{equation}\label{otto}
f(x)= x^n+a_{n-1}x^{n-1}+\dots+a_1x\pm 1\ ,
\end{equation}
with no roots on the unit circle such that 
$$X\cong X_{n,0}(\log |\nu_1|,\log |\nu_2|,\dots,\log |\nu_n|)\  ,$$
where $\nu_1,\nu_2,\dots,\nu_n$ are the roots of $f$.
\end{theo}

\proof Let $f\in\ZZ[x]$ be as in the theorem. Let $2s$ be the number of non-real roots of $f$.  We index the roots in a way such that $$\nu_1,\dots,\nu_{2s}\not\in\RR,\quad 
\bar\nu_{2k-1}=\nu_{2k},\  k=1,\dots,s\quad\mbox{and}\quad \nu_l\in\RR^* \mbox{ for }l=2s+1,\dots,n. $$ 
We consider the Cahen-Wallach space $X:=X_{n,0}(\log |\nu_1|,\dots,\log |\nu_n|)$ and want to construct data $V,\Lambda,t_0,\ph_0,h_0$ as required by Proposition~\ref{wulf}. 

We use the conventions of Example \ref{Exosc2}. In particular, $\fa=\CC^n$ with standard basis $e_1,\dots,e_n$ and $\theta_\fa$ is given by complex conjugation of the coordinates. We define $V:=(1+i)\RR^n$.
Then $V\cap\fa_+=0$.
Furthermore, the vectors $(1+i)e_l$ are eigenvectors of $L$ with corresponding eigenvalue $\log |\nu_l|$. In particular, $V$ is $L$-invariant.
Hence Condition~$(a)$ of Prop.~\ref{wulf} is satisfied.

We define a linear operator $\ph_0$ on $\CC^n$ by a real diagonal block matrix, where the blocks are
\begin{equation}\label{bella}
\left (
\begin{array}{cc}
\Re\left(\frac{\nu_{2k-1}}{|\nu_{2k-1}|}\right)&-\Im\left(\frac{\nu_{2k-1}}{|\nu_{2k-1}|}\right)\\
\Im\left(\frac{\nu_{2k-1}}{|\nu_{2k-1}|}\right)&\Re\left(\frac{\nu_{2k-1}}{|\nu_{2k-1}|}\right)
\end{array}
\right)\ ,\quad k=1,\dots, s,
\end{equation}
and the diagonal matrix $\displaystyle\diag\left(\frac{\nu_{2s+1}}{|\nu_{2s+1}|},\dots,\frac{\nu_{n}}{|\nu_{n}|}\right)$. Then $\ph_0\in K$, and $\ph_0$ leaves $V$ invariant. The complex eigenvalues of the semisimple operator $A:=e^L{\ph_0}|_{V}$ are precisely $\nu_1,\dots,\nu_n$. Thus $A$ has characteristic polynomial $f$. By Lemma~\ref{gigi} there is a lattice $\Lambda_0\subset V$ stabilised by $A$. Choose an arbitrary lattice $\Lambda_1\subset\fz$ and put 
$\Lambda:=\Lambda_1\oplus\Lambda_0\subset \fz\oplus V$. Since $V\subset(\fa,\omega)$ is a Lagrange subspace the group structure of $\fz\oplus V$ viewed as a subgroup of
$H$ is just given by vector space addition. Therefore, $\Lambda\subset \fz\oplus V$ is really a subgroup, hence a lattice. Now we see that Condition~$(b)$
of Prop.~\ref{wulf} is satisfied with $t_0=1$, $h_0=0$. It follows that $X$ admits a compact quotient.

Vice versa, let $X$ be an $(n+2)$-dimensional Cahen-Wallach space of real type that has a compact quotient. Let $V,\Lambda,t_0,\ph_0,h_0$ be as in Prop.~\ref{wulf}, Conditions~$(a)$ and $(b)$. 
For all Cahen-Wallach spaces of real type and arbitrary $0\ne t\in \RR$, $\ph\in K$ we have $H^{t\ph}=\fz$. 
In particular, $h_0$ acts trivially on $\fz\oplus V$ and we can assume $h_0=0$. 
We claim that $\Lambda\cap\fz\ne 0$. This is clear if $\fz\oplus V$ is non-abelian. In the abelian case we consider the linear operator $B:=1-e^{t_0L}\ph_0$ on $\fz\oplus V$. Its kernel is $\fz$. The lattice $\Lambda$ generates a $\QQ$-vector space $(\fz\oplus V)_{\Bbb Q}$ on which $B$ still acts with a one-dimensional kernel. A suitable rational multiple of a non-zero element of this kernel belongs to $\Lambda\cap\fz$.

By the claim, the projection $\Lambda_0$ of $\Lambda$ to $V$ is a lattice. Set $A:=e^{t_0L}{\ph_0}$ and $A_0:=A|_{V}$. Then $\Lambda_0$ is $A_0$-stable.
By Lemma~\ref{gigi} the characteristic polynomial $f_{A_0}$ is integral and the constant term of $f_{A_0}$ has absolute value $1$.
Let us denote its roots by $\nu_1,\dots,\nu_n$.

Since $L$ and $A$ commute the eigenspaces $E_\lambda$ of $L$ are $A$-invariant. All the complex eigenvalues of $A$ on $E_\lambda$ have absolute value
$e^{t_0\lambda}$. If $\lambda\ne \lambda'$, then $e^{t_0\lambda}\ne e^{t_0\lambda'}$. The $A$-invariant subspace $V\subset\fa$ can be decomposed
into a direct sum of invariant subspaces characterized by the absolute value of the $A$-eigenvalues appearing. The previous discussion implies that such a subspace is contained in a single $L$-eigenspace $E_\lambda$. We conclude that $V$ is $L$-invariant and that the $L$-eigenvalues $\lambda_i$ on $V$
are related to the $A_0$-eigenvalues $\nu_i$ by 
\begin{equation}\label{manfred}\lambda_i=\frac{\log |\nu_i|}{t_0}, \quad i=1,\dots,n.
\end{equation}
In particular, no eigenvalue of $A_0$ lies on the unit circle. Thus the polynomial $f_{A_0}$ has all the properties required by the theorem.

We consider the $L$ and $\theta_\fa$-invariant subspace $V\cap\theta_\fa V$. Since $V\cap\fa_+=\{0\}$ it is contained in $\fa_-$. But $L$ maps $\fa_-$ bijectively to $\fa_+$. Hence $V\cap\theta_\fa V=\{0\}$. We conclude that $\fa=V\oplus\theta_\fa V$. This implies that the eigenvalues of $L$ on $\fa$
are precisely 
$$\lambda_1,\dots,\lambda_n,-\lambda_1,\dots,-\lambda_n\ ,$$ 
where $\lambda_i$ is given by (\ref{manfred}). It follows that
$$X\cong  X_{n,0}(\lambda_1,\dots,\lambda_n)\cong X_{n,0}(\log |\nu_1|,\dots,\log |\nu_n|)\ .$$
\qed

The theorem has the following immediate consequences.

\begin{co} \label{countable}
The set of isometry classes of Cahen-Wallach spaces of real type admitting a compact quotient is countable.
\end{co}
\proof According to Theorem~\ref{BB} there is a subset of $\ZZ[x]$ that surjects to the set of isometry classes in question.
\qed

\begin{co}\label{spurnull}
Assume that $X_{n,0}(\lambda)$ admits a compact quotient. Then there is choice of signs such that
\begin{equation}\label{karl-heinz}
\displaystyle \sum_{i=1}^{n} \pm\lambda_i = 0\ .
\end{equation}
\end{co}

\proof Theorem~\ref{BB} implies that for some choice of signs $\pm\lambda_i= c\log(|\nu_i|)$, where the constant $c$ is independent of $i$, and
$\nu_i$ are the roots of a monic polynomial with constant term $\pm 1$. Thus $1=|\nu_1\nu_2\dots\nu_n|=|\nu_1|\dots|\nu_n|$. Taking logarithms the corollary follows.
\qed

For a given Cahen-Wallach space $X_{n,0}(\lambda)$ satisfying (\ref{karl-heinz}) it might be rather difficult to decide by a direct application of Theorem~\ref{BB} whether it admits a compact quotient. There is, however, a second necessary condition, which is relatively easy to check.

\begin{pr}\label{transc}
Assume that $X_{n,0}(\lambda)$ admits a compact quotient. Then the quotients $\lambda_i/\lambda_j$, $1\le i,j\le n$, are either
rational or transcendental.
\end{pr}

\proof By Theorem~\ref{BB} the number $\lambda_i/\lambda_j$ is a quotient of logarithms of algebraic numbers.
By  the Gelfond-Schneider Theorem (which solves Hilbert's seventh problem) such quotients are either rational
or transcendental, see e.g. \cite{Sch}, Satz~14 or \cite{Wa}, Thm.~1.4.
\qed
 
In Proposition~\ref{transc}, not too many of the ratios $\lambda_i/\lambda_j$ can be a rational number different from $\pm 1$.
For instance, we have

\begin{pr}\label{werwolf}
Let $X$ be a Cahen-Wallach space isomorphic to $$X_{n,0}(l_1,\dots,l_d,\lambda_{d+1},\dots,\lambda_n)\ ,$$ 
such that $l_i\in\QQ$, $|l_i|\ne |l_j|$ for $i\ne j$, $i,j=1,\dots,d$, and the sets $S_i:=\{k\ge d+1\mid |\lambda_k|= |l_i|\}$ have even cardinality for
$i=1,\dots,d$.
If $X$ admits a compact quotient, then $d\le \frac{n}{3}$.
\end{pr}

\begin{lm}\label{weswolf}
Let $\nu_1,\nu_2\in\RR\setminus\{-1,0,1\}$ be real roots of an irreducible polynomial $f\in\QQ[x]$. Assume that 
$\displaystyle \alpha:=\frac{\log |\nu_1|}{\log |\nu_2|}\in\QQ$. Then $\alpha=\pm 1$.
\end{lm}

\proof
We look at the Galois group $G(f)$ of $f$, i.e. $G(f)=\mathrm{Gal}(N/\QQ)$, where the field $N\subset\CC$ arises by adjoining all roots of 
$f$ to $\QQ$. Since $f$ is irreducible, there exists $g\in G(f)$ such that $g(\nu_1)=\nu_2$. We write $\alpha=\frac{p}{q}$, $p,q$ being coprime integers. We set $\rho_i:=\nu_i^2\in N\cap\RR^+\setminus\{1\}$. Then
$
g(\rho_1)=\rho_2
$
and $\rho_1^q=\rho_2^p$. We choose $m,n\in\ZZ$ such that $np+mq=1$ and define $\rho:=\rho_1^n\rho_2^m\in N\cap\RR^+\setminus\{1\}$.
This implies $\rho_1=\rho^p$, $\rho_2=\rho^q$, and
$$
g(\rho)^p=\rho^q\ .
$$
By induction we obtain for every $k\ge1$
$$ g^k(\rho)^{p^k}=\rho^{q^k}\ .$$
It follows that $\rho^{p^l}=\rho^{q^l}$, where $l$ is the order of $g$. Hence $p=\pm q$, i.e. $\alpha=\pm 1$.
\qed

{\sl Proof of Prop.~\ref{werwolf}. } Assume that $X_{n,0}(l_1,\dots,l_d,\lambda_{d+1},\dots,\lambda_n)$ admits a compact quotient.
Theorem \ref{BB} implies that there exist a polynomial $f\in\ZZ[x]$ of the form (\ref{otto}) having at least $d$ real roots 
$\nu_1,\dots,\nu_d$, a constant $c\in\RR^+$ and a choice of signs such that $\log |\nu_i|=\pm c l_i$ for $i=1,\dots,d$. 
In particular, $\log |\nu_i|/\log |\nu_j|\ne \pm 1$ for $i\ne j$. For $i=1,\dots,d$ let $f_i$ be the minimal polynomial of 
$\nu_i$ over $\QQ$.
Lemma \ref{weswolf} implies that $f_i(\nu_j)\ne 0$ for $j\in\{1,\dots,d\}\setminus\{i\}$. We conclude that $f$ is a product of at least $d$
irreducible factors. 
Hence it suffices to prove that $f_i$ has degree at least 3. Assume that $f_i$ were quadratic. Then the second root $\nu_i'$ of $f_i$ satisfies $\log |\nu_i'|=\log (1/|\nu_i|) = -\log |\nu_i|=\pm cl_i$, which contradicts the fact that $l_i$ (seen as a parameter of $X$) has odd multiplicity.
\qed
  
There are two other ways to look at the numbers $\nu_1,\dots,\nu_n$ appearing in Theorem \ref{BB}. Geometrically, they are just the eigenvalues of strictly hyperbolic matrices in $\GL(n,\ZZ)$. Conjugacy classes of hyperbolic elements of 
$\PGL(n,\ZZ)$ define closed geodesics on the 
Riemannian locally symmetric space $S:=\PGL(n,\ZZ)\backslash \PGL(n,\RR)/\PO(n)$. 
Therefore the set of all parameters of the form $(\log |\nu_1|,\log |\nu_2|,\dots,\log |\nu_n|)$ giving $(n+2)$-dimensional Cahen-Wallach of real type admitting compact quotients can be interpreted as a kind of multi-dimensional length spectrum
of the geodesic flow of $S$. The investigation of geodesic flows of locally symmetric spaces like $S$ is a classical and 
still very fruitful area of mathematical research at the edge between ergodic theory, number theory, and harmonic analysis,
see e.g. \cite{ELMV}.
On the other hand, 
from a purely number theoretic point of view,
$\nu_1,\dots,\nu_n$ are just (conjugate) units in algebraic number fields. We will now exploit this point of view to get more information on the set of isometry classes of Cahen-Wallach spaces of real type admitting compact quotients.


Let $K$ be an algebraic number field of degree $d$, i.e., a field extension of $\QQ$ of degree~$d$. Let $\fo_K\subset K$ be its ring of integers, and let
$\fo_K^*\subset\fo_K$ be the group of units of $K$. Note that $\fo_K^*$ is just the set of zeroes in $K$ of polynomials of the form
(\ref{otto}). Let $\Hom(K,\CC)=\{\tau_1,\dots,\tau_{2s},\tau_{2s+1},\dots,\tau_d\}$ be the set of embeddings $K\hookrightarrow\CC$, where we assume that $\tau_{2j}$ is the composition of $\tau_{2j-1}$ with complex conjugation for $j=1,\dots, s$ and that the
image of $\tau_i$ is contained in $\RR$ for $i\ge 2s+1$. If $K$ is given by the irreducible polynomial $f$, i.e., $K=\QQ[x]/(f)$, then
the embeddings correspond to the $d$ different roots $\nu_1,\dots,\nu_d$ of $f$ in $\CC$: $\tau_i(x)=\nu_i$.
A number field $K$ is called totally real if $s=0$. If $d=2,3$, one just says `real' quadratic or cubic field, instead.

For $0\le s\le [d/2]$ we define a vector subspace $\RR^d_s\subset\RR^d$ of dimension $d-1-s$ by
$$ \RR^d_s:=\left\{y=(y_1,\dots,y_d)\in\RR^d\mid \sum_{i=1}^{d} y_i=0,\ y_{2j-1}=y_{2j}\mbox{ for } j=1,\dots,s\right\}\ .$$
We consider the group homomorphism $l_K: \fo_K^*\rightarrow \RR^d_s$ given by
$$ l_K(\nu):=(\log |\tau_1(\nu)|,\dots,\log |\tau_d(\nu)|)\ .$$
It is a fundamental fact that the image of $l_K$ is a lattice in the $\RR$-vector space $\RR^d_s$ and that its kernel is precisely
the group $\mu_K$ of roots of unity in $K$ (see e.g. \cite{Neu} or \cite{FT}; this fact is the main ingredient of the modern proof
of Dirichlet's unit theorem stating that $\fo_K^*$ is the product of the finite cyclic group $\mu_K$ and a free abelian group of rank
$d-1-s$).

\begin{de}\label{blasebalg}
Let $K$ be a number field of degree $d$ admitting precisely $s$ pairs of complex conjugate embeddings $K\hookrightarrow\CC$.
We define a $\QQ$-form $\cH_K$ of the $\RR$-vector space $\RR^d_s$ as the vector subspace over $\QQ$ of $\RR^d_s$
generated by the lattice $\im l_K$.
\end{de}

For a subset $A\subset\RR^d$ we set $A^\reg:=A\cap (\RR^*)^d$.

The following theorem is essentially a reformulation of Theorem \ref{BB}.

\begin{theo}\label{BBB}
Let $X$ be an $(n+2)$-dimensional Cahen-Wallach space of real type. Then $X$ admits a compact quotient if and only if there exist a collection of number fields $K_1,K_2,\dots, K_r$ of degree $d_i$ (and corresponding data $s_i$ and $\cH_{K_i}$ as in Def.~\ref{blasebalg}) satisfying $d_i-s_i\ge 2$,
$ \sum_{i=1}^r d_i=n$
and vectors $\lambda^i\in\cH_{K_i}^\reg\subset\RR^{d_i}$ such that
$$ X\cong X_{n,0}(\lambda^1,\lambda^2,\dots,\lambda^r)\ .$$
\end{theo}

\proof
Assume that $X$ admits a compact quotient. By Theorem~\ref{BB} there exists a polynomial $f\in\ZZ[x]$ of the form (\ref{otto}) with roots $\nu_i$, $|\nu_i|\ne 1$, such that 
\begin{equation}\label{oberwolf}
X\cong X_{n,0}(\log |\nu_1|,\log |\nu_2|,\dots,\log |\nu_n|)\ .
\end{equation} 
We decompose $f=f_1\cdot\dots\cdot f_r$ into irreducibles over $\QQ$. Set $d_i:=\deg f_i$, $k_i:=d_1+d_2+\dots+d_{i-1}$, and let $s_i$ be the number of pairs of complex conjugate non-real roots of $f_i$. We order the roots of $f$ such that
$$  \nu_{k_i+1},\nu_{k_i+2},\dots,\nu_{k_i+2s_i},\nu_{k_i+2s_i+1},\dots,\nu_{k_i+d_i}$$
are the roots of $f_i$ and $\nu_{k_i+2j}=\overline{\nu_{k_i+2j-1}}$, $j=1,\dots,s_i$. We set $K_i:=\QQ(\nu_{k_i+1})$. Then 
$\nu_{k_i+1}$ is a unit in $K_i$, and the embeddings $\tau^i_j\in\Hom(K_i,\CC)$ are given by $\tau^i_j(\nu_{k_i+1})=\nu_{k_i+j}$.
We define
$$\lambda^i:=(\log |\nu_{k_i+1}|,\dots,\log |\nu_{k_i+d_i}|)=l_{K_i}(\nu_{k_i+1})\in\cH_{K_i}\ .$$
Since none of the roots lies on the unit circle, 
we have in fact $\lambda^i\in\cH_{K_i}^\reg$.
Now (\ref{oberwolf}) shows that $X\cong X_{n,0}(\lambda^1,\lambda^2,\dots,\lambda^r)$.

Vice versa, assume that for $i=1,\dots, r$ number fields $K_i$ of degree $d_i$,  $ \sum_{i=1}^r d_i=n$,
and vectors $\lambda^i\in\cH_{K_i}^\reg$ are given. Then there exists $m_i\in\NN$ such that $m_i\lambda^i\in\im l_{K_i}$.
Let $m$ be a common multiple of $m_1,\dots,m_r$. It follows that there are units $\rho_i\in K_i$ such that
$$ l_{K_i}(\rho_i)= m\lambda^i,\quad i=1,\dots r. $$

We define
$$f_i(x):=\prod_{\tau\in\Hom(K_i,\Bbb C)} (x-\tau(\rho_i)),\qquad f:=f_1\cdot f_2\cdot\dots\cdot f_r. $$
Since $m\lambda^i\in(\RR^{d_i})^\reg$ none of the zeroes $\nu_1,\dots,\nu_n$ of $f$ lies on the unit circle. By Galois theory, $f_i$ is just a power
of the minimal polynomial of $\rho_i$ over $\QQ$. Since $\rho_1,\dots,\rho_r$ are units 
it follows that $f$ is of the form (\ref{otto}). Now Theorem~\ref{BB} implies that
$$X_{n,0}(\lambda^1,\dots,\lambda^r)\cong X_{n,0}(m\lambda^1,\dots,m\lambda^r)=
X_{n,0}(\log |\nu_1|,\dots,\log |\nu_n|)$$
admits a compact quotient.
\qed

Now we consider the space $\cM_{n,0}$ of isometry classes of $(n+2)$-dimensional Cahen-Wallach spaces of real type. It comes with a continuous
surjection $\Phi_{n,0}: (\RR^*)^n\rightarrow  \cM_{n,0}$ that restricts to a homeomorphism
\begin{equation}\label{conv}
\{\lambda \in\RR^n\mid 1=\lambda_1\le\lambda_2\le\dots\le\lambda_n\}\cong  \cM_{n,0}\ ,
\end{equation}
see Subsection~\ref{class} for all that. Let $\cM_{n,0}^c\subset\cM_{n,0}$
be the subspace of isometry classes of spaces admitting a compact quotient.
We define $\cM_{n,0}^0:=\Phi_{n,0}((\RR^n_0)^\reg)\subset\cM_{n,0}$. 

\begin{co} \label{dense}
The set $\cM_{n,0}^c$ is a countable and dense subset of $\cM_{n,0}^0$.
\end{co}

\proof That $\cM_{n,0}^c$ is countable and contained in $\cM_{n,0}^0$ is the content of Corollaries~\ref{countable} and 
\ref{spurnull}, respectively. We have to prove the density statement.

Let $K$ be a totally real number field of degree $n$ (such number fields always exist). Then $\cH_K$ is a $\QQ$-form of the real vector space $\RR^n_0$. It follows that $\cH_K\subset \RR^n_0$
as well as $\cH_K^\reg\subset (\RR^n_0)^\reg$ are dense. Hence $\Phi_{n,0}(\cH_K^\reg)\subset\cM_{n,0}^0$ is dense.
But Theorem~\ref{BBB} tells us that $\Phi_{n,0}(\cH_K^\reg)\subset\cM_{n,0}^c$.
\qed

We remark that $\cH_K^\reg=\cH_K\setminus\{0\}$ for a totally real number field $K$.

Let us describe the spaces $\cM_{n,0}^c$ for $n\le 3$ more explicitly. Since $\cM_{1,0}^0=\emptyset$, we have $\cM_{1,0}^c=\emptyset$. The space $\cM_{2,0}^0$ consists of one point, the isometry class of $X_{2,0}(1,1)$. The latter space admits
a compact quotient (take an abitrary real quadratic field in Theorem~\ref{BBB} or the polynomial $x^2+kx+1$, $|k|>2$, in Theorem~\ref{BB}).
Thus $\cM_{2,0}^c$ is a singleton. The case $n=3$ is more interesting. Formula (\ref{conv}) provides a convenient parametrisation
$\phi$ of $\cM_{3,0}^0$ by the interval $[1,\infty)$:
$\phi (\lambda):=\Phi_{3,0}(1,\lambda,\lambda+1)$. We want to understand the subset $\phi^{-1}(\cM_{3,0}^c)\subset [1,\infty)$. 
We define a map $\psi:(\RR^3_0)^\reg\rightarrow [1,\infty)$ by $\psi(x):=\displaystyle\max_{i,j} \left(|x_i|/|x_j| \right)-1$.
Note that $\phi\circ\psi=\Phi_{3,0}|_{({\Bbb R}^3_0)^\reg}$.

A complete understanding of $\phi^{-1}(\cM_{3,0}^c)$ will depend on the validity of the still unproven {\bf four exponentials conjecture} in
transcendental number theory (see e.g. \cite{Wa}, Conj.~1.13), which we now state in a form convenient for our purposes:

Let $(\lambda_{ij})$ be a $(2\times 2)$-matrix of complex numbers such that $e^{\lambda_{ij}}$ is algebraic for all $i,j$ and
such that rows and columns are linearly independent over $\QQ$. Then its rows (and hence the columns) are linearly independent over $\CC$.

\begin{pr}\label{duli}
For a real cubic field $K$, we set $\Lambda_K:=\psi(\cH_K^\reg)$. Then $\Lambda_K$ is a countable dense subset of $[1,\infty)$ consisting of transcendental numbers and
\begin{equation}\label{holger}  \phi^{-1}(\cM_{3,0}^c)=\{1\}\cup\bigcup_K \Lambda_K\ ,
\end{equation}
where the union is taken over all isomorphism classes of real cubic fields. If the four exponentials conjecture is true,
then the union is disjoint.
\end{pr}

\proof For $n=3$ each collection of number fields appearing in Thm.~\ref{BBB} consists of a single cubic field. Thus
$\cM_{3,0}^c$ is the union of the sets $\Phi_{3,0}(\cH_K^\reg)$, $K$ running over all cubic fields. If $K$ is not real (i.e. $s=1$),
then $\Phi_{3,0}(\cH_K^\reg)$ consists of the single point $\phi(1)$, while for real fields we have by construction $\Phi_{3,0}(\cH_K^\reg)=\phi(\Lambda_K)$.  This proves (\ref{holger}). For the density of $\Lambda_K$ we refer to the proof of Cor.~\ref{dense}. Proposition~\ref{werwolf} implies $\Lambda_K\cap\QQ=\emptyset$. Thus, by Prop.~\ref{transc} the set $\Lambda_K$ consists entirely of transcendental numbers.

It remains to prove the last assertion of the proposition. 
Let $K_1$, $K_2$ be two non-isomorphic real cubic fields, and let $\lambda^i=(\lambda_{i1},\lambda_{i2},\lambda_{i3})\in \cH_{K_i}^\reg$, $i=1,2$. Reindexing the embeddings $K_i\hookrightarrow \RR$, if necessary, we may assume that $|\lambda_{i1}|\le|\lambda_{i2}|\le|\lambda_{i3}|$. Then the transcendental number $\psi(\lambda^i)$ equals $|\lambda_{i2}|/|\lambda_{i1}|$.
If $|\lambda_{11}|/|\lambda_{21}|$ were rational, then $K_1$ and $K_2$ would contain a common unit of infinite order. This is impossible since $K_1\cap K_2=\QQ$. We conclude that the matrix $(|\lambda_{ij}|)_{i,j=1,2}$ satisfies the assumptions of 
the four exponentials conjecture. Its validity would imply that the vectors $(|\lambda_{11}|,|\lambda_{12}|)$ and $(|\lambda_{21}|,|\lambda_{22}|)$ are linearly independent over $\CC$, i.e. $\psi(\lambda^1)\ne\psi(\lambda^2)$. 
\qed

A table of the first 100 real cubic fields $K$ (ordered by the size of their discriminant) including their fundamental units can be found
in \cite{Co}, Table B.4. Using these data, $\Lambda_K$ can be computed explicitly for those number fields.

\section{Good subspaces and the imaginary case} \label{pudelskern}
\subsection{Good subspaces and admissible tuples}\label{wagner}
Now we turn to spaces of imaginary type. In order to deal with the crucial condition (\ref{harald}) we introduce the following notion.

\begin{de}\label{johnny}
Let $(W,\theta)$ be a real $2n$-dimensional vector space with involution such that $\dim W_\pm=n$, and let 
$\rho: U(1)\rightarrow \GL(W)$ be a representation satisfying
\begin{equation}\label{b.good}
\rho(z)\circ \theta =\theta\circ \rho(z^{-1})\quad\mbox{for all } z\in U(1)\ .
\end{equation}
A subspace $V\subset W$ is called $\rho$-good (or just good if $\rho$ is understood) if $\dim V=n$ and
$$ \rho(z)V\cap W_+ =\{0\}\quad\mbox{for all } z\in U(1)\ .$$
\end{de}

Every representation of the form (\ref{b.good}) is equivalent to a direct sum of certain two-dimensional representations 
$(\rho_k, W,\theta)$, $k\in\ZZ$, where
\begin{itemize}
\item $W:=\CC$ and $\theta$ is given by complex conjugation and
\item $\rho_k(z)(w):= z^k w$ for $z\in U(1)\subset \CC^*$ and $w\in \CC$.
\end{itemize}
We remark that $\rho_k$ and $\rho_{-k}$ are equivalent as representations over the reals.

We will describe the representations appearing in Definition~\ref{johnny} by $n$-tuples 
$\uk=(k_1,\dots,k_n)\in\ZZ^n$,
where $\uk$ stands for the representation $\rho=\rho_{k_1}\oplus\dots\oplus\rho_{k_n}$ on $\CC^n$, and we will call the
corresponding $\rho$-good subspaces $\uk$-good. Formally, for $n=0$, we will describe the zero-representation $\rho$ by the empty tuple.

We will use the following description of real $n$-dimensional subspaces of $W\cong\CC^n$. Let $c_1,\dots,c_n\in \CC^n$ be linear independent over $\RR$ and let $C=(c_{jk})_{j,k=1,\dots,n}$ be the matrix with columns $c_1,\dots,c_n$. Then 
\begin{equation}\label{EV}
V_C:=\{Cr\mid r\in\RR^n\}
\end{equation}
is the $n$-dimensional real subspace of $W$ spanned by $c_1,\dots,c_n$. 

For $\uk=(k_1,\dots,k_n)$ and $z\in\CC$, let $z^\uk$ denote the diagonal matrix $\diag (z^{k_1}, \dots, z^{k_n})$.
\begin{lm}\label{Ldet}
The subspace $V:=V_C\subset W$ is $\uk$-good if and only if 
$$
\det \Im(z^\uk\, C)\not=0
$$
for all $z\in U(1)\subset\CC$. 
\end{lm}
\proof Obviously, $\rho(z)V\cap W_+\not=\emptyset$ holds if and only if there exists an element $r\in\RR^n$ such that $\rho(z)Cr\in\RR^n$. Furthermore, we have $\rho(z)Cr\in\RR^n$ if and only if 
$$0=\Im(\rho(z)Cr)= \Im(z^\uk\, C)\cdot r.$$
This equation has a non-trivial solution $r$ if and only if $\det \Im(z^\uk\, C)=0$. \qed

Let $G_n(W)$ be the Grassmannian of real $n$-dimensional subspaces of $W$ and, for $\uk \in\ZZ^n$, 
let $G_n^{\;\uk}\subset G_n(W)$
be the subset of $\uk$-good subspaces.

\begin{co}\label{open}
The subset $G_n^{\;\uk}(W)\subset G_n(W)$ is open.
\end{co}

\proof The subset $U:=\{V\in G_n(W)\mid V\cap W_+=\{0\}\}\subset G_n(W)$ is open and $G_n^{\;\uk}(W)\subset U$. There is a
continuous map $F: U\rightarrow M(n,\CC)$ characterised by $V_{F(H)}=H$ (see (\ref{EV})) and $\Im(F(H))=\id$ for all $H\in U$. 
For $C\in M(n,\CC)$ we consider the continuous function $f_C: S^1\rightarrow \RR$, 
$f_C(z):=\det \Im(z^\uk\, C)$. Then the map $G$ sending $H$ to $f_{F(H)}$ from
$U$ to the Banach space $C(S^1,\RR)$ is continuous. Functions without zeroes form an open subset $E\subset C(S^1,\RR)$. 
By Lemma \ref{Ldet} we have  $G_n^{\;\uk}(W)=G^{-1}(E)$.
\qed

The question is whether $G_n^{\;\uk}(W)$ is non-empty.
We are mainly interested in the case $W^\rho=\{0\}$, i.e. $k_i\ne 0$ for
$i=1,\dots,n$.

\begin{de}
An $n$-tuple $\uk=(k_1,\dots,k_n)\in(\ZZ_{\not=0})^n$ is called $\RR$-admissible if $G_n^{\;\uk}(\CC^{n})$ is non-empty.
\end{de}

For an $n$-tuple, the condition of $\RR$-admissibility is invariant under permutations and independent sign changes of the coordinates
as well as under multiplication with a common factor $m\in \ZZ_{\not=0}$. 

\begin{pr}\label{charles}
If $\uk=(k_1,\dots,k_n)\in(\ZZ_{\not=0})^n$ is $\RR$-admissible, then
\begin{equation}\label{Egood}
\sum_{j=1}^n \kappa_j k_j =0
\end{equation}
for suitable $\kappa_j\in\{1,-1\}$, $j=1,\dots,n$.
\end{pr}
\proof 
Since $k_1,\dots,k_n$ are integers, the determinant
$$\det \Im(z^\uk\, C)=\det \left(\textstyle{\frac1{2i}}(c_{lm} z^{k_l}-\overline{c_{lm}} z^{-k_l})_{l,m=1,\dots,n}\right)$$
is a rational function of the form $f_C(z)=\sum_\kappa d_\kappa z^{\kappa_ 1k_1+\dots +\kappa_n k_n}$, $\overline{d_\kappa}=d_{-\kappa}$, where the summation runs over all $\kappa=(\kappa_1,\dots,\kappa_n)\in\{1,-1\}^n$. We consider $f_C$ as a function $f_C:S^1\rightarrow \RR$. Since the integral of the function $t\mapsto e^{itk}$ over $[0,2\pi]$ vanishes for $k\in\ZZ_{\not=0}$, we obtain
\begin{equation} \label{Ei}
\int_0^{2\pi}f_C(e^{it})dt=2\pi\sum_{\kappa_1 k_1+\dots +\kappa_n k_n=0} \ d_\kappa.
\end{equation}
If there exists a $\uk$-good subspace $V=V_C$, then $f_C(z)\not=0$ for all $z\in S^1$ by Lemma~\ref{Ldet}. But if $f_C<0$ or $f_C>0$ on $S^1$,
then the sum on the right hand side of (\ref{Ei}) cannot vanish, which gives the assertion.
\qed 

The following theorem reduces the classification of Cahen-Wallach spaces of imaginary type admitting a compact quotient to that of $\RR$-admissible tuples.

\begin{theo}\label{CC}
Let $X$ be an $(n+2)$-dimensional Cahen-Wallach space of imaginary type. Then $X$ admits a compact quotient if and only if there exists an $\RR$-admissible $d$-tuple $(k_1,\dots,k_d)\in(\ZZ_{\not=0})^d$, $0\le d\le n$, $d \equiv n\ (2)$, such
that
$$ X\cong X_{0,n}(k_1,\dots, k_d, \mu_{d+1},\dots,\mu_n)\,,
$$
where the remaining parameters $\mu_i\in\RR^*$, $i=d+1,\dots,n$, all appear with even multiplicity.
\end{theo}

\proof {\sl Construction part: Suppose that $X=X_{0,n}(k_1,\dots, k_d, \mu_{d+1},\dots,\mu_n)$, where $(k_1,\dots,k_d)$ is $\RR$-admissible 
and $\mu_{d+2j-1}=\mu_{d+2j}=:\hat\mu_j$ for $j=1,\dots,(n-d)/2$.} Using the identification of $\fa$  with $\CC^n$ (see Example~\ref{Exosc2}) 
we make the splitting
$$\fa=\CC^d\oplus \CC^{n-d}=:\fa_1\oplus\fa_2\ .$$
Let $V_1\subset \CC^d$ be a $(k_1,\dots,k_d)$-good subspace. Moreover, we define $\phi:=\phi_{\hat\mu}\in\so(\fa_2)$ for $\hat\mu=(\hat \mu_1,\dots,\hat\mu_{(n-d)/2} )$ as in (\ref{phi}). Then  $\phi$ commutes with $L|_{\fa_2}$ and with complex conjugation. We put $V_2:=\fa_2^{L+\phi}$~and 
$$V:=V_1\oplus V_2,\  t_0:=2\pi,\ \ph_0:={\id}_{\fa_1}\oplus e^{2\pi\phi}\in K\,.$$
These data satisfy Conditions $(a)$ and $(b')$ in Proposition~\ref{wulf}. Hence $X$ admits a compact quotient, which finishes this part of the proof. For the reader who is interested in the explicit construction of lattices, we remark that the proof of Proposition~\ref{wulf} shows that such a compact quotient $\Gamma\setminus X$ can be obtained by taking $\Gamma=\langle \Lambda, \gamma_0\rangle$, where $\gamma_0:=(0,0,t_0)\cdot \ph_0\in \hat G\rtimes K$ and $\Lambda$ is an arbitrary lattice of the group $\fz\oplus V$.

{\sl Classification part:} Assume that $X=X_{0,n}(\tilde\mu_1,\dots, \tilde\mu_n)$ has a compact quotient. Let $V,t_0,\ph_0$ be as in Prop.~\ref{wulf}, Conditions~$(a)$ and $(b')$. We order the parameters $\tilde\mu_i\in\RR\setminus\{0\}$ such that
$$  \tilde\mu_1,\dots,\tilde\mu_d\in \frac{\pi}{t_0}\ZZ \quad\mbox{and}\quad \tilde\mu_{d+1},\dots,\tilde\mu_n\not\in \frac{\pi}{t_0}\ZZ $$
for some $0\le d\le n$. As above, this induces a splitting 
$$\fa=\CC^d\oplus \CC^{n-d}=:\fa_1\oplus\fa_2\ .$$
For $\mu\in\RR^*$ let $\fa(\mu)$ be the eigenspace of $L^2$ with eigenvalue $-\mu^2$. Assume that the parameter $\tilde\mu_i$ appears with odd multiplicity. We claim that then
$i\le d.$
Indeed, $\fa(\tilde\mu_i)_-$ is odd dimensional, hence $\fa(\tilde\mu_i)_-^{\ph_0^2}\ne\{0\}$. We consider the natural projection $P:V\rightarrow\fa(\tilde\mu_i)^{\ph_0^2}$. It is non-trivial and equivariant under $e^{t_0L}\ph_0$. Choose $v\in V$ with $P(v)\ne 0$. Then
$$ e^{2it_0\tilde\mu_i}P(v)=(e^{t_0L}\ph_0)^2 P(v)=P((e^{t_0L}\ph_0)^2 v)=P(v)\ .$$
Thus $e^{2it_0\tilde\mu_i}=1$, and the claim follows. The claim implies that the parameters $\tilde\mu_i$ for $i>d$ appear with even multiplicity.

Let now $i>d$. We consider the natural $e^{t_0L}\ph_0$-equivariant projection $Q:V\rightarrow \fa(\tilde\mu_i)$ and claim that
\begin{equation}\label{irr}
Q(V)\cap\fa(\tilde\mu_i)_+ =\{0\}\ .
\end{equation}
Indeed, for $v\in V$ we have 
$$Q(v)=Q(e^{t_0L}\ph_0 v)= e^{it_0\tilde\mu_i}\ph_0 Q(v)\quad \mbox{ and }\quad\theta_\fa Q(v)=e^{-it_0\tilde\mu_i}\ph_0 \theta_\fa Q(v)\ .$$
Thus $Q(v)=\theta_\fa Q(v)$ implies $e^{2it_0\tilde\mu_i}\ph_0 Q(v)=\ph_0 Q(v)$. Since $i>d$ we conclude that $Q(v)=0$.

Let $R:V\rightarrow \fa_2$ be the natural projection. Then (\ref{irr}) implies that $\dim R(V)\le{n-d}$. Hence $\dim (V\cap\fa_1)\ge{d}$.
We conclude that $\dim (V\cap\fa_1)={d}$  and that $V\cap\fa_1$ is a $(k_1,\dots, k_d)$-good subspace of $\fa_1$, where $k_i=t_0\tilde\mu_i/\pi$,
$i=1,\dots, d$. Thus the $d$-tuple $(k_1,\dots,k_d)$ is $\RR$-admissible. For $i>d$ we set $\mu_i:=t_0 \tilde\mu_i/\pi$. Then
$$ X\cong X_{0,n}(k_1,\dots, k_d, \mu_{d+1},\dots,\mu_n)\ ,$$
and the parameters have the required properties.
\qed

\begin{co}\label{spurnull2}
Assume that $X_{0,n}(\mu)$ admits a compact quotient. Then there is choice of signs such that
$$
\displaystyle \sum_{i=1}^{n} \pm\mu_i = 0\ .
$$
\end{co}

\proof Combine Thm.~\ref{CC} with Prop.~\ref{charles}.
\qed

Now we want to construct examples of Cahen-Wallach spaces of imaginary type that admit compact quotients. Theorem \ref{CC} reduces this task to the construction of $\uk$-good subspaces of the $2n$-dimensional real vector space $W =\CC^{n}$, where $\uk=(k_1,\dots,k_n)\in(\ZZ_{\not=0})^n$. Recall that we consider the involution $\theta$ on $W$ given by the complex conjugation on $\CC^n$, hence $W_+=\RR^n\subset\CC^n$, and that $\rho=\rho_{k_1}\oplus\dots\oplus\rho_{k_n}$ is the $U(1)$-representation on $\CC^n$ introduced in the beginning of this section. 

In the remaining part of this subsection we will construct explicit examples of $\uk$-good subspaces, mainly in small dimensions. For $n=2$ this is rather easy. The subspace $V=\CC\cdot (1,i)\subset \CC^2$ is $(k,k)$-good for every $k\in\ZZ_{\not=0}$.  In order to make the calculations for higher dimensions more readable we introduce the following notations:
$$I(k):=\Im(z^k),\ I(k,\omega):=\Im(z^k\omega),\ R(k):=\Re(z^k),\ R(k,\omega)=\Re(z^k\omega), $$
where $z\in S^1$, $\omega\in\CC$ and $k\in\ZZ$.
\begin{ex}[n=3]\label{Exdim3}{\rm Let $\uk=(k_1,k_2,k_3)\in(\ZZ_{\not=0})^3$ be such that $k_1=k_2+k_3$. We consider $$C=\mbox{\begin{small}$\left(\begin{array}{ccc} i&-1&0\\ 1&i&-\omega\\ i\omega&-\omega& i\end{array}\right)$ \end{small}}$$
and we claim that  $V_C$ is a good subspace for a suitable choice of $\omega$. Indeed, we compute
$$f_C(z):=\det \Im(z^\uk\, C)=\mbox{\begin{small}$\left|\begin{array}{ccc} R(k_1)&-I(k_1)&0\\ I(k_2)&R(k_2)&-I(k_2,\omega)\\ R(k_3,\omega)&-I(k_3,\omega)& R(k_3)\end{array}\right| $\end{small}}$$
expanding the determinant along the first row. We obtain
\begin{eqnarray}\label{Edev} f_C(z)&=& R(k_1)\cdot \mbox{\begin{small}$\left|\begin{array}{cc} R(k_2)&-I(k_2,\omega)\\ -I(k_3,\omega)& R(k_3)\end{array}\right|$\end{small}}
+I(k_1)\cdot\mbox{\begin{small}$\left|\begin{array}{cc} I(k_2)&-I(k_2,\omega)\\ R(k_3,\omega)& R(k_3)\end{array}\right|$\end{small}}\nonumber\\[1ex]
&=&  \mbox{\begin{small}$\left|\begin{array}{cc} R(k_1)R(k_2)+I(k_1)I(k_2)&-I(k_2,\omega)\\ -R(k_1)I(k_3,\omega)+I(k_1)R(k_3,\omega)& R(k_3)\end{array}\right|$\end{small}}\\[1ex]
&=& \mbox{\begin{small}$\left|\begin{array}{cc} R(k_2-k_1)&-I(k_2,\omega)\\ I(k_1-k_3,\omega)& R(k_3)\end{array}\right|$\end{small}}\nonumber\\[1ex]
&=&R(k_3)^2+I(k_2,\omega)^2. \nonumber
\end{eqnarray}
If we take $\omega=e^{ir}$ such that $r$ is an irrational multiple of $\pi$, then $R(k_3)$ and $I(k_2,\omega)$ do not vanish at the same time. Consequently, $V_C$ is a $\uk$-good subspace for this choice of~$\omega$.
If we take $\omega=1$,
then $V_C$ is always $(k_1,k_2,k_3)$-good or $(k_1,k_3,k_2)$-good. Indeed, if the multiplicity of 2 in the prime factorisation of $k$ is less or equal to that of $2$ in the prime factorisation of $k'$, then $R(k')^2+I(k)^2\ne 0$.
}\end{ex}
\begin{ex}[n=4]\label{Exdim4}{\rm Let $\uk=(k_1,k_2,k_3,k_4)\in (\ZZ_{\not=0})^4$ be a quadruple satisfying $k_2-k_1=k_4-k_3$. We consider  $V_C$ for 
$$C= \mbox{\begin{small} $\left(\begin{array}{cccc} i&-1&0&0\\ 1&i&-\omega&-i\omega\\ i\omega&-\omega& i&-1\\
0&0&1&i\end{array}\right)$\end{small}}.$$
We claim that $V_C\subset \CC^4$ is a good subspace for a suitable choice of $\omega$. Indeed, expanding the determinant $f_C(z):=\det \Im(z^\uk\, C)$ along the first and the last row proceeding in the same way as in (\ref{Edev}) we obtain 
\begin{eqnarray*}f_C(z)&=&\mbox{\begin{small}$\left|\begin{array}{cccc} R(k_1)&-I(k_1)&0&0\\ I(k_2)&R(k_2)&-I(k_2,\omega)&-R(k_2,\omega)\\ R(k_3,\omega)&-I(k_3,\omega)& R(k_3)& -I(k_3)\\ 0&0&I(k_4)&R(k_4)\end{array}\right|$\end{small}}\\[1ex]
&=&  \mbox{\begin{small} $\left|\begin{array}{cc} R(k_1)R(k_2)+I(k_1)I(k_2)&-I(k_2,\omega)R(k_4)+R(k_2,\omega)I(k_4)\\- R(k_1)I(k_3,\omega)+I(k_1)R(k_3,\omega)& R(k_3)R(k_4)+I(k_3)I(k_4)\end{array}\right|$\end{small}}\\[1ex]
&=& \mbox{\begin{small} $\left|\begin{array}{cc} R(k_2-k_1)&-I(k_2-k_4,\omega)\\ I(k_1-k_3,\omega)& R(k_4-k_3)\end{array}\right|$\end{small}}\\[1ex]
&=&R(k_2-k_1)^2+I(k_1-k_3,\omega)^2.
\end{eqnarray*}
As in Example \ref{Exdim3} we conclude that $V_C$ is $\uk$-good if $\omega$ is in $\pi\cdot (\RR\setminus \QQ)$.  If we take $\omega=1$, 
then $V_C$ is $(k_1,k_2,k_3,k_4)$-good or $(k_1,k_3,k_2,k_4)$-good by the same argument as in the previous example.

By the way, note that $V_C$ is a complex subspace. This will become of interest in Subsection~\ref{trans}. 
}\end{ex}
\begin{ex}[n=6]\label{Exdim6}{\rm Suppose that $\uk=(k_1,\dots,k_6)\in(\ZZ_{\not=0})^6$ satisfies $k_1+k_5+k_6=k_2+k_3+k_4$. We consider $C:=(c_1,i c_1, c_2, i c_2, c_3, i c_3)$, where
$$c_1:=(1,0,0,2,-i,1)^\top,\ c_2:=(0,1,0-i,-1,-i)^\top,\ c_3:=(0,0,1,0,-1,i)^\top.$$
We want to find examples of 6-tuples $\uk$ for which $V_C$ is $\uk$-good. Expanding the determinant $f_C(z):=\det \Im(z^\uk\, C)$ along the first three rows we obtain
$$
f_C(z)= \left| \mbox{\begin{small} $\begin{array}{rrc} 2I(k_1-k_4) & -R(k_2-k_4) &0\\R(k_1-k_5) & I(k_2-k_5) &-I(k_3-k_5)\\ I(k_1-k_6) &-R(k_2-k_6) &-R(k_3-k_6)
\end{array}$\end{small}}\right|.
$$
With
$$\alpha:= k_1-k_2,\ \beta:=k_1-k_4, \ \gamma:= k_5-k_3,\ \delta:= k_6-k_3$$
this gives
\begin{eqnarray*}
f_C(z)&=& \left| \mbox{\begin{small} $\begin{array}{crc} 2I(\beta) & -R(\beta-\alpha) &0\\R(\alpha+\beta+\delta) & I(\beta+\delta) &I(\gamma)\\ I(\alpha+\beta+\gamma) &-R(\beta+\gamma) &-R(\delta)
\end{array}$\end{small}}\right|\\[1ex]
&=& -I(\gamma)(-2I(\beta)R(\beta+\gamma)+I(\alpha+\beta+\gamma)R(\beta-\alpha))\\
&&-R(\delta)(2I(\beta)I(\beta+\delta)+R(\alpha+\beta+\delta)R(\beta-\alpha))\\[1ex]
&=& -I(\gamma)(I(\gamma)-\textstyle{\frac12}I(2\beta+\gamma) +\textstyle{\frac12}I(2\alpha+\gamma))\\
&&-R(\delta)(R(\delta)-\textstyle{\frac12}R(2\beta+\delta)+\textstyle{\frac12}R(2\alpha+\delta))\\[1ex]
&=&-\textstyle{\frac14}(4+2R(2\alpha)-2R(2\beta)-2R(2\gamma)+2R(2\delta)+R(2\alpha+2\delta)\\
&&-R(2\alpha+2\gamma)+R(2\beta+2\gamma)-R(2\beta+2\delta))\\[1ex]
&=&-\textstyle{\frac14} \Re( (2+z^{2\alpha}-z^{2\beta})(2-z^{2\gamma}+z^{2\delta})).
\end{eqnarray*}
If $k_2=k_4$ (which is equivalent to $\alpha=\beta$) and both $\gamma= k_5-k_3$ and  $\delta= k_6-k_3$ are odd, then  $f_C(z)>0$ for all $z\in \U(1)$, hence $V_C$ is $\uk$-good for those $\uk$.

Now suppose that $k_1-k_2= k_5-k_3$ and $k_1-k_4=k_6-k_3$, which is equivalent to $\alpha=\gamma$ and $\beta=\delta$. Then we have  $4f_C(z)=-4 + \Re((z^{2\alpha}-z^{2\beta})^2)$. If we now choose $\uk$ such that both $\alpha= k_1-k_2$ and $\beta=k_1-k_4$ are odd, then  $f_C(z)>0$ for all $z\in \U(1)$, hence $V_C$ is $\uk$-good. 

Besides these two series of examples there are many special choices of $\alpha,\beta,\gamma$ and $\delta$ such that $f_C(z)>0$ for all $z\in \U(1)$. Each of these choices gives us an infinite series of admissible 6-tuples $\uk$, namely
$$k_1=k+\alpha+\beta,\ k_2=k+\beta,\ k_3=k-\gamma-\delta,\ k_4=k+\alpha,\ k_5=k-\delta,\ k_6=k-\gamma$$
for any $k\in\ZZ\setminus\{-\alpha,\,-\beta,\,\gamma,\,\delta,\,\-\alpha-\beta,\,\gamma+\delta\}$. 
For instance, one can check numerically that $f_C(z)>0$ for $(\alpha,\beta,\gamma,\delta)=(1,5,3,12)$. If we put $k=-16$ we get $\uk=-(10,11,31,15,28,19)$. We will get back to this example a the end of this subsection.}
\end{ex}
\begin{ex}[Examples by induction]\label{Bind}{\rm Now we want to construct $\tilde \uk$-good subspaces of $\CC^{n+2}$ starting from $\uk$-good subspaces in $\CC^n$, where $\uk$ and $\tilde \uk$ are related in the following way. Let $\uk=(k_1,\dots,k_n)$ be $\RR$-admissible. We take $k,\hat k$ such that $k_n=2k-\hat k$ and put $\tilde \uk=(k_1,\dots,k_{n-1},k,k,\hat k)$. Let $V_C$ be $\uk$-good. We may assume that the last row of $C$ equals $(0,\dots,0,1,i)$. We claim that  $V_{\tilde C}$ is $\tilde \uk$-good for
$$\tilde C:=\left( \begin{array}{c|c}
(c_{lm})_{l=1,\dots,n-1,\, m=1,\dots,n}&0\\[0.5ex]
\hline
\begin{array}{c|c} 0\ \ & \mbox{\begin{small} $\begin{array}{cc} 1 &i\\i& -1\\ 0&0\end{array}$\end{small}}
\end{array}&\mbox{\begin{small} $\begin{array}{cc} -1&- i\\i&-1\\1&i\end{array}$\end{small}}
\end{array}\right).$$
Indeed, for $f_{\tilde C}(z):=\det\Im(z^{\tilde \uk}\, \tilde C)$, we have
\begin{eqnarray*} f_{\tilde C}(z)&=& \det \left( \begin{array}{c|c}
\Im(z^{k_l}c_{lm})_{l=1,\dots,n-1,\, m=1,\dots,n}&0\\[0.5ex]
\hline
\begin{array}{c|c} 0\ \ & \mbox{\begin{small} $\begin{array}{cc} I(k) &R(k)\\ R(k)& -I(k)\\ 0&0\end{array}$\end{small}}
\end{array}&\mbox{\begin{small} $\begin{array}{cc} -I(k)& -R(k)\\R(k)&-I(k)\\I(\hat k)&R(\hat k)\end{array}$\end{small}}
\end{array}\right) \\[1ex]
&=&\det \left( \begin{array}{c|c}
\Im(z^{k_l}c_{lm})_{l=1,\dots,n-1,\, m=1,\dots,n}&0\\[0.5ex]
\hline
\begin{array}{c|c} 0\ \ & \mbox{\begin{small} $\begin{array}{cc} I(k) &R(k)\\ R(k)& -I(k)\end{array}$\end{small}}
\end{array}&\mbox{\begin{small} $\begin{array}{c} -I(k-\hat k)\\R(k-\hat k)\end{array}$\end{small}}
\end{array}\right) \\[1ex]
&=&\det \left( \begin{array}{c}
\Im(z^{k_l}c_{lm})_{l=1,\dots,n-1,\, m=1,\dots,n}\\[0.5ex]
\hline
\begin{array}{c|c} 0\ \ & \mbox{\begin{small} $\begin{array}{cc} I(2k-\hat k) &R(2k-\hat k)\end{array}$\end{small}}
\end{array}
\end{array}\right)\\[1ex]
&=&\det\Im(z^\uk\, C)\ \not=\ 0
\end{eqnarray*}
since $V_C$ is $\uk$-good by assumption. Hence $V_{\tilde C}$ is $\tilde \uk$-good. Note that $V_{\tilde C}$ is complex if $V_{C}$ is complex.}\end{ex}

The induction gives rise to many examples of $\uk$-good subspaces in arbitrary dimensions. For instance, we can show that $\uk$ is $\RR$-admissible if at most four entries of $\uk$ have odd multiplicity and if $\uk$ satisfies the necessary condition (\ref{Egood}).

\begin{co} \label{co} For $l\le 4$, let $\uk=(k_1,\dots,k_{2m},k_{2m+1},\dots, k_{2m+l}) \in(\ZZ_{\not=0})^{2m+l}$ be such that $k_{2j-1}=k_{2j}$ for $j=1,\dots,m$. Then $\uk$ is $\RR$-admissible if and only if it satisfies {\rm (\ref{Egood})} for a suitable choice of signs.
\end{co}
\proof We have to show that Condition (\ref{Egood}) is sufficient. We may assume that $-2\sum_{j=1}^m(-1)^jk_j+(-1)^mk_{2m+1}+k_{2m+2}+\dots+k_{2m+l}=0$. Examples~\ref{Exdim3} and~\ref{Exdim4} show that $(k'_{2m+1},k_{2m+2},\dots,k_{2m+l})$ is $\RR$-admissible if $k'_{2m+1}+k_{2m+2}+\dots+k_{2m+l}=0$. 
Now the assertion follows by repeated application of the induction step described in Example~\ref{Bind}. \qed

In particular, Cor.~\ref{co} shows that for $n\le 4$ the trace condition (\ref{Egood}) is not only necessary but also sufficient for the existence of a $\uk$-good subspace. It is not clear whether this is also true for higher dimensions. Unfortunately, the described induction only yields examples of $\RR$-admissible $\uk$ whose entries $k_j$ of $\uk$ satisfy various linear equations of the form $\sum_{j\in J} \kappa_j k_j=0$, where $\kappa_j\in\{1,-1\}$ and $J$ is a subset of $\{1,\dots,n\}$. Moreover, $k_i=k_j$ for at least one pair $(i,j)$, $i\not=j$. This leads us to the following

{\bf Question:} Does the existence of a $\uk$-good subspace imply any further condition for $\uk$ besides Condition (\ref{Egood})?

We have some hope that the answer to this question is `no' also in dimension $n>4$. For instance, $\uk=(10,11,31,15,28,19)$ is an example  of an $\RR$-admissible 6-tuple (see Example~\ref{Exdim6}), that does not satisfy any linear equation of the kind mentioned above except the trace condition.  

\subsection{Special subspaces}\label{SSS} 
In this subsection we will study a special class of good subspaces. We want to do this in a more general context, which will become of importance in Sections~\ref{hans-werner} and~\ref{last}. 

Let $(W,\ip)$ be a Euclidean space of dimension $2n$ and let $\theta:W\rightarrow W$ be an involutive isometry. As usual, we denote the eigenspaces of $\theta$ by $W_+$ and $W_-$. We suppose that $\dim W_+=\dim W_-$. Furthermore, let be given a map $L\in\so(W)$  that anticommutes with $\theta$. 
\begin{de}\label{henze} Let $\phi\in\so(W)$ be a map that commutes with $L$ and $\theta$.
A subspace $V\subset W$ is called $(L,\phi)$-special  if  $W=V\oplus W_+$ and if $(L+\phi)(V)=V$. 

A subspace $V\subset W$ is called $L$-special if there exists a map $\phi\in\so(W)$ commuting with $\theta$ and $L$ such that $V$ is $(L,\phi)$-special.
\end{de} 
Every $L$-special subspace $V\subset W$ satisfies $e^{tL}(V)\cap W_+=0$ for all $t\in\RR$. Indeed, choose $\phi$ such that $V$ is $(L,\phi)$-special, then $e^{tL}e^{t\phi}(V)=V$ for all $t\in\RR$. Hence $e^{tL}(V)\cap W_+=e^{-t\phi}(V)\cap W_+=e^{-t\phi}(V \cap e^{t\phi}(W_+))=e^{-t\phi}(V \cap W_+)=0$ for all $t\in\RR$ since $e^{-t\phi}$ commutes with $\theta$. 

In particular, suppose that the eigenvalues of $L$ on $W_{\Bbb C}$ are in $i\ZZ$, and let $\pm ik_1,\dots,\pm ik_n$ be these eigenvalues. Then every $L$-special subspace is $(k_1,\dots,k_n)$-good. 

In the remainder of this section we will suppose that $L\in\so(W)$ is invertible. Then $\dim W_+=\dim W_-$ holds automatically.

Our aim is to give a criterion for the existence of $L$-special subspaces in terms of the eigenvalues of $L$ on $W_{\Bbb C}$. Let $\pm i\mu  _1,\dots,\pm i\mu_n$ be these eigenvalues. We may assume that the absolute values of $\mu_k=:\alpha_k$, $k=1,\dots,p$, are pairwise distinct and that $\mu_{p+2j-1}=\mu_{p+2j}=:\beta_j$, $j=1,\dots,q$. 

Let us consider the following examples before we will formulate a general criterion.
\begin{ex}\label{B0}{\rm
Take $W=\RR^4\cong\CC^2$ with standard scalar product, $\theta$ the complex conjugation and let $L\in\so(4)$ be defined by $L(z_1,z_2)=i(\beta z_1,\beta z_2)$. Then 
$V=\CC\cdot (1,i)$ considered as a real vector space is a complement of $W_+=\RR^2$. Furthermore, $V$ is $L$-invariant, thus $L$-special.}\end{ex}

\begin{ex}\label{B1}{\rm Consider $W=\CC^2\oplus\CC^{2q}$ as a real vector space endowed with its standard scalar product. Let $\theta$ be the complex conjugation. According to the decomposition $W=\CC^2\oplus\CC^{2q}$ we define 
$$L:=L_\alpha\oplus L_\beta,\quad \alpha=(\alpha_1,\alpha_2),\ \beta=(\beta_1,\beta_1,\dots,\beta_q, \beta_q)$$
using (\ref{L}). Suppose that 
\begin{equation} \label{Eab} 
\alpha_1-(-1)^q \alpha_2=2\sum_{j=1}^q (-1)^{j+1}\beta_j.
\end{equation} 
We are going to show that there exists an $L$-special subspace $V\subset W$. 

We define inductively $\gamma_1,\dots, \gamma_q\in \RR$ by 
\begin{equation}\label{Eabc}
\gamma_q:=\beta_q-\alpha_2,\quad \gamma_j=-\gamma_{j+1}+\beta_j-\beta_{j+1}, \quad j=q-1,\dots, 1.
\end{equation}
Then 
\begin{equation}\label{Eag}
\alpha_1=\beta_1+\gamma_1
\end{equation} 
by (\ref{Eab}). According to the decomposition $W=\CC^2\oplus\CC^{2q}$ we define $\phi\in\so(W)$ by
$$\phi=0_2\oplus \phi_\gamma,\quad \gamma=(\gamma_1,\dots,\gamma_q)$$
using (\ref{phi}), where $0_2$ denotes the zero map.
Let us denote the (complex) standard basis of $\CC^{2+q}$ by $e'_1,e'_2,e_1,\dots,e_q$. Then $L+\phi$ maps
\begin{eqnarray*}
b_0:=e'_1+ie_1+e_2&\longmapsto & i\alpha_1e'_1+i(\beta_1+\gamma_1)(ie_1+e_2)\\
b_q:=e'_2+e_{2q-1}+ie_{2q}&\longmapsto& i\alpha_2e'_2+i(\beta_q-\gamma_q)(ie_{2q-1}+ie_{2q})\\
b_j:=e_{2j-1}+ie_{2j}+e_{2j+1}-ie_{2j+2}&\longmapsto& i(\beta_j-\gamma_j)(e_{2j-1}+ie_{2j})\\
&&+i(\beta_{j+1}+\gamma_{j+1})(e_{2j+1}-ie_{2j+2})
\end{eqnarray*} for $j=1,\dots,q-1$. Equations (\ref{Eabc}) and (\ref{Eag}) show that $b_0,\dots,b_q$ are eigenvectors of $L+\phi$ considered as a complex linear map on $\CC^2\oplus\CC^{2q}$. In particular, $V:=\Span\{b_0,\dots,b_q\}$ is a $(q+1)$-dimensional complex subspace of $W$, which is invariant under $L+\phi$. We claim that $V\subset W$ considered as a real vector space is $L$-special. To verify this it remains to show that the projection of $V$ to $W_-$ is surjective. But this is true since this projection is spanned by the projections of $b_0,\dots,b_q$ and $ib_0,\dots, ib_q$, which are 
\begin{eqnarray*}
&ie_1,\quad      ie_{2j}-ie_{2j+2} \ (j=1,\dots,q-1),      \quad  ie_{2q},&\\
&ie'_1+ie_2,\quad    ie_{2j-1}+ie_{2j+1} \ (j=1,\dots,q-1),\quad     ie'_2+ie_{2q-1}.&
\end{eqnarray*}
}\end{ex}

\begin{ex}\label{B2}{\rm Similar to the previous example we consider the real vector space $W=\CC\oplus\CC^{2q}$ endowed with the Euclidean standard scalar product and complex conjugation~$\theta$. Let $L$ be defined by
$$L=L_\alpha \oplus L_\beta,\quad \alpha\in\RR, \ \beta= (\beta_1,\beta_1,\dots,\beta_q, \beta_q)\,.$$
Suppose now that 
$$\alpha=2\sum_{j=1}^q (-1)^{j+1}\beta_j.$$

We define $b_1,\dots,b_{q-1}$ as in Example~\ref{B1}. Moreover, we set $b_0:=e'+ie_1+e_2$, where $e':=1\in\CC$. 
Then the real vector space $V:=\Span_{\Bbb C}\{b_0,\dots,b_q\}\oplus\,\RR (e_{2q-1}+i e_{2q})$ is an $L$-special subspace of $W$. Indeed, $\phi$ can be chosen similar to the map $\phi$ in Example~\ref{B1}, where $\gamma_1,\dots,\gamma_q$ are again defined as in (\ref{Eabc}) but now with $\alpha_2=0$.
}\end{ex}
\begin{pr}\label{special} 
There is an $L$-special subspace $V\subset W$ if and only if there exist
\benur 
\item pairwise disjoint subsets $I_1,\dots, I_{p_0}$ and $I_{2p_0+1},\dots,I_p$ of $\{1,\dots,q\}$; 
\item numbers $\delta_k, c_i\in \{1,-1\}$;
\item and a permutation $\sigma$ of $\{1,\dots,p\}$
\end{enumerate}
such that 
\begin{eqnarray}
\alpha_{\sigma(2k-1)}+\delta_k \alpha_{\sigma(2k)} &=& 2\sum _{i\in I_k} c_i\beta_i,\quad k=1,\dots,p_0,\label{E*} \\
\alpha_{\sigma(l)}&=& 2\sum _{i\in I_l} c_i\beta_i,\ \quad l=2p_0+1,\dots,p. \label{E**}
\end{eqnarray}
\end{pr}
\proof Suppose that Equations (\ref{E*}) and (\ref{E**}) are satisfied for suitable $\delta_k, c_i\in\{\pm1\}$. Since we may change the signs of $\alpha_1,\dots,\alpha_p$ and $\beta_1,\dots,\beta_q$ there is an orthogonal $\theta$-invariant decomposition $W=W_1\oplus \dots \oplus W_r$ such that $L=L_1\oplus \dots \oplus L_r$ with $L_\nu :W_\nu\rightarrow W_\nu$ and $(W_\nu, L_\nu)$ is isomorphic to one of the Examples \ref{B0} -- \ref{B2} for $\nu=1,\dots, r$. Hence there exists an $L$-special subspace. 

Now suppose that there is an $L$-special subspace $V\subset W$. Then $W=V\oplus W_+$ and there exists a map $\phi\in\so(W)$ commuting with $\theta$ and $L$ such that $(L+\phi)(V)=V$. We can identify $W$ with $\CC^p\oplus \CC^{2q}$ such that $L$ and $\phi$ are given with respect to this decomposition by 
\begin{eqnarray*}
L=L_\alpha\oplus L_\beta, &&\alpha=(\alpha_1 ,\dots,\alpha_p),\ \beta= (\beta_1, \beta_1,\dots, \beta_q,\beta_q),\\
\phi=0_p\oplus \phi_\gamma& & \gamma=(\gamma_1,\dots, \gamma_q) 
\end{eqnarray*}
according to (\ref{L}) and (\ref{phi}), where $0_p$ denotes the zero map.

We denote by $e'_1,\dots, e'_p, e_1,\dots, e_{2q}$ the (complex) standard basis of $\CC^p\oplus \CC^{2q}$. Considered as a complex linear map, $L+\phi$ has eigenvalues $i\alpha_k$ on 
\begin{equation}\label{EW1}
W^1_k:=\CC\cdot e'_k,\quad  k=1,\dots,p, 
\end{equation}
and $i(\beta_j-\gamma_j)$, $i(\beta_j+\gamma_j)$
on 
\begin{equation}\label{EW2}
W^2_j:=\Span{}_{\Bbb C}\{e_{2j-1}, e_{2j}\},\quad  j=1,\dots,q.
\end{equation}
Now we consider $W^1_k$ and $W^2_j$ as real vector spaces. 
Then their complexifications $(W^1_k)_{\Bbb C}$ and $(W^2_j)_{\Bbb C}$ are $(L+\phi)$-invariant subspaces of $W_{\Bbb C}$ on which $L+\phi$ has eigenvalues $\pm i\alpha_k$  and $\pm i(\beta_j+\gamma_j)$, $\pm i(\beta_j-\gamma_j)$, respectively. We put
$$
\rho_j:=\beta_j+\gamma_j,\quad \rho_j':=\beta_j-\gamma_j,\quad j=1,\dots,q.
$$
To $(W,L+\phi)$, we assign a graph $G$. In general, this graph has multiple edges and loops. The set ${\cal V}$ of vertices and the set ${\cal E}$ of edges are 
\begin{equation} \label{EVE}
{\cal V}:= \{ |\rho|\mid i\rho \mbox{ is an eigenvalue of $L+\phi$ on $W_{\Bbb C}$}\}, \quad {\cal E}:=\{1,\dots,q\},
\end{equation}
where $j$ is an edge between $|\rho_j|$ and $|\rho_j'|$. 
Let $G'$ be a connected component of $G$ with set of vertices ${\cal V}'\subset {\cal V}$ and set of edges ${\cal E}'\subset {\cal E}$.  
We define
$$W':= \bigoplus_{|\alpha_k|\in {\cal V}'} W^{1}_k\,\oplus\, \bigoplus_{j\in {\cal E}'}W_j^{2},\quad W'':= (W')^\perp.$$
Then $V':=W'\cap V$ and $V'':=W''\cap V$ are special subspaces in $W'$ and $W''$, respectively. Indeed, $V_{\Bbb C}$ decomposes into eigenspaces of $L+\phi$. Since $G'$ is a connected component of $G$, these eigenspaces are subspaces of either $W'$ or $W''$. Hence $V=(V\cap W')\oplus (V\cap W'')$, which proves the assertion. Thus it suffices to prove the assertion for the special subspace $V'\subset W'$ that corresponds to the connected component $G'$ of $G$.  For simplicity of notation, we will write again $G$ instead of $G'$ and $V\subset W$ instead of $V'\subset W'$ and we now denote by $p$ the number of eigenvalues of odd multiplicity of $L+\phi$ on $W'$.

If $p$ is odd, we put $\alpha_{p+1}:=0$. Note, that in this case $0$ is already a vertex of $G$ since $\dim V$ is odd if $p$ is odd and thus $(L+\phi)|_V$ has a non-trivial kernel. Furthermore, we put 
$\tilde p:= p$ if $p$ is even and $\tilde p:=p+1$ if $p$ is odd.

We need the following property of a connected graph. Given a set ${\cal V}'$ of $2m$ vertices of $G$, we find $m$ edge-disjoint paths such that ${\cal V}'$ equals the set of endpoints of these paths. In order to verify it, we choose a tree $T$ in $G$ that contains all elements of ${\cal V}'$.  We fix a vertex $P$ of $T$ and consider $T$ as an out-tree rooted at $P$. Now we choose a pair $\{v,v'\}\subset {\cal V}'$ such that the lowest common ancestor of $v$ and $v'$ has maximal distance to $P$. Let $l$ denote the path in the undirected tree $T$ joining $v$ and $v'$. If we remove the edges of $l$ from $T$, then all elements of ${\cal V}'':={\cal V}'\setminus \{v,v'\}$ are contained in the same connected component $T'$ of the remaining graph and we can proceed in the same way with $T'$ and ${\cal V}''$, etc. 

We apply this property to ${\cal V}':=\{|\alpha_1|,\dots, |\alpha_{\tilde p}|\}$ and conclude that there exist a permutation $\sigma$ of $\{1,\dots,\tilde p\}$ and pairwise edge-disjoint paths $l_k$, $k=1,\dots, \tilde p/2$, such that $l_k$ joins $|\alpha_{\sigma(2k-1)}|$ and $|\alpha_{\sigma(2k)}|$. If $p$ is odd, we may assume that $\sigma(\tilde p) =\tilde p$.

We fix $1\le k\le  \tilde p/2$ and consider the path $l_k$ consisting of the sequence $(j_1,\dots,j_r)$ of edges. By construction of $G$ we have 
\begin{eqnarray*}
\alpha_{\sigma(2k-1)} &=& c_{j_1}\beta_{j_1} +\eps_{j_1}\gamma_{j_1},\\
c_{j_\nu}\beta_{j_\nu} - \eps_{j_\nu}\gamma_{j_\nu}&=&-(c_{j_{\nu+1}}\beta_{j_{\nu+1}} + \eps_{j_{\nu+1}}\gamma_{j_{\nu+1}}),\quad \nu=1,\dots,r-1,\\
c_{j_{r}}\beta_{j_{r}} -\eps_{j_{r}}\gamma_{j_{r}}&=&\delta_k \alpha_{\sigma(2k)}
\end{eqnarray*}
for suitable $c_{j_1},\dots, c_{j_r},\eps_{j_1},\dots, \eps_{j_r},\delta_k \in \{1,-1\}$.

We put $I_k:=\{j_1,\dots,j_r\}$.  We obtain (\ref{E*}) if $2k\not=p+1$. Otherwise we get (\ref{E**}).\qed

\section{The general case}\label{hans-werner}
\subsection{Preliminaries}

We will need the following refinement of Lemma \ref{gigi}.

 
 \begin{lm}\label{nono}
Let $(V,\omega)$ be a symplectic vector space over $\RR$, and let $A\in \Sp(V,\omega)$ be semisimple such that its characteristic polynomial $f_A$ has integer coefficients.
Then there exists a lattice $\Lambda\subset V$ with $A\Lambda=\Lambda$ and $\omega(\Lambda\times \Lambda)\subset\ZZ$.
\end{lm}

\proof 
Let us first reduce the assertion to the case that $f_A$ is irreducible over $\QQ$. For a polynomial $p$ let $p^*$ be the corresponding reciprocal polynomial. Then $f_A^*=f_A$.
We decompose $f_A=f_1\cdot\dots\cdot f_{2s}\cdot f_{2s+1}\cdot\dots\cdot f_{2s+r}$
into irreducible monic factors such that $f_{2k-1}^*=\kappa_k f_{2k}$, $k=1,\dots,s$, where $\kappa_k=\pm 1$ is the constant term
of $f_{2k-1}$, and $f_l^*=f_l$, $l>2s$. There is a corresponding decomposition of $V$ into $A$-invariant pairwise orthogonal symplectic subspaces
$$ V= V_1\oplus\dots \oplus V_s\oplus V_{s+1}\oplus\dots\oplus V_{s+r} $$
and a further decomposition
$$  V_k= W_k\oplus W_k'\ ,\quad k=1,\dots, s\ ,$$
into $A$-invariant Lagrange subspaces such that
$$ f_{A|_{W_k}}=f_{2k-1},\ f_{A|_{W_k'}}=f_{2k},\  1\le k\le s,\quad f_{A|_{V_{s+l}}}=f_{2s+l},\ 1\le l\le r.$$
According to Lemma~\ref{gigi} we find $A$-stable lattices $\Lambda_k\subset W_k$, $k\le s$. Let 
$$\Lambda_k'=\{ w'\in W_k'\mid \omega(w,w')\in\ZZ\mbox{ for all } w\in \Lambda_k\}\subset W_k'$$
be the corresponding dual lattice. Assuming the lemma for irreducible characteristic polynomials there are $A$-stable lattices
$\Lambda_{s+l}\subset V_{s+l}$, $1\le l\le r$ satisfying $\omega(\Lambda_{s+l}\times \Lambda_{s+l})\subset\ZZ$.
Then
$$ \Lambda=\Lambda_1\oplus\Lambda_1'\dots \oplus \Lambda_s\oplus\Lambda_s'\oplus \Lambda_{s+1}\oplus\dots\oplus \Lambda_{s+r}\subset V$$
has the desired properties.

It remains to prove the lemma for irreducible $f_A$.
For every monic self-reciprocal polynomial $p$ of degree $2n$ over $\ZZ$ one can find a matrix $M\in\Sp(n,\ZZ)$ whose characteristic polynomial equals $p$, cf.~\cite{Ki}, see also \cite{Ri}, Theorem A1.  Hence there exists a matrix $A_0\in\Sp(n,\ZZ)$, where $\dim V=2n$, having the same characteristic polynomial as $A$. Since $f_A$ is irreducible over $\QQ$ all roots of $f_A$ are simple. Hence $A=T ^{-1} A_0T$ for some isomorphism $T:V\to\RR^{2n}$. In particular, if $\omega_0$ denotes the standard symplectic form on $\RR^{2n}$, then $\omega':=T^*\omega_0$ is an $A$-invariant symplectic form on $V$, the lattice $\Lambda':=T^{-1}(\ZZ^{2n})\subset V$ is $A$-stable and $\omega'(\Lambda'\times\Lambda')\subset \ZZ$.

Let $F\subset \CC^*$ be the set of eigenvalues of $A$. It is stable under inversion and complex conjugation.  We split it into a disjoint union $F=F_+\cup F_0\cup F_-$ of eigenvalues having modulus greater than $1$, equal to $1$, and smaller than $1$, respectively. We fix a basis $\{e_\nu\mid \nu\in F\}$ of the complexification $V_\Bbb{C}$ consisting of eigenvectors of $A$ and satisfying $e_{\bar\nu}=\overline{e_\nu}$. Then $A$-invariant symplectic forms are in one-to-one correspondence to functions $c: F\rightarrow \CC^*$ satisfying
$$ c(\nu)=-c(\nu ^{-1})=\overline{c(\bar\nu)}\,.$$
The correspondence sends a form $\tilde \omega$ to the function $c_{\tilde \omega}$ given by
$$ c_{\tilde \omega}(\nu):= \tilde \omega(e_\nu, e_{\nu^{-1}})\ ,
$$
where we have extended $\tilde \omega$ to a $\CC$-bilinear form on $V_\Bbb{C}$.

In order to relate the original form $\omega$ to the form $\omega'$ constructed above we associate to a polynomial $r\in \ZZ[x]$
a function $d_r: F\rightarrow \CC$ by
$$ d_r(\nu):=(r(\nu)+r(\nu^{-1}))\frac{c_{\omega'}(\nu)}{c_\omega(\nu)}\ .$$
Note that $d_r$ takes real values on $F_0$ and satisfies $d_r(\bar\nu)=\overline{d_r(\nu)}$. 

We claim that there exists a polynomial $r\in \ZZ[x]$ such that
\begin{equation}\label{max} d_r(\nu)\ne 0 \quad\mbox{ for all }\nu \in F 
\end{equation}
and (the crucial condition)
\begin{equation}\label{moritz} d_r(\nu)>0 \quad\mbox{ for }\nu \in F_0 \ .
\end{equation}
If $F_0$ is empty, then we can take $r=1$. Otherwise we find coefficients $a_k\in\RR$ such that 
$$ \sum_{k=0}^{N} a_k (\nu^k+\nu^{-k}) =  \frac{c_{\omega}(\nu)}{c_{\omega'}(\nu)}\quad\mbox{ for all }\nu \in F_0\subset S^1 \ .$$
Here $N$ can be taken such that $N+1$ is the number of elements of the set $F_0$ modulo complex conjugation.
We choose $b_k\in\QQ$ sufficiently close to $a_k$ such that 
$$ f(\nu):=  \sum_{k=0}^{N} b_k (\nu^k+\nu^{-k})$$
satisfies $ f(\nu)c_{\omega'}(\nu)/c_\omega(\nu)>0$ for all $\nu \in F_0 \ .$

Now take $q\in\ZZ$, $q>0$, such that $qb_k\in\ZZ$ for all $k$. We define
$$ r(x):=\sum _{k=0}^{N} qb_k x^k\ .$$
Then for $\nu\in F_0$,
$$ d_r(\nu)=qf(\nu)\frac{c_{\omega'}(\nu)}{c_\omega(\nu)}>0\ .$$
It remains to prove (\ref{max}). 
We consider the ring $\ZZ[x]/(f_A)$,
where $(f_A)$ denotes the ideal generated by the characteristic polynomial $f_A$ of $A$. Since $\det(A)=1$, the constant term of $f_A$ is equal to one, thus the element $x$ is invertible in $\ZZ[x]/(f_A)$. Hence, there is a polynomial $g\in\ZZ[x]$ such that $\nu^{-1}=g(\nu)$ for all $\nu\in F$. Set $h:=r+r\circ g\in \ZZ[x]$. Assume that $d_r(\nu)=0$ for some $\nu\in F$. Then $h(\nu)=r(\nu)+r(\nu^{-1})=0$.
Since 
$f_A$ is irreducible over $\QQ$ this implies $f_A|h$. But then $h(\nu)=0$ for all $\nu\in F$. This contradicts~(\ref{moritz}).

Let us fix $r\in\ZZ[x]$ such that (\ref{max}) and (\ref{moritz}) are satisfied. We define an element $D\in \GL(V_{\Bbb C})$
leaving $V\subset V_{\Bbb C}$ invariant and commuting with $A$ by
$$ De_\nu:=\left \{ \begin{array} {cl} d_r(\nu)^{-1}e_\nu\,,&\nu\in F_+\\
d_r(\nu)^{-1/2}e_\nu\,,&\nu\in F_0\\
  e_\nu\,,&\nu\in F_-\,.
  \end{array}\right.
 $$
 For $\nu\in F$ we compute
 $$ \omega'((r(A)+r(A^{-1}))De_\nu,De_{\nu^{-1}})=d_r(\nu)^{-1}(r(\nu)+r(\nu^{-1}))c_{\omega'}(\nu)=c_\omega(\nu)\ .$$
 It follows that for all $v,w\in V$
 $$ \omega(v,w)=\omega'((r(A)+r(A^{-1}))Dv,Dw)\ .$$
 Set $\Lambda:=D^{-1}\Lambda'$. Then $\Lambda\subset V$ is an $A$-stable lattice and
 $$ \omega(\Lambda\times\Lambda)=\omega'((r(A)+r(A^{-1}))\Lambda'\times\Lambda')\subset 
 \omega'(\Lambda'\times\Lambda')\subset \ZZ\ .
 $$
 \qed 
 
 In the following we will view the complex spectrum of a linear operator on a real vector
 space as well as the collection of all roots of a polynomial as multisets, i.e. sets equipped with a
 multiplicity function. Multisets will be denoted by $\{\dots\ms$, where the elements will be repeated
 according to their multiplicity. For instance, the roots of the polynomial $x^4+2x^2+1$ form the
 multiset $\{i,i,-i,-i\ms$.
 
We denote the spectrum of a linear operator $B:W\rightarrow W$ by $\spec(B)$. For any submultiset
$\rho\subset\spec(B)$ of cardinality $m$ let $G_m(W)^{B,\rho}\subset G_m(W)$ be the subset
of all $B$-invariant $m$-dimensional subspaces $V\subset W$ such that $\spec(B|_V)=\rho$. The following lemma will be needed later.

\begin{lm}\label{lehm}
Let $(W,\omega)$ be a real symplectic vector space, and let $B\in \fsp(W,\omega)$ be semisimple with purely imaginary spectrum. Then, for $\rho\subset\spec(B)$, the set 
$$G_m(W)^{B,\rho}_\mathrm{reg}:=\{V\in G_m(W)^{B,\rho}\mid \mathrm{rad}(\omega|_V)\subset\ker B\}$$
is dense in $G_m(W)^{B,\rho}$.
\end{lm}

\proof For a positive real number $\tau$ let  $W(\tau)$ denote the eigenspace
of $B^2$ with eigenvalue $-\tau^2$. Then we have a decomposition of $W$ into orthogonal $B$-invariant symplectic subspaces 
$$ W=\ker B \oplus \bigoplus_\tau W(\tau)\ .$$
On $W(\tau)$ we introduce the structure of a complex vector space equipped with a non-degenerate hermitian form $h_\tau$ as follows.
The complex structure is given by $J_\tau:=\frac{1}{\tau} B|_{W(\tau)}$, and we define $h_\tau(w_1,w_2):=\omega(w_1,J_\tau w_2)-i\omega(w_1,w_2)$. Now let $\rho\subset\spec(B)$. We may assume that $\rho$ is invariant under complex conjugation.
Otherwise $G_m(W)^{B,\rho}$ would be empty, and there is nothing to show. Let $m(\tau)$ be the multiplicity of $i\tau$ in $\rho$. The above decomposition of $W$ induces a homeomorphism
$$\Psi: G_m(W)^{B,\rho}\longrightarrow G_{\Bbb R,m(0)}(\ker B)\times\prod_{\tau} G_{\Bbb C, m(\tau)}(W(\tau))\ .$$
Here the subscripts $\RR$, $\CC$ indicate that we consider Grassmannians of real or complex subspaces, respectively.
Let $U_\tau$ be the subset of $G_{\Bbb C, m(\tau)}(W(\tau))$ consisting of subspaces that are non-degenerate with respect to $h_\tau$. Then 
$$\Psi\left( G_m(W)^{B,\rho}_{\mathrm{reg}}\right) =G_{\Bbb R,m(0)}(\ker B)\times\prod_{\tau} U_\tau\ .$$
The lemma now follows from the fact that $U_\tau\subset G_{\Bbb C, m(\tau)}(W(\tau))$ is dense.
\qed

\subsection{Special constellations}
Now we introduce certain parameters $P$, called blocks, that describe commuting linear operators $L,\phi$ on $\CC^d$ 
together with the spectrum $\rho$ of $L+\phi$ on some $(L,\phi)$-special subspace.
In particular, each block $P$ determines a dimension $d=d(P)$ and a vector $\mu(P)\in\RR^d$ such that $L=L(P)$ is given by $L_{\mu(P)}$ (see (\ref{L})).
Here, as in the Examples~\ref{B0} -- \ref{B2}, the vector space $\CC^d$ is equipped with the standard Euclidean inner product and with the involution given by complex conjugation. There are four types of blocks.

{\bf Type I.} $P=(0)$:  
$$\begin{array}{ll} d(P)=1&\rho(P)=\{0\ms\\
\mu(P)=0& \phi(P)=0\ .
\end{array}$$

{\bf Type II.} $\displaystyle P={\rho\choose \rho'}\in\RR^2\setminus\{0\}$:  
$$\begin{array}{ll} d(P)=2\qquad\rho(P)=\{\rho,-\rho\ms\\
\mu(P)=
\left(\frac{\rho+\rho'}{2},\frac{\rho+\rho'}{2}\right)\\ 
\phi(P)(z_1,z_2)=\frac{\rho-\rho'}{2}(-z_2,z_1)\ .
\end{array}$$

{\bf Type III.} $P=(\rho_1,\rho_2,\dots,\rho_{r-1},\rho_r)\in(\RR^*)^r,\ r\ge 2$, 
such that $|\rho_i|\ne |\rho_j|$ for $i\ne j$:  (for $r=2$ this row vector has to be distinguished from the column vector in Type II)
$$\begin{array}{ll} d(P)=2r\qquad\rho(P)=\{\rho_1,-\rho_1,\dots,\rho_r,-\rho_r\ms\\
\mu(P)=
\left(\rho_1,\frac{\rho_1+\rho_2}{2},\frac{\rho_1+\rho_2}{2},\dots,\frac{\rho_{r-1}+\rho_r}{2},\frac{\rho_{r-1}+\rho_r}{2},\rho_r\right)\\
\gamma(P)=(\gamma_1,\dots,\gamma_{r-1}),\ \mbox{ where } \gamma_i=(\rho_i-\rho_{i+1})/2\\
\phi(P)=0_1\oplus \phi_{\gamma(P)}\oplus 0_1.
\end{array}$$

{\bf Type IV.} $P=(\rho_1,\rho_2,\dots,\rho_{r-1},0)\in(\RR^*)^{r-1}\subset \RR^r,\ r\ge 2$, 
such that $|\rho_i|\ne |\rho_j|$ for $i\ne j$:  (for $r=2$ this row vector has to be distinguished from the column vector in Type II)
$$\begin{array}{ll} d(P)=2r-1\qquad\rho(P)=\{\rho_1,-\rho_1,\dots,\rho_{r-1},-\rho_{r-1},0\ms\\
\mu(P)=
\left(\rho_1,\frac{\rho_1+\rho_2}{2},\frac{\rho_1+\rho_2}{2},\dots,\frac{\rho_{r-2}+\rho_{r-1}}{2},\frac{\rho_{r-2}+\rho_{r-1}}{2},\frac{\rho_{r-1}}{2},\frac{\rho_{r-1}}{2}\right)\\ 
\gamma(P)= (\gamma_1,\dots,\gamma_{r-2},\rho_{r-1}/2 ), \mbox{ where } \gamma_i=(\rho_i-\rho_{i+1})/2,\ i=1,\dots, r-2,\\
\phi(P)=0_1\oplus\phi_{\gamma(P)}\,.
\end{array}$$
A block of Type II is called of {\bf Type II.a} if $\rho=0$, of {\bf Type II.b} if $\rho=\rho'$, and of {\bf Type II.c} if $\rho=-\rho'$.

A {\bf special constellation} $\cP=(P_1|P_2|\dots|P_l)$ is a direct sum of blocks having the following property:
For each non-negative real number $\rho_0$ there is at most one block of Type III or Type IV appearing in $\cP$, which starts or ends with 
$\pm \rho_0$. In particular, there appears at most one block of Type IV. We set 
\begin{eqnarray*} 
d(\cP)=\sum_{k=1}^l d(P_k)\qquad&&\rho(\cP)=\rho(P_1)\uplus\rho(P_2)\uplus\dots\uplus\rho(P_l)\\
\mu(\cP)=(\mu(P_1),\dots,\mu(P_l))\in\RR^{d(\cP)}\\
L(\cP)=L_{\mu(\cP)}=\bigoplus_{k=1}^l L(P_k)\qquad&& \phi(\cP)=\bigoplus_{k=1}^l \phi(P_k)\ .
\end{eqnarray*}
The empty constellation is considered as a special constellation of dimension $d(\cP)=0$.

The relevance of special constellations is explained by the following close relative of Proposition~\ref{special}.

\begin{pr}\label{veryspecial}
Let $(W,\theta)$ be a Euclidean space of dimension $2n$ with involution as in Definition~\ref{henze}. Let $L,\phi\in\so(W)$ be two commuting elements, where $L$ anticommutes while $\phi$ commutes with $\theta$. Let $\sigma$ be a multiset of cardinality $n$
supported on $\RR$. Then there exists an $(L,\phi)$-special subspace $V\subset W$ such that
$\spec((L+\phi)|_V)=i\sigma$ if and only if $(W,\theta,L,\phi)$ is isometrically isomorphic to
$(\CC^{d(\mathcal P)}, \mathop{\rm conj}, L(\cP),\phi(\cP))$ for some special constellation $\cP$ satisfying $\rho(\cP)=\sigma$. 
\end{pr}
\proof Examples \ref{B0}, \ref{B1} and \ref{B2} show that if $(W,\theta,L,\phi)$ is isometrically isomorphic to a block $P$ of Type II, III or IV, then $W$ contains a special subspace $V\subset W$ such that $\spec(L+\phi)|_V=i\rho(P)$. Obviously, the same is true for blocks of Type I. Since each special constellation $\cP$ is a direct sum of such blocks, this proves the existence of a special subspace with spectrum $i\rho(\cP)$. 

Now suppose that $W$ admits an $(L,\phi)$-special subspace $V$ and that $\spec((L+\phi)|_V)=i\sigma$. First we will show that $(W,\theta,L,\phi)=(W_0,\theta_0,L_0,\phi_0)\oplus(W_1,\theta_1,L_1,\phi_1)$, where 
\benur
\item $(W_0,\theta_0,L_0,\phi_0)$ is isometrically isomorphic to
$(\CC^{d(\cP_0)}, \mathop{\rm conj}, L(\cP_0),\phi(\cP_0))$ for a special constellation $\cP_0$ that is a direct sum of blocks of type  I,
\item $W_1$ admits an $(L_1,\phi_1)$-special subspace $V_1$, such that 
\begin{equation}\label{EZerl}
i\sigma=i\rho(\cP_0)\uplus\, \spec((L_1+\phi_1)|_{V_1}),
\end{equation}
\item $\Ker L_1\cap \Ker \phi_1=0$.
\end{enumerate}
We define $W_0:=\ker L \cap\ker \phi$ and $W_1:=W_0^\perp\subset W$. For $j=0,1$, let $\theta_j,L_j,\phi_j$ be the linear maps induced by $\theta,L$ and $\phi$ on $W_j$.  Obviously, (i) and (iii) are satisfied. Let $p_j:W\rightarrow W_j$, $j=0,1$, denote the orthogonal projections. We define $p:V\rightarrow (W_0)_-$ by $p(v)=p_0(v)_-$, where $p_0(v)_-$ denotes the projection of $p_0(v)\in W_0$ to $(W_0)_-$. Then $p$ is surjective since the projection of $V$ to $W_-$ is surjective. Now we put $\tilde V_1:=\ker p$. Since $p$ is surjective, we get $\dim \tilde V_1= \dim V-\dim (W_0)_-=(\dim W-\dim W_0)/2=\dim W_1/2$. Note that $\tilde V_1$ is $(L+\phi)$-invariant because of $p_0((L+\phi)(v)) =0\in (W_0)_-$. Now observe that $p_1|_{\tilde V_1}$ is injective. Indeed, if $p_1(v)=0$ for $v\in \tilde V_1$, then $v\in (W_0)_+$, thus $v\in V\cap W_+=0$. We will show that $V_1:=p_1(\tilde V_1)\subset W_1$ is $(L_1,\phi_1)$-special. Indeed, $(L_1+\phi_1)(V_1)=(L_1+\phi_1)p_1(\tilde V_1)= p_1(L+\phi)(\tilde V_1)\subset p_1(\tilde V_1)=V_1$ shows that $V_1$ is $(L+\phi)$-invariant. We have $\dim V_1=\dim \tilde V_1=\dim W_1/2$. Furthermore, $p_0(\tilde V_1)\subset (W_0)_+$ implies $V_1\cap (W_1)_+=p_1(\tilde V_1\cap W_+)=0$. Hence, $W_1=V_1\oplus (W_1)_+$. Finally, Equation (\ref{EZerl}) holds since $L+\phi$ acts trivially on $V/\ker p$. 

By the above considerations it remains to consider the case $\ker L \cap\ker \phi=0$.  
We will use the notation introduced in Subsection~\ref{SSS}. We may assume  $W= \CC^p\oplus\CC^{2q}$, $L= L_\alpha\oplus L_\beta$ and $\phi= 0_p\oplus \phi_\gamma$. Then $\alpha_k\not=0$ for all $k=1,\dots,p$. We have a further decomposition 
\begin{equation}\label{Edec1}
W=\bigoplus_{k=1}^p W^1_k\oplus \bigoplus_{j=1}^q W^2_j,
\end{equation}
 where $W^1_k$ and $W^2_j$ are defined by (\ref{EW1}) and (\ref{EW2}). 

We want to weaken the notion of a special subspace in the following way. Let $\cZ$ be a decomposition $W=\bigoplus_{\nu=1}^r W_\nu$ of the $(L+\phi)$-module $W$ into $(L+\phi)$-invariant subspaces. Typically, $W_\nu$ will be a direct sum of spaces of type $W^1_k$ and $W^2_j$. For a subset $J\subset \{1,\dots,r\}$ we define $W_J:= \bigoplus_{\nu\in J} W_\nu$. An $(L+\phi)$-invariant subspace $V\subset W$ is called {\it pseudo-special with respect to $\cZ$} if $\dim V= \frac 12 \dim  W$ and if 
$$
\dim (V\cap W_J) \le {\textstyle \frac12} \dim W_J 
$$
for all $J\subset \{1,\dots,r\}$. 
Let $\pro_J$ denote the projection from $W$ to $W_J$ with respect to $W=W_J\oplus W_{\bar J}$, where $\bar J=\{1,\dots,r\}\setminus J$. 

Every $(L,\phi)$-special subspace is pseudo-special with respect to~(\ref{Edec1}).

{\sl Claim 1. } Suppose that $p$ is even and that $\cZ'$ is a decomposition of $W$ into $(L+\phi)$-invariant subspaces $W_\nu$, $\nu=1,\dots, r$, of (real) dimension four. Let $V\subset W$ be pseudo-special with respect to $\cZ'$. Then, for each $\nu$, we can choose one of the two pairs of conjugate eigenvalues of $L+\phi$ on $W_\nu$ such that the multiset of the chosen pairs $\pm i\tau_\nu$ equals the spectrum of $L+\phi$ on $V$:
$$
\spec((L+\phi)|_V)=\{\pm i \tau_1,\dots, \pm i \tau_r\ms\,.
$$
{\sl Proof of Claim 1. } We prove the claim by induction on $\dim_{\Bbb C}W$. For $\dim_{\Bbb C}W=2$ the assertion is obvious. For $\dim_{\Bbb C}W>2$, we consider non-empty complementary subsets $J_1$ and $J_2$ of $\{1,\dots,r\}$. The decomposition $\cZ$ of $W$ induces decompositions of $W_{J_1}$  and $W_{J_2}$. Thus we can speak of pseudo-special subspaces of $W_{J_1}$  and $W_{J_2}$ with respect to these decompositions. In order to prove the claim, it suffices to show that we can decompose $\{1,\dots,r\}=J_1\cup J_2$ and find pseudo-special subspaces $V_1\subset W_{J_1}$ and $V_2\subset W_{J_2}$ such that $\sigma(V)=\sigma(V_1)\uplus \sigma(V_2)$. 

We consider the following three cases, where $J$ always denotes a non-empty proper subset of $\{1,\dots,r\}$ and $\bar J$ its complement:

Case 1: $\dim (V\cap W_J) = \frac12 \dim W_J$ for some $J$. In this case we can choose $J_1=J$, $J_2=\bar J$, $V_1=V\cap W_{J_1}$, $V_2=\pro_{J_2}(V)$. 

Case 2: $\min_J (\frac12 \dim W_J-\dim (V\cap W_J) )= 1$. Let $J_1$ be such that $\dim (V\cap W_{J_1}) = \frac12 \dim W_{J_1}-1$.
Then $\dim (V\cap W_{J_1})$ is odd, thus $(L+\phi)|_{V\cap W_{J_1}}$ has an odd-dimenional kernel. On the other hand, $\dim W_{J_1}$ is even, hence the kernel of $(L+\phi)|_{W_{J_1}}$ also contains an element $v\not\in V$, $v\not=0$.  
We put $J_2=\overline{J_1}$, $V_1=(V\cap W_{J_1})\oplus \RR\cdot v$. Since $\dim (W_{J_2}\cap \pro_{J_2}V)=\frac 12 \dim W_{J_2}+1$ is odd, $L+\phi$ has a non-trivial kernel on  $V':=W_{J_2}\cap \pro_{J_2}V$. We choose $V_2\subset V'$ to be an $(L+\phi)$-invariant complement of a one-dimensional subspace of this kernel.

Case 3: $\dim (V\cap W_J) \le \frac12 \dim W_J-2$ for all $J$. Then we put $J_1=\{1\}$, $J_2=\overline{J_1}$, thus $W_{J_1}=W_1$. By assumption, $\dim (V\cap W_{J_2}) \le \frac12 \dim W_{J_2}-2$, thus $\dim \pro_{J_1}V=4$. This implies that both pairs of eigenvalues $\pm i\tau_1$, $\pm i\tau'_1$ of $L+\phi$ on $W_1$ belong to the spectrum of $V$. We choose $V_1\subset W_{1}$ to be a two-dimensional $(L+\phi)$-invariant subspace (i.e., we choose a pair of eigenvalues $\pm i\tau_1$). Let $V'\subset W_1$ be the $(L+\phi)$-invariant subspace that is complementary to $V_1$ (i.e., the one with eigenvalues $\pm i\tau_1'$). Then we put $V_2:=\pro_{J_2}((W_1\oplus V')\cap V)$. In order to verify that $V_2$ is pseudo-special, we use that $\pro_{J_2}$ is injective on $(W_1\oplus V')\cap V$. 
\qed\\
Let us proceed with the proof of Prop.~\ref{veryspecial}. Let us first consider the case $p=0$. Since any $(L,\phi)$-special subspace of $W$ is pseudo-special with respect to $W=\bigoplus_{j=1}^q W^2_j$, we can apply Claim 1. We obtain a choice of eigenvalues $\pm i\tau_j$, $j=1,\dots,q$.  We set $\rho_j:=\tau_j$. 
Recall that $\tau_j=\kappa \beta_j+\kappa'\gamma_j$ for suitable $\kappa,\kappa'\in\{1,-1\}$ and put $\tau_j'=\kappa \beta_j-\kappa'\gamma_j$. Then $\pm i\tau'_j$ is the second pair of eigenvalues of $L+\phi$ on $W^2_j$ besides $\pm i\tau_j$.  We put $\rho_j':=\tau'_j$.
Then $P_j=$ {\small ${\displaystyle {\rho_j\choose\rho_j'}}$} is a block of Type II and, for the special constellation $\cP=(P_1|\dots|P_q)$, we have $(W,\theta,L,\phi)\cong(\CC^{d(\mathcal P)}, \mathop{\rm conj}, L(\cP),\phi(\cP))$ and $\rho(\cP)=\sigma$.

Next we turn to the case $p>0$, which is more involved. Let us first assume that $p$ is even. In this case, we can combine the spaces $W^1_k$, $k=1,\dots,p$, pairwise and we obtain a decomposition $\cZ'$ of $W$ into 4-dimensional $(L+\phi)$-invariant subspaces: 
$$W_\nu := \left\{\begin{array}{ll} W^1_{2\nu-1}\oplus W^1_{2\nu}\,, & \nu=1,\dots,p/2\,,\\ W^2_{\nu-p/2}\,, & \nu=p/2+1,\dots,p/2+q \,.\end{array}\right.$$
Since $V$ is pseudo-special with respect to the finer decomposition (\ref{Edec1}), it is also pseudo-special with respect to $\cZ'$. 
Hence we may apply Claim 1 and obtain a choice of eigenvalues $\pm i\tau_\nu$, $\nu=1,\dots,p/2+q$. Let $\pm i\tau'_{\nu}$ denote the second pair of eigenvalues of $L+\phi$ on $W_{\nu}$. Obviously, we have $\tau_\nu,\tau'_\nu\in\{\pm\alpha_{2\nu-1},\pm\alpha_{2\nu}\}$ for $\nu =1,\dots,p/2$.
For $j=1,\dots,q$, we define 
\begin{equation}\label{Err}
\rho_j:=\tau_{p/2+j},\quad \rho_j':=\tau'_{p/2+j}.
\end{equation}

As in the proof of Prop.~\ref{special}, we assign to $(W,L+\phi)$ the graph $G$ defined by (\ref{EVE}), where now $\rho_j$ and $\rho'_j$ are given by (\ref{Err}).  Using Claim 1 we will $G$ endow with an orientation and a charge. A charge $w$ on a graph is a vertex-weight taking values in $\{1,-1,0\}$. We will say that a vertex is positively charged if its weight is $1$ and that it is negatively charged if its weight is $-1$. Positively and negatively charged vertices will be visualised by $\pv$ and $\mv$, respectively.

Orientation and charge of $G$ are defined in the following way.  Each edge $j$ is directed from $|\rho_j'|$ to $|\rho_j|$. A vertex $|\rho|$ is charged if and only if $|\rho|=|\alpha_k|$ for some $k=1,\dots,p$, and for these vertices the charge $w$ is defined by
$$w(|\rho|)=\left\{ \begin{array}{rl} 1\,, &\mbox{if $|\rho|=|\tau_{\nu}|$ for some $\nu\in\{1,\dots,p/2\}$\,,}\\-1\,, & \mbox{if $|\rho|=|\tau'_{\nu}|$ for some $\nu\in\{1,\dots,p/2\}.$}  \end{array}\right.$$

Now $\spec((L+\phi)|_V)$ can be read off of the charged directed graph $G$ in the following way. 
For $\rho\ne 0$, let $m(|\rho|)$ be the multiplicity
of $i\rho$ in $\spec((L+\phi)|_V)$. In order to get uniform formulas, we define $m(0)$ to be half of the multiplicity of $0$ in $\spec((L+\phi)|_V)$. Then $m(|\rho|)$ is related to the in-degree $\deg^-(|\rho|)$ of $|\rho|$ by 
\begin{equation}\label{Ein} 
m(|\rho|)=\left\{ \begin{array}{ll} \deg^-(|\rho|)+1\,, &\mbox{if $|\rho|$ is positively charged,}\\ \deg^-(|\rho|)\,, &\mbox{if $|\rho|$ is not positively charged.} \end{array}\right.
\end{equation}
The term on the right hand side of this equation will be called {\em charged in-degree}.

A path in a directed graph is an alternating sequence $v_0e_1v_1e_2\dots e_rv_r$ of pairwise distinct vertices $v_0,\dots,v_r$  and edges $e_1,\dots, e_r$ beginning and ending with a vertex such that the vertex that precede an edge is the starting point and the vertex that follows an edge is the end point of that edge.
A charged path in a charged oriented graph is a path starting with a positively charged vertex and ending with a negatively charged vertex. Usually, we identify a path with its underlying graph. Below we will draw a path without its starting point and endpoint as a bold arrow 
\begin{tikzpicture}
  \node (1) {};
  \node (2) [right of=1] {};
  \draw[->, very thick] (1) to (2);
\end{tikzpicture}  
omitting all inner vertices.

The plan is to arrange the edges into disjoint sets such that every set defines a block of the special constellation $\cP$ we are looking for. Some of these sets will be edge sets of charged paths, the remaining sets consist of one edge only. 
To do so, we wish to join all positively charged points to negatively charged ones by edge-disjoint paths such that each charged vertex is the starting point or the end point of exactly one of these paths. Each of these paths will define a block of Type III. 
Every edge that does not belong to one of these paths will define a �block of Type II. 

However, in general, joining the positive charged vertices to negatively charged ones by edge-disjoint paths is not possible for our original graph $G$. To reach our goal, we will have to modify charge and orientation of $G$ according to the following procedure, which preserves the charged in-degree of each vertex. Let be given an oriented charged graph and paths $t_1,\dots,t_j$ in this graph. An {\em orientation modification} of this graph is obtained by changing the direction of all edges of $t_1,\dots,t_j$ and changing the charge of all starting and end points of these paths.  

We will say that a set $A\subset{\cal V}$ is positively charged if $w(A):=\sum_{v\in A} w(v)>0$.

{\sl Claim 2. } If $A\subset{\cal V}$ is positively charged, then there exists an edge that is directed from a vertex in $A$ to a vertex in ${\cal V}\setminus A$. 

{\sl Proof of Claim 2. } Let $W'\subset W$ be an $(L+\phi)$-invariant subspace. For an $(L+\phi)$-invariant subspace $V'\subset W'$ of $W$ we denote by $m_{W'}(V',|\rho|)$ the multiplicity of $i\rho$ as an eigenvalue of $(L+\phi)|_{V'}$, and by $m'_{W'}(V',|\rho|)$ the multiplicity of $i\rho$ as an eigenvalue of $L+\phi$ on $W'/V'$. 
For $\rho=0$ we modify these definitions by a factor $\frac{1}{2}$ as above.
The relation of $m'(|\rho|):=m'_{W}(V,|\rho|)$ to the out-degree $\deg^+(|\rho|)$ of $|\rho|$ is similar to that of $m(|\rho|)=m_{W}(V,|\rho|)$ to the in-degree:
$$ m'(|\rho|)=\left\{ \begin{array}{ll} \deg^+(|\rho|)+1\,, &\mbox{if $w(|\rho|)<0$,}\\ \deg^+(|\rho|),& \mbox{if $w(|\rho|)\ge0$.} \end{array}\right.
$$
We set 
\begin{eqnarray*}
\ind(|\rho|)&:=&m(|\rho|)-m'(|\rho|)\\
\mbox{i-}\ind_{W'}(|\rho|)&:=&\,m_{W'}(W'\cap V,|\rho|)-m_{W'}'(W'\cap V,|\rho|)\,\\
\mbox{p-}\ind_{W'}(|\rho|)&:=& \,m_{W'}(\pro_{W'}V,|\rho|)-m_{W'}'(\pro_{W'}V,|\rho|)\,,
\end{eqnarray*}
where $\pro_{W'}\in\End(W)$ is the orthogonal projection to $W'$.
Then we have
\begin{eqnarray}
\ind(|\rho|)&=&\deg^-(|\rho|)-\deg^+(|\rho|)+w(|\rho|), \label{luemmel}\\
&=&\mbox{i-}\ind_{W'}(|\rho|)+ \mbox{p-}\ind_{(W')^\perp}(|\rho|).\label{Eind}
\end{eqnarray}
Moreover, for $R':=\{|\rho|\mid \mbox{$i\rho$ is eigenvalue of } (L+\phi)|_{W'}\, \}$, we get
$$\sum_{|\rho| \in R'}\mbox{i-}\ind_{W'}(|\rho|)=\mbox{i-}\ind(W'):=\textstyle{\frac12}(\,\dim(W'\cap V)-\dim (W'/(W'\cap V)\,).$$
Let $e_+$ be the number of edges going out of $A$ and let $e_-$ be the number of edges going into $A$.  
Equation (\ref{luemmel}) yields 
\begin{equation}\label{Einda}
\sum_{|\rho|\in A}\ind(|\rho|)=w(A)+e_--e_+\,.
\end{equation}
Let ${\cal E}'$ denote the set of edges between vertices in $A$. We put
$$W':= \bigoplus_{k,|\alpha_k|\in A}W^1_k\ \oplus\ \bigoplus_{j\in{\cal E}'} W^2_j, \quad W'':= (W')^\perp$$
and apply (\ref{Eind})  to the left hand side of (\ref{Einda}). Since all eigenvalues of $L+\phi$ on $W'$ are in $A$, this gives
\begin{equation}\label{Einda1}
\sum_{|\rho|\in A}\ind(|\rho|)= \mbox{i-}\ind(W')+\sum_{|\rho|\in A}\mbox{p-}\ind_{W''}(|\rho|)\,.\\
\end{equation}
Furthermore, 
\begin{equation}\label{Einda2}
\sum_{|\rho|\in A}\mbox{p-}\ind_{W''}(|\rho|)\le \sum_{|\rho|\in A} m_{W''}(\pro_{W''}V,|\rho|)\le e_++e_-\,. 
\end{equation}
Indeed, let $j_1,\dots,j_e$, $e=e_++e_-$, be the edges between $A$ and the complement of $A$. 
Numbers $i\rho$ with $|\rho|\in A$ can be eigenvalues of $(L+\phi)|_{W''}$ only if 
$$\pm\rho\in\{\rho_{j_\nu}, \rho'_{j_\nu}\mid \nu=1,\dots,e\}\, .$$ 
Moreover, for fixed $\nu$, $1\le\nu \le e$, only one of the numbers $|\rho_{j_\nu}|, |\rho'_{j_\nu}|$ can be in $A$ since otherwise $j_\nu$ would be an inner edge of $A$. 

Equations (\ref{Einda}), (\ref{Einda1}) and (\ref{Einda2}) yield
$$2e_+\ge w(A)-\mbox{i-}\ind(W').$$
Since $V$ is pseudo-special, we have $\mbox{i-}\ind(W')\le 0$ and the assertion follows.
\qed\\
{\sl Corollary.} Let $A$ and $B$ be disjoint sets of vertices of $G$ such that $A\cup B$ contains all charged vertices of $G$. If $A$ is positively charged, then there exists a path starting in $A$ and ending in $B$.

\proof Consider the set $\tilde A$ of vertices that can be reached by a path starting in $\tilde A$. Then no edge goes out of $\tilde A$.  Thus $\tilde A$ is not positively charged. Hence $\tilde A$ contains a vertex of~$B$.\qedohne

A connected charged directed graph is called of Type $\cN_j$ if it is the union of pairwise edge-disjoint charged paths $t_\nu$, $\nu=1,\dots,j$ with distinct starting points and distinct end points and if it does not contain any charged vertex that is neither a starting nor an end point of one of the paths $t_\nu$. We will say that a graph is of Type $\cN$ if it is of Type $\cN_j$ for some $j\in\NN$. We will write $G\in\cN_j$ and $G\in\cN$, respectively. The figure below shows an example of a graph of Type $\cN_3$.
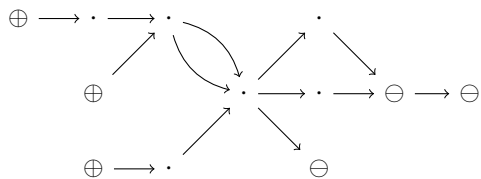
\begin{center}
\begin{tikzpicture}
  \node (1) {$\pv$}; \node (1a)[right of=1] {$\cdot$}; \node (1b) [right of=1a] {$\cdot$}; \node (1c) [right of=1b] {}; \node (1d) [right of=1c] {$\cdot$};
  \node (2a) [below of=1a] {$\pv$}; \node (2b) [right of=2a] {}; \node (2c) [right of =2b] {$\cdot$}; \node (2d) [right of=2c] {$\cdot$};\node (2e) [right of=2d] {$\mv$};\node (2f) [right of=2e] {$\mv$};
  \node (3a) [below of=2a] {$\pv$}; \node (3b) [right of=3a] {$\cdot$}; \node (3c) [right of=3b] {}; \node (3d) [right of=3c] {$\mv$};
  \draw[->] (1) to (1a);\draw[->] (1a) to (1b);\draw[->,bend left] (1b) to (2c);\draw[->,bend right] (1b) to (2c); \draw[->] (2c) to (2d); 
  \draw[->] (2a) to (1b);\draw[->] (2c) to (3d);  \draw[->] (2c) to (1d);  \draw[->] (2d) to (2e); \draw[->] (2e) to (2f); 
  \draw[->] (3a) to (3b);\draw[->] (3b) to (2c);  \draw[->] (1d) to (2e);
\end{tikzpicture}
\end{center}
 
A charged directed graph is called of Type $\cQ$ if it belongs to one of the sets $\cQ_j$ defined by the following induction. We start with $\cQ_1$, whose only element is the graph that consists of just one positively charged vertex and has no edges. Elements of $\cQ_j$, $j>1$, are graphs of the form $Q \cup t \cup N$, where $Q\in \cQ_{l}$ ($1\le l<j$) and $N\in\cN_{j-l}$ are disjoint and $t$ is a path starting in $Q$ and ending in $N$ such that $Q\cap t$ and $N\cap t$ consist of exactly one vertex and such that none of its inner vertices is charged:
$$\cQ_1=\{ \pv \},\quad \cQ_j=\left\{ \ \left.
\begin{tikzpicture} [baseline={([yshift=-7.5ex]current bounding box.north)}, node distance=4em, state/.style ={circle, draw}]
  \node[state] (1) {$Q$};
  \node[state] (2) [below of=1] {$N$};
  \draw[->, very thick] (1) edge node [right] {$t$} (2); 
\end{tikzpicture}
 \ \ \right| \ Q\in \cQ_{l}, \ N\in\cN_{j-l},\  1\le l<j \right\}\,.$$
Each graph $Q\in\cQ_j$ contains exactly $j+1$ positively charged and $j$ negatively charged vertices. In particular, $w(Q)=1$ for all graphs $Q$ of type $\cQ$.

{\sl Claim 3.} Let be given a charged oriented graph that equals the union $Q\cup t_-$ of a graph $Q$ of type $\cQ$ and a path $t_-$ with starting point in $Q$ and negatively charged end point not containing any inner charged vertex. Moreover, suppose that $Q$ and $t_-$ have in common only the starting point of $t_-$:
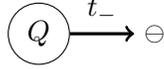
\begin{center}
\vspace{-1ex}
\begin{tikzpicture} [baseline={([yshift=-7.5ex]current bounding box.north)}, node distance=4em, state/.style ={circle, draw}]
  \node[state] (1) {$Q$};
  \node (2) [right of=1] {$\mv$};
  \draw[->, very thick] (1) edge node [above] {$t_-$} (2); 
\end{tikzpicture}
\end{center}
\vspace{-2ex}
Then, possibly after an orientation modification, this graph contains a graph whose connected components are of type $\cN$ and which contains all charged points of $Q\cup t_-$.

{\sl Proof of Claim 3. } We will prove the claim for $Q\in\cQ_j$ by induction on $j$. For $j=1$ the assertion is obvious, since $\cQ_1=\{\pv\}$. Now assume $j>1$. Let $Q_j$ be in $\cQ_j$, i.e., $Q_j=Q\cup t_+ \cup N$, where $Q\in \cQ_l$, $N\in \cN_{j-l}$ for some $1\le l<j$ and $t_+$ is the connecting path. If $t_-$ starts from $Q$, then the assertion for $j$ follows immediately from the induction hypothesis for $l<j$ applied to $Q\cup t_-$. If $t_-$ starts from an inner point of $t_+$, we proceed in a similar way.  Let us now consider the case that $t_-$ starts from a vertex in $N$. Let $i$ be the endpoint of $t_+$ and let $o$ be the starting point of $t_-$. By definition, $N$ is the union of edge-disjoint charged paths $t_1,\dots,t_{j-l}$. If there exists a path $t_k$, $1\le k \le j-l$, such that the vertices $i$ and $o$ belong to the same connected component $N_o$  of $N- t_k$, then we can apply the induction hypothesis to $Q\cup t_+\cup N_o$. Here $N- t_k$ denotes the charged oriented graph $\bigcup_{\nu\not=k} t_\nu$ (without the charges of the endpoints of $t_k$ if these endpoints belong to some $t_\nu$ for $\nu\not=k$).
%
Now we assume that there does not exist such a path $t_k$. 

The following figure illustrates the situation for $N=\{t_1\}$:
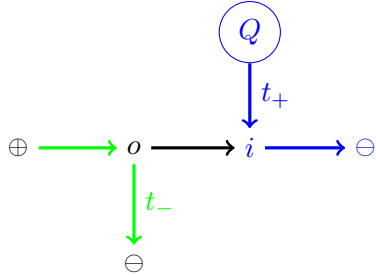
\begin{center}
\begin{tikzpicture} [node distance=4em, state/.style ={circle, draw}]
  \node (1a){}; \node (1b) [right of=1a] {}; \node[state, blue] (1c) [right of=1b] {$Q$}; \node (1d) [right of=1c] {};
  \node (2a) [below of=1a] {$\pv$}; \node (2b) [right of=2a] {$o$}; \node[blue] (2c) [right of =2b] {$i$}; \node[blue] (2d) [right of=2c] {$\mv$};
  \node (3a) [below of=2a] {}; \node (3b) [right of=3a] {$\mv$};
  \draw[->, very thick, blue] (1c) edge node [right] {$t_+$} (2c);
  \draw[->, very thick, green] (2a) to (2b);\draw[->,very thick] (2b) to (2c);  \draw[->, very thick, blue] (2c) to (2d);  
  \draw[->,very thick, green] (2b) edge node [right] {$t_-$} (3b);
\end{tikzpicture}
\end{center}
In this picture, $t_1$ is the horizontal path. The green edges define a charged path without inner charged points. To the blue subgraph we can apply the induction hypothesis. Here we have assumed that $o$ lies before $i$ on $t_1$ or is equal to $i$. If this is not the case, we modify the orientation along~$t_1$. 
Now let $N$ consist of at least two edge-disjoint paths $t_1,\dots, t_{j-l}$.  We choose a path $t\in\{t_1,\dots,t_{j-l}\}$ that contains $i$. By assumption, $o\not\in t$. Moreover, the intersection of $t$ and the connected component $N_o$ of $o$ in $N- t$ contains at least one vertex. After possibly modifying the orientation along $t$, we may assume that  the first vertex $n_0$ of $t\cap N_o$ on $t$ is located before $i$ on $t$. We denote by $t'$ the segment of $t$ from $\pv$ to $n_0$ and by $t''$ the segment from $i$ to $\mv$. Now we can apply the induction hypothesis to $t'\cup N_o\cup t_-$ (green subgraph in the figure below) ignoring all charges of inner points of $t'$ and to $Q\cup t_+\cup t_-$ (blue subgraph) ignoring all charges of inner points of $t_+\cup t''$:
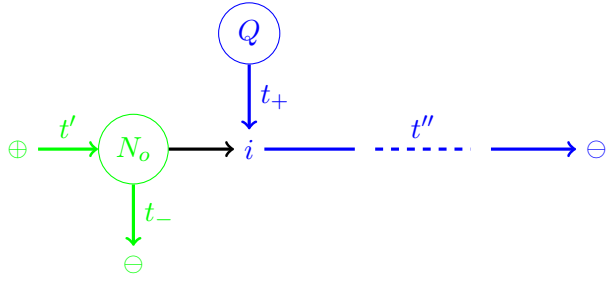
\begin{center}
\begin{tikzpicture} [node distance=4em, state/.style ={circle, draw}]
  \node (1a){}; \node (1b) [right of=1a] {}; \node[state, blue] (1c) [right of=1b] {$Q$}; \node (1d) [right of=1c] {};
  \node[green] (2a) [below of=1a] {$\pv$}; \node[state, green] (2b) [right of=2a] {$N_o$}; \node[blue] (2c) [right of =2b] {$i$}; \node (2d) [right of=2c] {}; \node (2e) [right of=2d] {}; \node[blue] (2f) [right of=2e] {$\mv$};
  \node (3a) [below of=2a] {}; \node[green] (3b) [right of=3a] {$\mv$};
  \draw[->, very thick, blue] (1c) edge node [right] {$t_+$} (2c);
  \draw[->, very thick, green] (2a) edge node [above] {$t'$} (2b);\draw[->,very thick] (2b) to (2c);  \draw[-, very thick, blue] (2c) to (2d);  \draw[-, very thick, blue, dashed] (2d) edge node [above] {$t''$} (2e); \draw[->, very thick, blue] (2e) to (2f); 
  \draw[->,very thick, green] (2b) edge node [right] {$t_-$} (3b);
\end{tikzpicture}
\end{center}
Here, $t$ equals the horizontal line. The dashed line indicates that $t''$ and $N_o$ may have vertices in common.  The connected components of $N-t$ different from $N_o$ remain unchanged.
\qed\\
We want to define sequences $(A_j)_{j\in{\Bbb N}}$ and $(B_j)_{j\in{\Bbb N}}$ of subgraphs of a graph $G_j$ that arises from $G$ by orientation modification. These graphs will have the following properties:
\benur
\item  $A_j$ is a disjoint union of graphs of type $\cQ$;
\item  $B_j$ is a disjoint union of graphs of type $\cN$ and graphs that consist only of one negatively charged vertex;
\item  $A_j\cap B_j=\emptyset$ and $A_j\cup B_j$ contains all charged vertices of $G_j$.
\end{enumerate}
The graph $A_1$ consists of all positively charged vertices of $G_1:=G$ and $B_1$ consists of all negatively charged vertices. Both graphs do not have edges. Clearly, (i) -- (iii) are satisfied. Suppose, $A_j$, $B_j$ and $G_j$ are already defined. If $A_j=\emptyset$, put $A_{j+1}=A_j$, $B_{j+1}=B_j$ and  $G_{j+1}=G_j$. Now assume $A_j\not=\emptyset$. Then $w(A_j)>0$ by (i). The corollary of Claim~2 says that there is a path $t$ from $A_j$ to $B_j$. We may assume that $A_j\cap t$ and $B_j\cap t$ contain only the starting point  and the end point of $t$, respectively. Then $t$ does not contain inner charged vertices since all charged vertices are in $A_j\cup B_j$.
Let $A$ be the connected component of $A_j$ that contains the starting point of $t$ and let $B\subset B_j$ be the connected component containing the end point. We put $\hat G:= A\cup t\cup B$. 

Let us first consider the case that $B$ is of Type $\cN$. Then $\hat G$ is of Type $\cQ$ and we put 
\begin{equation}\label{EA1}
G_{j+1}:=G_j,\quad A_{j+1}:=(A_j\setminus A)\cup \hat G,\quad B_{j+1}:= B_j\setminus B.
\end{equation}
Then (i)-(iii) are obviously satisfied. 

Now suppose that $B$ consists of a single vertex. Claim~3 implies that, after an orientation modification, $\hat G$ contains a graph $N$ whose connected components are of Type $\cN$ and which contains all charged vertices of $\hat G$. We consider the  orientation modification of $\hat G$ as a modification of $G_j$ and denote the resulting graph by $G_{j+1}$. We define
\begin{equation}\label{EA2} 
A_{j+1}:=A_j\setminus A,\quad B_{j+1}:= (B_j\setminus B)\cup N.
\end{equation}
Then the conditions (i)--(iii) are satisfied by construction.   

We claim that the sequences $(A_j)_{j\in{\Bbb N}}$ and $(B_j)_{j\in{\Bbb N}}$ stabilise. Indeed, let us determine the charge $w(A_j)$, which equals the number of connected components of $A_j$.  If $A_{j+1}$ is defined by (\ref{EA2}), then $w(A_{j+1})=w(A_j)-1$. If it is defined by (\ref{EA1}), then $w(A_j)=w(A_{j+1})$. However, in the latter case the step from $B_{j}$ to $B_{j+1}$ reduces the number of connected components of Type $\cN$. Thus, after finitely many steps of this kind, we have to apply again (\ref{EA2}). Hence, there exists an index $j_0$ such that $w(A_{j})=0$ for all $j\ge j_0$. In particular, $A_{j_0}=\emptyset$. Thus $B_{j_0}$ consists of $p/2$ edge-disjoint paths $t_1,\dots,t_{p/2}$ that join all positively charged vertices of $G_{j_0}$  to distinct negatively charged ones. 

Let $t_\nu$ be one of the charged paths constituting $B_{j_0}$. Suppose that $t_\nu$ starts at the positively charged vertex $|\alpha_k|$ and ends at the negatively charged vertex $|\alpha_{k'}|$ and let $j_2,j_3,\dots, j_r$ be the sequence of its edges. Then, by construction, 
\begin{equation} \label{Ear} 
|\alpha_k|=|\rho'_{j_2}|,\ |\rho_{j_2}|=|\rho'_{j_3}|,\dots, |\rho_{j_{r-1}}|=|\rho'_{j_r}|,\ |\rho_{j_r}|=|\alpha_{k'}|\,.
\end{equation}

Clearly, the parameters of the block $P_\nu$ we are aiming at have absolute values $|\alpha_k|, |\rho_{j_2}|,$ $\dots, |\rho_{j_r}|$ since these are the moduli of the eigenvalues of $L+\phi$ we want to realise within this block. It remains to carefully choose signs such that the derived parameters $\mu(P_\nu)$ and $\gamma(P_\nu)$ coincide up to sign with the given parameters $(\alpha_k,\beta_{j_2},\beta_{j_2},\dots, \beta_{j_r},\beta_{j_r},\alpha_{k'})$ and $(\gamma_{j_2},\dots, \gamma_{j_r})$, respectively.

We define the block $P_\nu=(\rho_1^\nu,\dots,\rho_r^\nu)$ of Type III by the following induction. We put $\rho_1^\nu:= \alpha_k$. Now suppose that $\rho^\nu_{l-1}$ is already defined for some $l$ with $2\le l\le r$ and suppose that  $\rho^\nu_{l-1}$ has the property $|\rho^\nu_{l-1}|=|\rho_{j_{l-1}}|$ if $l\ge 3$. By (\ref{Ear}), $\rho^\nu_{l-1}=\pm\rho'_{j_l}=\pm \beta_{j_l}\pm\gamma_{j_l}$ for a suitable choice of signs. Hence we can choose $\beta^\nu_l$, $\gamma^\nu_l$ such that $\rho^\nu_{l-1}=\beta_l^\nu+\gamma^\nu_l$ and $|\beta_l^\nu|=|\beta_{j_l}|$, $|\gamma^\nu_l|=|\gamma_{j_l}|$. Now we put $\rho^\nu_l:=\beta^\nu_l-\gamma_l^\nu$. Then we have 
$|\rho_l^\nu|=|\rho_{j_l}|$ since $\pm\rho'_{j_l}=\beta^\nu_l+\gamma_l^\nu$. Furthermore, $|\rho^\nu_r|=|\rho_{j_r}|=|\alpha_{k'}|$.

Since the vertices of $t_\nu$ are pairwise distinct, $|\rho_i^\nu|\not=|\rho_j^\nu|$ for $i\not=j$. Hence $P_\nu$ is a block. Its derived parameters are 
\begin{eqnarray}\label{Ederp1} 
\mu(P_\nu)=(\alpha_k,\beta_2^\nu,\beta_2^\nu,\dots,\beta_r^\nu,\beta_r^\nu,\pm\alpha_{k'}),\quad
\gamma(P_\nu)=(\gamma_2^\nu,\dots,\gamma_{r}^\nu ),
\end{eqnarray}
where $|\beta_\iota^\nu|=|\beta_{j_\iota}|$ and $|\gamma_\iota^\nu|=|\gamma_{j_\iota}|$ for $\iota=2,\dots,r$.

Let $j_1,\dots,j_s$ be those edges that do not belong to one ot the paths $t_1,\dots,t_{p/2}$. For $\iota=1,\dots,s$ we define a block $P_{p/2+\iota}=\displaystyle{\rho_\iota \choose \hat\rho_\iota}$ of Type II by  $\rho_\iota:=\rho_{j_\iota}$ and  $\hat \rho_\iota:=\kappa_\iota\beta_{j_\iota}-\kappa'_\iota\gamma_{j_\iota}$, if $\rho_{j_\iota}=\kappa_\iota\beta_{j_\iota}+\kappa'_\iota\gamma_{j_\iota}$ for $\kappa_\iota,\kappa_\iota'\in\{1,-1\}$. Then 
\begin{equation}\label{Ederp2} 
\mu(P_{p/2+\iota})=(\kappa_\iota\beta_{j_\iota},\kappa_\iota\beta_{j_\iota}),\quad \phi(P_{p/2+\iota})(z_1,z_2)=\kappa_j'(-z_2,z_1). \end{equation} 

If we now define $\cP:= (P_1|\dots|P_{p/2+s})$, then, by (\ref{Ederp1}) and (\ref{Ederp2}),  $(W,\theta,L,\phi)$ is isometrically isomorphic to
$(\CC^{d(\mathcal P)}, \mathop{\rm conj}, L(\cP),\phi(\cP))$. By construction, $\rho(\cP)$ contains every $\rho\in\sigma$, $\rho\not=0$, with multiplicity equal to the charged in-degree of $|\rho|$ and $0$ with multiplicity equal to twice the charged in-degree of $0$. By (\ref{Ein}), this is equal to the multiplicity of $\rho$ in $\sigma$.

It remains to discuss the case of odd $p$. Instead of $W$ we will work with the space $\tilde W:=W\oplus \CC$ with the trivial actions of $L$ and
$\phi$ on $\CC$ and consider the $(L,\phi)$-special subspace $\tilde V:=V\oplus i\RR$ of $\tilde W$. In particular, we replace $p$ by $\tilde p:=p+1$, $\sigma$ by $\tilde\sigma:=\sigma\uplus\{0\ms$,
and we have $\alpha_{p+1}=0$. We can now apply the above proof to $(\tilde W,\tilde V)$; the non-triviality of $\ker L\cap \ker\phi$ does not disturb that. Note that the vertex $0$ of the corresponding charged directed graph $G$ is always charged (positively or negatively). Therefore, after all the necessary orientation modifications, we get a particular charged path, say $t_1$,
starting or ending at the vertex $0$. Doing an additional orientation modification along that path, if necessary, we may assume
that it ends at $0$. The corresponding `block' $P_1$ (which looks like a block of Type III except for $\rho^1_r=0$) is isomorphic
to the direct sum of a block $P_1'$ of Type IV (having the same parameters as $P_1$) and a block $Q$ of Type I. All other blocks
of the resulting constellation $\cP=(Q|P_1'|P_2|\dots|P_{\tilde p/2+s})$ are blocks of Type~III and Type II, respectively.
We obtain $(\tilde W,\theta,L,\phi)\cong(\CC^{d(\mathcal P)}, \mathop{\rm conj}, L(\cP),\phi(\cP))$ and $\tilde\sigma=\rho(\cP)$.
Let $\cP'=(P_1'|P_2|\dots|P_{\tilde p/2+s})$. 
Then $(W,\theta,L,\phi)\cong(\CC^{d(\mathcal P')}, \mathop{\rm conj}, L(\cP'),\phi(\cP'))$ and $\sigma=\rho(\cP')$.
\qed

We will use special constellations in order to ensure the existence of certain subspaces $V\subset W$ invariant under operators of the form $A=\exp(2\pi(L+\phi))$ with prescribed $\spec(A|_V)$.
This motivates the following 

\begin{de}
For a special constellation $\cP$ we define a multiset $\nu(\cP)$ 
by
$$ \nu(\cP):=\{ e^{2\pi i \rho'}\mid \rho'\in \rho(\cP)\ms\ .$$
 A special constellation is called minimal if the one-parameter group
\begin{equation}\label{nase}\{\exp(t(L(\cP)+\phi(\cP))\mid t\in\RR\} 
\end{equation}
is contained in the closure of the group generated by 
$\exp(2\pi(L(\cP)+\phi(\cP))$. 
\end{de}

The parameters $\rho,\rho',\rho_i$ describing the blocks contained in a special constellation $\cP$ form a
vector $\tilde\rho\in\RR^{d'}$ for some $d'\le d(\cP)$. Then $\cP$ is minimal if and only if  
\begin{equation}\label{stuber}
\langle \tilde\rho,\ZZ^{d'}\rangle \cap \ZZ=\{0\}\ .
\end{equation}
Indeed, $\RR^{d'}$ parametrises a torus $T^{d'}\subset \U(d(\cP))$ via 
$\tilde\rho\mapsto A_{\tilde\rho}:=\exp(2\pi(L(\cP)+\phi(\cP))$. The character group of that torus is naturally isomorphic
to $\ZZ^{d'}$. A character corresponding to $\xi\in\ZZ^{d'}$ vanishes on $A_{\tilde\rho}$ if and only if 
$\langle \tilde\rho,\xi\rangle \in\ZZ$, while it vanishes on the full one-parameter group (\ref{nase}) if and only if $\langle \tilde\rho,\xi\rangle=0$.

Condition (\ref{stuber}) implies in particular that all non-zero parameters of a minimal special constellation $\cP$ and all non-zero coordinates of $\mu(\cP)$ are irrational.

\begin{de}\label{DPadm} Let $\cP$ be a special constellation of dimension $n=d(\cP)$. Let $\mu(\cP)=(\mu_1,\dots,\mu_n)$. An $n$-tuple $\uk=(k_1,\dots,k_n)\in\ZZ^n$ is called $\cP$-admissible if the following conditions hold:
\begin{enumerate}
\item $L_\uk$ commutes with $\phi(\cP)$.
\item $\mu_i+k_i\ne 0$ for $i=1,\dots,n$.
\item $G_n(\CC^n)^{L(\cP)+\phi(\cP),i\rho(\cP)}$ contains a $\uk$-good subspace $V\subset \CC^n$.
\end{enumerate}
\end{de}

If $\cP$ is the trivial special constellation of dimension $n$, i.e. it consists of blocks of Type~I, only,
then $\cP$-admissibility coincides with $\RR$-admissibility. 

The definition of special constellations of dimension $n$ yields a sign vector $\kappa\in\{-1,1\}^n$ such that
\begin{equation}\label{anton}
\langle \kappa, \mu(\cP)\rangle = 0\ .
\end{equation}
If $\uk$ is a $\cP$-admissible $n$-tuple, then the argument of the proof of Proposition~\ref{charles} shows that there exists 
(a possibly different) $\kappa\in\{-1,1\}^n$ such that
\begin{equation}\label{puenktchen}
\langle \kappa, \uk\rangle = 0\ .
\end{equation}
However, we have more:

\begin{pr}\label{fourier}
If $\uk$ is a $\cP$-admissible $n$-tuple, then there exists
$\kappa\in\{-1,1\}^n$ satisfying {\rm (\ref{anton})} and {\rm (\ref{puenktchen})}.
\end{pr}
\proof We argue similarly as in the proof of Proposition~\ref{charles}.
With the real $n$-dimensional subspace $V=V_C\subset \CC^n$ we associate a function of two variables
$$f_C(z,s):=\det\Im (z^{\uk}\exp(s L(\cP)) C)=\sum_\kappa d_\kappa 
z^{\langle\kappa,\uk\rangle} e^{is\langle\kappa,\mu(\cP)\rangle}\,,$$
where the summation runs over all $\kappa\in\{1,-1\}^n$.
Now suppose that $V$ is $(L(\cP),\phi(\cP))$-special. We claim that then $f_C$ is constant with respect to $s$. Indeed,
the invariance of $V$ with respect to $L(\cP)+\phi(\cP)$ implies that 
\begin{equation}\label{furunkel}
\exp(s L(\cP)) C = \exp(-s \phi(\cP)) CA(s)
\end{equation}
for some continuous family $s\mapsto A(s)\in \GL(n,\RR)$. Since $A(0)=\id$, we have $\det A(s)>0$. Let $C_{\Bbb R}$ be the real linear map from $\RR^n$ to $\CC^n$ induced by $C$ (or, equivalently, the $(2n\times n)$-matrix $(\Re C, \Im C))$.
Equation (\ref{furunkel}) implies that $C_{\Bbb R}^\top C_{\Bbb R}=A(s)^\top C_{\Bbb R}^\top C_{\Bbb R}A(s)$.
The map $C_{\Bbb R}$ is injective.
Taking determinants
eventually implies
that $A(s)\in \SL(n,\RR)$. Since $\exp(-s \phi(\cP))\in\SO(n)$ commutes with $z^\uk$ 
the claim follows. Thus
$$f_C(z,s)=\sum_{\langle \kappa,\mu(\cP)\rangle=0} d_\kappa 
z^{\langle\kappa,\uk\rangle} \ .$$
The constant term in the Fourier series of $f_C$ (considered as a function on $S^1$) is equal to 
$$ \sum_{\langle \kappa,\mu(\cP)\rangle= \langle\kappa,\uk\rangle=0} d_\kappa\ .$$
Lemma \ref{Ldet} implies that it is non-zero whenever $V$ is $\uk$-good. The proposition now follows.
\qed

\begin{pr}\label{Pktyp}
Let $\cP=(0|\dots|0|P_1|\dots|P_{d'})$ be a special constellation consisting of $d$ blocks of Type {\rm I} and blocks $P_1,\dots,P_{d'}$ of Type {\rm II.a}. Then $\uk =(\uk_0,\uk')\in \ZZ^d\times \ZZ^{2d'}$ is $\cP$-admissible if and only if $\uk_0$ is $\RR$-admissible and  $\uk'=(k_1',k_1',\dots,k'_{d'},k'_{d'})$.
\end{pr}
\proof Let $\uk =(\uk_0,\uk')\in \ZZ^d\times \ZZ^{2d'}$ be $\cP$-admissible. By assumption there exists a $\uk$-good subspace $V\subset \CC^{d+2d'}$ that is contained in the kernel of $L(\cP)+\phi(\cP)$.
The kernel of $L(\cP)+\phi(\cP)$ equals $\CC^d\oplus V'$, where $V'\subset \CC^{2d'}$ and $\dim_{\Bbb R}V'=2d'$.    The orthogonal projection $p'(V)$ of $V$ to $\CC^{2d'}$ is contained in $V'$. On the other hand, the projection of a good subspace to $\CC^{2d'}$ has real dimension at least $2d'$. Thus $p(V)=V'$ has dimension $2d'$. Consequently, $V_0:=V\cap \CC^d\subset \CC^d$ is $\uk_0$-good. Since $L(\cP)|_{{\Bbb C}^d}=0$, the $d$-tuple $\uk_0$ has only non-vanishing components. Hence $\uk_0$ is $\RR$-admissible. Since $L_\uk$ commutes with $\phi(\cP)=0_d\oplus \phi_\gamma$ for some $\gamma\in(\RR^*)^{d'}$, we have $\uk'=(k_1',k_1',\dots,k'_{d'},k'_{d'})$. The converse direction is easy to verify.
\qed

The following fact will become important in Section \ref{last}.

\begin{pr}\label{lagrange}
Let $\uk$ be a $\cP$-admissible
$n$-tuple. We define a real symplectic form $\omega:=\langle (L_\uk+L(\cP))(\,\cdot\,),\cdot \rangle$ on $W:=\CC^n$, where $\ip$ is the standard Euclidean inner product on $\CC^n$. If a real subspace $V\subset W$ is $(L(\cP),\phi(\cP))$-special and $\uk$-good, then $\omega|_{V\times V}\ne 0$, i.e. $V$ is not Lagrangian.
\end{pr}

\proof By  Def.~\ref{DPadm}, 2., $L:=L_\uk+L(\cP)$ is bijective. We consider the polar decomposition of this complex linear operator $L= J |L|$, 
where $|L|:=\sqrt{L^*L}=\sqrt{-L^2}$. Then $J^2=-\id$, thus $J$ defines a new complex structure on $W$. Moreover,
$$h(v,w):=\langle |L|v,w\rangle +i\langle |L|v,Jw\rangle
=\langle |L|v,w\rangle - i\omega(v,w),\quad v,w\in\CC^n,$$
is a positive definite Hermitian form on $(W,J)$. We denote the corresponding unitary group by 
$\U(h)$ and its subgroup leaving $\RR^n\subset W$ invariant by $\grO(h)$. Let  $\fu(h)$, $\fo(h)$ be the corresponding Lie algebras. Then $L_\uk, L(\cP)\in\fu(h)$, $\phi(\cP)\in\fo(h)$. Note that $\U(h)\subset \Sp(W,\omega)$.

Let $\cL\subset G_n(W)$ be the Grassmannian of Lagrangian subspaces of $(W,\omega)$. 
The map $\U(h)\ni A \mapsto A(\RR^n)\in\cL$
yields an isomorphism $\U(h)/\grO(h)\cong\cL$. Moreover, the square of the determinant $\det_J$ of endomorphisms of the complex
vector space $(W,J)$
\begin{equation}\label{paul} \U(h)/\grO(h)\ni [A] \mapsto ({\det} _J A)^2\in S^1\subset\CC 
\end{equation}
induces an isomorphism of fundamental groups
\begin{equation}\label{ernst} \pi_1({\mathcal L},\RR^n)\stackrel{\sim}\longrightarrow\pi_1(S^1,1)\cong\ZZ\ .
\end{equation}
These facts can be found in \cite{GS77} for example.

Now let $V\subset W$ be an  $(L(\cP),\phi(\cP))$-special Lagrangian subspace of $(W,\omega)$. 
The closed curve $[0,1]\ni t\mapsto \exp(2\pi tL_\uk)V=\exp(2\pi t(L+\phi(\cP)))V$ defines an element $c_\uk\in
\pi_1({\mathcal L},V)\cong\ZZ$. According to (\ref{paul}) and (\ref{ernst}) the corresponding integer is given
by the winding number of the closed curve
$$
[0,1]\ni t\mapsto ({\det}_J\exp(2\pi t(L+\phi(\cP))))^2=e^{4\pi t\,{\tr}_J(L+\phi(\cP))}=e^{4\pi t\,{\tr}_J L}=e^{4\pi it\,{\tr}|L|}\in S^1 \ .$$
Here $\tr_J$ is the trace on $(W,J)$, while $\tr$ denotes the usual trace on $\CC^n$. Since $\tr |L|$ is positive we conclude that $c_\uk\ne 0$.

We consider the open subset $\cL_0:=\{V'\in\cL\mid V'\cap\RR^n=\{0\}\}\subset\cL$. The map $\Phi\mapsto$ graph$(\Phi)$ is a diffeomorphism from the real vector space
$$ \{\Phi\in\Hom(i\RR^n,\RR^n)\mid L\circ\Phi\in\End(i\RR^n)\mbox{ is symmetric w.r.t. }\ip|_{i\Bbb R^n} \}$$
to $\cL_0$.
Hence $\cL_0$ is contractible. The non-triviality of  $c_\uk\in \pi_1({\mathcal L},V)$ implies that
$$\{\exp(2\pi tL_\uk)V\mid t\in [0,1]\}\not\subset\cL_0\ ,$$ 
thus $V$ is not $\uk$-good.
\qed

\subsection{The main theorem}

For a polynomial $f$ in one variable we denote the multiset of those roots of $f$ lying on the unit circle
by $\nu_{c}(f)$.

\begin{theo}\label{DD}
Let $X$ be a Cahen-Wallach space of type $(p,q)$. Then $X$ admits a compact quotient
if and only if there exist
\begin{enumerate}
\item[(a)] a polynomial $f\in\ZZ[x]$ of degree $p+q$ of the form {\rm(\ref{otto})} having precisely $q$ roots on the unit circle (counted with multiplicity),
\item[(b)] a special constellation $\cP$ of dimension $q$ with $\nu(\cP)=\nu_c(f)$,
\item[(c)] a $\cP$-admissible $q$-tuple $\uk\in\ZZ^q$
\end{enumerate}
such that
\begin{equation}\label{auf} X\cong X_{p,q}\Big(\log |\nu_1|,\dots,\log |\nu_p|;\,2\pi(\mu(\cP)+\uk)\,\Big)\ ,
\end{equation}
where $\nu_1,\dots,\nu_p$ are the roots of $f$ of modulus different from $1$.

The assertion remains true if we require that $\nu_{c}(f)$ contains no roots of unity except $1$ and that $\cP$ is minimal.
\end{theo}

\begin{re}\label{Rri}{\rm 
For $q=0$ (real case) the theorem specialises to Thm.~\ref{BB}. Let us discuss the less obvious relation to 
Thm.~\ref{CC} (imaginary case). Let $X$ be given as in Thm.~\ref{CC}, i.e., $X\cong X_{0,q}(\uk_0,\mu')$, where $\uk_0\in (\ZZ_{\not=0})^d$ is $\RR$-admissible and $\mu'=(\mu_1,\mu_1,\dots,\mu_{d'},\mu_{d'})\in(\RR^*)^{2d'}$ for $q=d+2d'$. Then it is not hard to find data $f,\cP,\uk$ satisfying (a), (b), (c) of Thm.~\ref{DD} such that $X\cong X_{0,q}(2\pi(\mu(\cP)+\uk))$. We just put $f(x)=(x-1)^q$ and consider the special constellation $\cP=(0|\dots|0|P_1|\dots | P_{d'})$, where $P_j={0\choose {2\mu_j}}$, $j=1,\dots,d'$, is a block of Type~II.a.  Then $\mu(\cP)=(0,\mu')\in \RR^d\times(\RR^*)^{2d'}$ and $\nu(\cP)=\{1,\dots,1\ms=\nu_c(f)\in\RR^q$.  Furthermore, the
$q$-tuple $\uk:=(\uk_0,0)$ is $\cP$-admissible by Prop.~\ref{Pktyp}.
Conversely, suppose that $X$ is given as in Thm.~\ref{DD} by the data $f,\cP,\uk$. If all roots of a polynomial of the form (\ref{otto}) have modulus $1$, then they are roots of unity.
Thus Thm.~\ref{DD} tells us that in the imaginary case it suffices to consider minimal special constellations $\cP$ for  which $\rho(\cP)$ is supported at ${0}$, see (\ref{stuber}). Such constellations consist of blocks of Type I and Type II.a, only. Now Prop.~\ref{Pktyp} yields 
$X=X_{0,q}(2\pi(\mu(\cP)+\uk))\cong X_{0,q}(\uk_0,\mu_1+k_1',\mu_1+k_1',\dots,\mu_{d'}+k'_{d'},\mu_{d'}+k'_{d'})$, where $\uk_0$ is $\RR$-admissible. Consequently, in the imaginary case, Thm.~\ref{DD} specialises to Thm.~\ref{CC}.
}\end{re}
\begin{re}\label{Reasy}{\rm
For the convenience of the reader, let us recall two elementary properties of polynomials of the form (\ref{otto}). We will frequently use these properties in the remainder of this section. Let $f\in\ZZ[x]$ be an irreducible polynomial of the form~(\ref{otto}), $f(x)\not=x\pm1$, and let $f$ have a root on the unit circle. Then the following holds.
\begin{itemize}
\item[(a)] The roots of $f$ come in pairs $\nu$, $\nu^{-1}$. In particular, the degree of $f$ is even.
\item[(b)] If one of the roots of $f$ is a root of unity, then all roots are roots of unity. 
\end{itemize}
}\end{re}
{\sl Proof of Thm.~\ref{DD}. }
We assume that we are given data $f,\cP,\uk$ satisfying (a), (b), (c). 
We decompose $f=f_0f_1$, $f_i\in\ZZ[x]$, where $f_0$ has no roots on the unit circle and
all $\QQ$-irreducible factors of $f_1$ do have such a root. 
We index the roots of $f$ of modulus different from $1$
in a way such that $$\Im(\nu_{2k-1})>0,\quad 
\bar\nu_{2k-1}=\nu_{2k},\  k=1,\dots,s\quad\mbox{and}\quad \nu_l\in\RR^* \mbox{ for }l=2s+1,\dots,p $$ 
and such that
$\nu_1,\dots,\nu_{2s_0}$, $\nu_{2s+1},\dots,\nu_{2s+r_0}$, $2s_0+r_0=\deg(f_0)$, are the roots of $f_0$.
We consider the Cahen-Wallach space 
$$X_{p,q}\Big(\log |\nu_1|,\dots,\log |\nu_p|;\,2\pi(\mu(\cP)+\uk\,)\Big)\ .$$ 
We want to show that it admits a compact quotient by constructing data $V,\Lambda,t_0,\ph_0,h_0$ as required by Proposition~\ref{wulf}. 

Again we follow the conventions of Example \ref{Exosc2}. There is a splitting
$$ \fa=\CC^p\oplus \CC^q=:\fa_{\Bbb{R}}\oplus\fa_I,\qquad  L=L_{\Bbb R}\oplus L_I.$$

Note that $L_I=2\pi(L(\cP)+L_\uk)$. We set 
$$A_I:=\exp(L_I+2\pi\phi(\cP))=\exp(2\pi(L(\cP)+\phi(\cP)))\ .$$
We set $B:=L(\cP)+\phi(\cP)\in\fsp(\fa_I, \omega|_{\fa_I})$. Then $G_q(\fa_I)^{B,i\rho(\cP)}\subset G_q(\fa_I)^{A_I,\nu_c(f)}$.

Since $\uk$ is $\cP$-admissible the space 
$G_q(\fa_I)^{B,i\rho(\cP)}$ contains a $\uk$-good subspace.
Taking Corollary~\ref{open} into account it follows that the space $G_q^{\;\uk}(\fa_I)^{B,i\rho(\cP)}$ of $\uk$-good subspaces belonging to $G_q(\fa_I)^{B,i\rho(\cP)}$
is a non-empty open subset of $G_q(\fa_I)^{B,i\rho(\cP)}$. Lemma~\ref{lehm} implies that the intersection of $G_q^{\;\uk}(\fa_I)^{B,i\rho(\cP)}$ with $G_q(\fa_I)^{B,i\rho(\cP)}_\mathrm{reg}$ is non-empty.
We choose an element $V_I$ of that intersection. We decompose the vector space $V_I$
into $V_I^0:=\ker B\cap V_I$ and a $B$-invariant complement $V_I^1$ of $V_I^0$. Then the restriction of the symplectic form $\omega$ to  $V_I^1$ is
non-degenerate.

Now we consider $\fa_{\Bbb{R}}=\CC^p$. We define a complex linear operator $\ph_{\Bbb{R}}$ on $\fa_{\Bbb{R}}$ by a block diagonal matrix as in the proof of
Thm.~\ref{BB}, (\ref{bella}) (with $n$ replaced by $p$). Let $V_{\Bbb{R}}^0\subset\fa_{\Bbb{R}}$ be the real subspace spanned by
the vectors $(1+i)e_l$, $l=1,\dots,2s_0,2s+1,\dots,2s+r_0$. 
The space $V_{\Bbb{R}}^0$ is totally isotropic, invariant under $L_{\Bbb{R}}$ and $\ph_{\Bbb{R}}$, and satisfies 
$V_{\Bbb{R}}^0\cap\fa_+=\{0\}$.  

The roots $\nu_l$ with the remaining indices (i.e., the roots of $f_1$ outside the unit circle) come in pairs $\nu_l$, $\nu_l^{-1}$. We can therefore group these indices into blocks of the form $(2k-1,2k,2k'-1,2k')$ and $(l,l')$ such that $\nu_{2k-1}^{-1}=\nu_{2k'}$,  
$\overline{\nu_{2k-1}}=\nu_{2k}$, $\overline{\nu_{2k'-1}}=\nu_{2k'}$, and $\nu_{l}^{-1}=\nu_{l'}$, $\nu_l\in\RR$.
Let $V_{\Bbb{R}}^1\subset\fa_{\Bbb{R}}$ be the real subspace spanned by
the vectors
$$ (1+i)e_{2k-1},\ (1+i)e_{2k},\ (1-i)e_{2k-1}+(1+i)e_{2k'},\ (1-i)e_{2k}-(1+i)e_{2k'-1} $$
and
$$ (1+i)e_l,\ (1-i)e_l+(1+i)e_{l'}, $$
where the indices $k$, $l$ run over all blocks.
The space $V_{\Bbb{R}}^1$ is symplectic, invariant under $L_{\Bbb{R}}$ and $\ph_{\Bbb{R}}$, and satisfies 
$V_{\Bbb{R}}^1\cap\fa_+=\{0\}$.

We now set 
\begin{equation}\label{didi}
V:= V_{\Bbb{R}}^0\oplus V_{\Bbb{R}}^1 \oplus V_I, \ t_0:=1,\ \ph_0:=\ph_{\Bbb{R}}\oplus\exp({2\pi\phi(\cP)}),\ h_0:=0
\end{equation}
and $A:=e^{ L}\ph_0=e^{ L_{\Bbb{R}}}\ph_{\Bbb{R}}\oplus A_I$. It remains to construct an ($\id_\fz\oplus A$)-stable lattice
$\Lambda$ in the nilpotent group $\fz\oplus V$. 

By construction, $A$ is semisimple. The characteristic polynomial of $A|_{V_{\Bbb{R}}^0}$ is equal to $f_0$. By Lemma~\ref{gigi}
there exists an $A$-stable lattice $\Lambda_{\Bbb{R}}^0\subset V_{\Bbb{R}}^0$.
We set $V^1:=V_{\Bbb{R}}^1 \oplus V_I^1$. Then $V^1$ is $A$-invariant, symplectic and $A|_{V^1}\in\Sp(V^1,\omega|_{V^1})$ has characteristic
polynomial $f_1/(x-1)^{q_0}$, $q_0=\dim V_I^0$. By Lemma~\ref{nono}
there exists an $A$-stable lattice $\Lambda^1\subset V^1$ satisfying $\omega(\Lambda^1\times\Lambda^1)\subset\ZZ$.
We choose a lattice $\Lambda_I^0\subset V_I^0$ such that $\omega(\Lambda_I^0\times\Lambda_I^0)\subset\ZZ$. It is $A$-stable 
since $A$ acts trivially on $V_I^0$.

Since  $V_{\Bbb{R}}^0$ is isotropic and the spaces  $V_{\Bbb{R}}^0$, $V^1$ and $V_I^0$ are pairwise orthogonal with respect to
$\omega$ the $A$-stable lattice $\Lambda_0:= \Lambda_{\Bbb{R}}^0\oplus\Lambda^1\oplus \Lambda_I^0\subset
V_{\Bbb{R}}^0\oplus V^1 \oplus V_I^0=V$ satisfies $\omega(\Lambda_0\times\Lambda_0)\subset\ZZ$. Therefore
$$ \Lambda:= \mbox{$\frac{1}{2}$}\ZZ\oplus \Lambda_0\subset \fz\oplus V$$
is an ($\id_\fz\oplus A$)-stable lattice as desired.

For the opposite direction we assume that a Cahen-Wallach space $X$ of type $(p,q)$ admits a compact quotient, i.e. there are
objects $V,\Lambda,t_0,\ph_0,h_0$ satisfying Conditions $(a)$ and $(b)$ of Proposition~\ref{wulf}. We want to find
a polynomial $f\in\ZZ[x]$ of the form (\ref{otto}) having precisely $q$ roots on the unit circle, none of them a root of unity different from $1$, a minimal special constellation $\cP$ such that $\nu(\cP)=\nu_c(f)$, and a $\cP$-admissible $q$-tuple $\uk$
such that (\ref{auf}) holds.

There is a splitting $\fa=\fa_{\Bbb{R}}\oplus\fa_I$ such that $L|_{\fa_\Bbb{R}}$ has real and $L|_{\fa_I}$ has imaginary eigenvalues. By assumption, the operator $A:=e^{t_0L}\ph_0$ leaves invariant $V$. We decompose $V=V_{\Bbb{R}}\oplus V_I$ into $A$-invariant subspaces such that the eigenvalues of $A$ have modulus different from one on $V_{\Bbb{R}}$ and equal to one on $V_{I}$. Obviously, $V_*:=V\cap \fa_*$, $*=\RR, I$.

The closure $T$ in $\grO(\fa_I)$ of the group generated by
$A_I:=A|_{\fa_I}$ 
has finitely many connected components since it is compact. By replacing $t_0$ by $kt_0$ 
and $\ph_0$ by $\ph_0^k$ (and $h_0$ by $h_0^k$) for a certain $k\in\NN$, if necessary, we can assume that $T$ is connected. In particular, no root
of unity except $1$ is an eigenvalue of $A_I$ on $\fa_I$. In the following, it will turn out to be useful to work with $A^2$ instead of $A$.

Let $f$ be the characteristic polynomial of $A^2|_V$. Then no root of $f$ is a root of unity different from $1$. Moreover, 
$\tilde f:=(x-1)f$ is the characteristic polynomial of $\Ad(h_0t_0\ph_0)^2|_{\fz\oplus V}$. The latter operator stabilises the $\ZZ$-module
$\tilde \Lambda$ generated by $\Lambda\subset  \fz\oplus V$. Since $\Lambda$ is a lattice in the 1-connected nilpotent Lie group
$\fz\oplus V$, the  $\ZZ$-module
$\tilde \Lambda$ is a lattice in $\fz\oplus V$ considered as a vector space (see e.g. the remark after Thm.~2.12 in \cite{Ragh}). By Lemma~\ref{gigi} the polynomial $\tilde f$
is integral of the form (\ref{otto}). Hence the same is true for $f$.

The roots $\nu_1,\dots,\nu_r$ of $f$ that lie outside the unit circle are the eigenvalues of $A^2|_{V_{\Bbb{R}}}$. In particular, $r=p$.
As in the second part of the proof of Thm.~\ref{BB} we now see that the eigenvalues of $L|_{\fa_{\Bbb{R}}}$ are precisely
\begin{equation}\label{tauto}  
\frac{\log |\nu_1|}{2t_0},\dots,\frac{\log |\nu_p|}{2t_0},-\frac{\log |\nu_1|}{2t_0},\dots,-\frac{\log |\nu_p|}{2t_0}\ .
\end{equation}

We set $T^\theta:=\{\theta_\fa \psi \theta_\fa\mid \psi\in T\}$. Elements of $T^\theta$ commute with $L_I:=L|_{\fa_I}$ as well as with elements of $T$.
Let $\ft$ and $\tilde\ft$ be the Lie algebras of the tori $T$ and $\tilde T:=TT^\theta$, respectively.
Conjugation by $\theta_\fa$ defines an involution on $\tilde\ft$. Let $\tilde\ft=\ft_+\oplus\ft_-$ be the corresponding eigenspace decomposition. Since $A_I\in T\subset \tilde T$ there exist elements $L_0\in\ft_-$, $\phi\in\ft_+$ such that $L_0+\phi\in\ft$ and 
$A_I=\exp(\pi(L_0+\phi))$. 
We have $T(V_I)\subset V_I$. Hence, $V_I$  is an $(L_0,\phi)$-special subspace of $\fa_I$. 
By Proposition~\ref{veryspecial} we can identify $\fa_I$ with $\CC^q$ in a way such that
$$  L_0=L(\cP)\ \mbox{ and }\ \phi=\phi(\cP) $$
for some special  constellation $\cP$ of dimension $q$ satisfying
$$\nu(\cP)=\spec(A_I^2|_{V_I})=\nu_c(f)\ .$$ 
The group generated by $A_I^2$ is dense in $T$. Therefore $\cP$ is minimal.

We have 
$$ \exp(2\pi L(\cP))=A_I \left( \theta_\fa A_I^{-1} \theta_\fa \right)= \exp(2t_0 L_I)\ .$$
Since $L_I$ and $L(\cP)$ are commuting semisimple operators anticommuting with $\theta_\fa$ this implies that we can adapt
the basis of $\fa_I$ further such that
\begin{equation}\label{schluck}  \frac{t_0}{\pi} L_I = L(\cP)+ L_\uk  
\end{equation}
for some $\uk\in\ZZ^q$. The operators $L(\cP)$ and $L_I$ commute with $\phi(\cP)$. Hence so does $L_\uk$.
Since
$$ e^{t L_\uk} V_I\cap\fa_+=e^{t L_\uk}\exp\left( t(L(\cP)+\phi(\cP))\right) V_I\cap\fa_+=e^{t\phi(\cP)}\left( e^{\frac{tt_0}{\pi}L_I}V_I\cap\fa_+\right) =\{0\}$$
for all $t\in\RR$ the $q$-tuple $\uk$ is $\cP$-admissible. Now (\ref{tauto}) combined with (\ref{schluck}) implies (\ref{auf}).
\qed 

\begin{co}\label{spurnull3}
Assume that $X_{p,q}(\lambda,\mu)$ admits a compact quotient. Then there are choices of signs such that
$$
\displaystyle \sum_{i=1}^{p} \pm\lambda_i = 0\quad\mbox{ and }\quad \sum_{j=1}^{q} \pm\mu_j = 0\ .
$$
In particular, spaces of type $(1,q)$ or $(p,1)$ do not admit compact quotients.
\end{co}

\proof The equation for $\lambda$ follows as in the proof of Cor.~\ref{spurnull}, whereas the equation for $\mu$ is a consequence
of Prop.~\ref{fourier}.
\qed
\begin{re}\label{Rpol}{\rm Let $f\in\ZZ[x]$ be a polynomial of the form (\ref{otto}) which has exactly $q\not=1$ roots on the unit circle. Suppose that the multiplicity of $-1$ in $\nu_{c}(f)$ is even. Then there exists a special constellation $\cP$ and a $q$-tuple $\uk\in \ZZ^q$ such that the data $f,\cP,\uk$ satisfy the Conditions (b) and (c) of Thm.~\ref{DD}. 

This is clear for $q=0$. Let us verify the assertion for $q>1$. It suffices to find a special constellation $\cP$ such that $\nu(\cP)=\nu_c(f)$ and $\mu(\cP)$ has non-vanishing components. Then $\uk=0$ is $\cP$-admissible since according to Thm.~\ref{veryspecial} there exists an $(L,\phi)$-special subspace $V$ such that $\spec((L+\phi)|_V)=i\rho(\cP)$. Obviously, $V$ is $\uk$-good and Conditions 1.\ and 2.\ of Def.~\ref{DPadm} are satisfied. 
Let us show that such a special constellation does exist. We put $\nu_c(f)':=\nu_c(f)\setminus\{1\ms$ if $q$ is odd and $\nu_c(f)':=\nu_c(f)$ if $q$ is even. Then  $\nu_c(f)'=\{e^{\pm 2\pi i \rho_1},\dots,e^{\pm 2\pi i \rho_r}\ms$ for suitable $\rho_1,\dots, \rho_r\in \RR^*$, where $r=[q/2]$. We define $P_1={\rho_1\choose 0}$ if $q$ is even and  $P_1=(\rho_1,0)$ if $q$ is odd. Moreover, $P_j:={\rho_j\choose 0}$ for $j=2,\dots,d$. Now we set $\cP=(P_1|\dots|P_r)$.  Then $\nu(\cP)=\nu_c(f)$ and $\mu(\cP)$ has only non-vanishing components. 
}\end{re}

Regarding the construction of Cahen-Wallach spaces of type $(2,q)$, of particular interest are  the irreducible polynomials of degree $2k$ in $\ZZ[x]$ that have exactly $2k-2$ roots on the unit circle. Let us denote the set of these polynomials by $F_{2k}$.  Each $f\in F_{2k}$ is reciprocal and the two roots of modulus different from 1 are real, see Rmk.~\ref{Reasy}. If $r$ is the root of $f\in F_{2k}$ of modulus greater than 1, then $|r|$ is called {\it Salem number of degree $2k$}. This is not the original definition but note that there also exists a polynomial in $F_{2k}$ that has $|r|$ as a zero. In the following examples we will be especially interested in $F_4$ and $F_6$. It is not hard to prove that 
$$ F_4=\{ x^4-ax^3+bx-ax+1\mid 2|a|>|b+2|, \ b\not=2,\,b\not=\pm a+1 \},$$
see also \cite{B1}. Salem numbers of degree 6 are studied, e.g., in \cite{B2}. In the Supplement to \cite{B2} one can find tables listing examples of polynomials contained in $F_6$. For instance, $f(x)=x^6-x^4-x^3-x^2+1$ is in $F_6$.

\begin{ex}[compositions]\label{Excomp}{\rm 
Suppose we are given two Cahen-Wallach spaces  $X_1$, $X_2$ of type $(p_1,q_1)$ and $(p_2,q_2)$ admitting compact quotients. According to Thm.~\ref{DD} these spaces are given by data $f_1,\cP_1,\uk_1$ and $f_2,\cP_2,\uk_2$, respectively. Then $f:=f_1f_2$ together with the direct sums $\cP:=(\cP_1|\cP_2)$ and $\uk:=(\uk_1,\uk_2)$ satisfy (a), (b), (c) of Thm.~\ref{DD}. Hence, the data  $f,\cP,\uk$ define a new Cahen-Wallach space $X$ of type $(p_1+p_2,q_1+q_2)$ which has compact quotients. 

We will say that $X$ is {\it composed} of $X_1$ and $X_2$ if there are data $f_i,\cP_i,\uk_i$, $i=1,2$, such that $X_i$ is given by $f_i,\cP_i,\uk_i$ and such that $X$ is isometric to the space constructed from $f_i,\cP_i,\uk_i$, $i=1,2$, in the above way.  

In particular, we can construct examples of Cahen-Wallach spaces of mixed type admitting compact quotients composing examples of purely real and imaginary type.  Let us formulate this in the notation of the Theorems~\ref{BB} and~\ref{CC}. Let $X_1$ be a Cahen-Wallach space of type $(p,0)$ and let $X_2$ be one of type $(0,q)$, both  with compact quotients.  Then $X_1\cong X_{p,0}(\log |\nu_1|,\dots,\log |\nu_p|)$ according to Thm.~\ref{BB} and $X_2\cong X_{0,q}(\uk_0,\mu')$, where $\uk_0\in (\ZZ_{\not=0})^d$ is $\RR$-admissible and $\mu'=(\mu_1,\mu_1,\dots,\mu_{d'},\mu_{d'})\in(\RR^*)^{2d'}$ for $q=d+2d'$ according to Thm.~\ref{CC}. Composing $X_1$ and $X_2$ yields the Cahen-Wallach space $X_{p,q}(\log |\nu_1|,\dots,\log |\nu_p|;2\pi\uk_0,2\pi \mu')$, see Rmk.~\ref{Rri}. 

Suppose that $X=X_{p,q}(\log |\nu_1|,\dots,\log |\nu_p|;\,2\pi(\mu(\cP)+\uk)\,)$ is given as in Thm.~\ref{DD}. We claim that $X$ is composed of spaces of real and imaginary type having compact quotients if $f_1:=f_{\Bbb R}:=(x-\nu_1)\dots(x-\nu_p)$ is in $\ZZ[x]$. Indeed, under this assumption $f$ decomposes as $f=f_1f_2$ with $f_2\in \ZZ[x]$. Furthermore, we have trivial decompositions  $\cP=(\cP_1|\cP_2)$ and $\uk=(\uk_1,\uk_2)$, where $\cP_1$ and $\uk_1$ are empty and $\cP_2=\cP$, $\uk_2=\uk$. If we now define $X_i$ by $f_i,\cP_i,\uk_i$, $i=1,2$, then $X_1$ has real type, $X_2$ has imaginary type and $X$ is composed of $X_1$ and $X_2$.
}\end{ex}

Now we want to show that there are Cahen-Wallach spaces of mixed type admitting compact quotients that are not composed of spaces of purely real and imaginary type having compact quotients. Thus Thm.~\ref{DD} really yields more Cahen-Wallach spaces than we can get by such compositions. Let us study the situation in small dimensions.
Recall that spaces of type $(1,q)$ or $(p,1)$ do not admit compact quotients. Thus the smallest type where we can find examples is $(p,q)=(2,2)$. However, let us start with type $(3,2)$, which is even easier.

\begin{ex}[type (3,2)]\label{Ex32}{\rm
Each space of type $(3,2)$ admitting a compact quotient is isometric to a composed one. Indeed, according to Example~\ref{Excomp} it suffices to show that for any polynomial $f$ that satisfies the assumptions of Thm.~\ref{DD},  the polynomial $f_{\Bbb R}=(x-\nu_1)(x-\nu_2)(x-\nu_3)$ is in $\ZZ[x]$. If the remaining roots $\nu_4$, $\nu_5$ of $f$, i.e., those on the unit circle, are in 
$\{1,-1\}$, then this is obviously true. If not, $f_{\Bbb R}\in \ZZ[x]$ follows from Rmk.~\ref{Reasy} (a). }\end{ex}

The next example will show that for given data $f, \cP, \uk$ satisfying (a), (b), (c) of Thm.~\ref{DD}, $X\cong X_{p,q}(\log |\nu_1|,\dots,\log |\nu_p|;\,2\pi(\mu(\cP)+\uk)\,)$ can be isometric to a composed example even if $f$ is irreducible. In particular, that $f_{\Bbb R}$ is in $\ZZ[x]$ is not a necessary condition for being isometric to a composed example. Suppose that $f_{\Bbb R}\not\in\ZZ[x]$ but $X$ is isometric to a composed example $X'$.  Let $f', \cP', k'$ be the data that are associated with $X'$ according to Example~\ref{Excomp}. Applying the construction in the proof of Thm.~\ref{DD} to the data $f, \cP, k$ and the data $f', \cP', k'$ then yield essentially different discrete cocompact subgroups on $X\cong X'$.
 
\begin{ex}[type (2,2)]\label{Ex22}{\rm
Each Cahen-Wallach space $X$ of type $(2,2)$ admitting a compact quotient is isometric to one of the composed ones constructed in Example~\ref{Excomp}. Indeed, $X\cong X_{2,2}(r,r;s,s)$ by Cor.~\ref{spurnull3}.
Let $\nu,\nu^{-1}$, $|\nu|\not=1$, be the roots 
of an irreducible quadratic polynomial $f_0\in\ZZ[x]$ of the form (\ref{otto}) and choose $\mu$ such that 
$(r,s)$ and $(\log|\nu|,2\pi \mu)$ are proportional. Then $X_0=X_{2,0}(\log |\nu|,-\log |\nu|)$ and $X_1=X_{0,2}(\mu,\mu)$ have compact quotients and are parametrised according to Thm.~\ref{BB} and Thm.~\ref{CC}, respectively. Hence $X\cong X_{2,2}(r,r;s,s)
\cong X_{2,2}(\log|\nu|,-\log |\nu|,2\pi \mu,2\pi \mu)$ arises by composing $X_0$ and $X_1$ as in  Example \ref{Excomp}. Nonetheless, if $f,\cP, \uk$ are the data describing $X$ according to Thm~\ref{DD}, then $f$ can be irreducible. Indeed, 
take $f\in F_4$, choose $\cP$ and $\uk$ according to Rmk.~\ref{Rpol} and consider the space $X$ defined by $f,\cP,\uk$. }\end{ex}

\begin{ex}[type (2,3)]\label{Ex23}{\rm
Here we will get for the first time examples that are not composed. Each $X$ of type $(2,3)$ that admits a compact quotient is isometric to $X_{2,3}(\log |\nu_1|,\log |\nu_2|;\,2\pi(\mu(\cP)+\uk\,))$, where $\nu_1,\nu_2, \cP$ and $\uk$ are as in Thm.~\ref{DD}. We may assume that $\cP$ is minimal and that $\nu_c(f)$ contains no roots of unity except 1. If $f_{\Bbb R}(x)=(x-\nu_1)(x-\nu_2)$ is in $\ZZ[x]$, then Example~\ref{Excomp} shows that $X$ is composed from $X_{2,0}(\log |\nu_1|,\log |\nu_2|)$ and $X_{0,3}(\uk)$.

Suppose now that $f_{\Bbb R}\not\in \ZZ[x]$. Then $f(x)=\hat f(x)(x- 1)$, where $\hat f$ is irreducible and has two conjugate roots $\nu_3=\bar \nu_4$ on the unit circle and two roots $\nu_1,\nu_2$ with $|\nu_1|, |\nu_2|\not=1$. Hence $\hat f\in F_4$.  We claim that $X$ is not composed provided the four exponentials conjecture is true. Let us first show that $\cP$ is of Type IV. Assume not, then $\cP=(0|P_1)$, where $P_1$ is of Type II. By Prop.~\ref{Pktyp}, $P_1$ cannot be of Type II.a since there do not exist $\RR$-admissible 1-tuples $k_1\in \ZZ_{\not=0}$. If $P_1$ is not of Type II.a, then $\rho(\cP)=\{0,\pm\rho\}$ for some $\rho\not=0$. By assumption, $G_3(\CC^3)^{L(\cP)+\phi(\cP),i\rho(\cP)}$ contains a $\uk$-good subspace $V\subset \CC^3$. Since $\rho\not=0$, we have $\dim _{\Bbb R}\ker(L(\cP)+\phi(\cP))|_V=1$. Hence the projection of the real subspace $V\subset\CC^3$ to the first component is one-dimensional. Since $V$ is $\uk$-good, this contradicts Lemma~\ref{Ldet}. Hence $\cP$ consists of a block $(\rho,0)$ of Type IV, which implies $X\cong X_{2,3}(\log |\nu_1|,-\log |\nu_1|;\,2\pi (\rho+2k, \rho/2+k,\rho/2+k)\,)$ for some $k\in\ZZ$. Assume that $X$ is composed. We have already seen that composed examples are isometric to $X_{2,3}(\log |\nu'|,-\log |\nu'|;\,2\pi \uk')$, where $\nu'$ is the root of a quadratic polynomial over $\ZZ$ and $\uk'$ is $\RR$-admissible. Hence there exists an integer $k'\not=0$ such that the vectors $(\log |\nu_1|, \rho/2 +k)$ and $(\log |\nu'|, k')$ are proportional. Recall that $\nu_1\in\RR$ and $e^{2\pi i\rho}$ are algebraic since they are roots of $f$. 
Hence we can apply the four exponentials conjecture in the form stated before Prop.~\ref{duli} to $(\lambda_{11},\lambda_{12})=(\log |\nu_1|,2\pi i(\rho/2 +k))$ and $(\lambda_{21},\lambda_{22})=(\log |\nu'|, 2\pi ik')$. Obviously, $\lambda_{i1}\in\RR$ and $\lambda_{i2}\in i\RR$ are independent over $\QQ$ for $i=1,2$. Moreover, $\lambda_{12}$ and $\lambda_{22}$ are independent over $\QQ$. Indeed, $\rho$ is irrational since $e^{2\pi i\rho}$ is a root of $f$, thus it is not a root of unity. If the conjecture is true, $(\log |\nu_1|, 2\pi i(\rho/2 +k))$ and $(\log |\nu'|, 2\pi ik')$ are linearly independent over $\CC$, which contradicts our assumption. Hence $X$ is not composed of spaces of real and and imaginary type having compact quotients.
}\end{ex}

\begin{ex}[type (2,4)]{\rm We want to show that there exist examples of Cahen-Wallach spaces of type (2,4) admitting compact quotients that are not composed of spaces of real and and imaginary type. Take a polynomial $f\in F_6$ and let $e^{\pm 2\pi i \rho_1}$ and $e^{\pm 2\pi i \rho_2}$ be its roots on the unit circle. Let $\cP$ consist of the block $(\rho_1,\rho_2)$ of Type III and put $\uk=0$. Then $\nu (\cP)=\nu_c(f)$ and the components of $\mu(\cP)= (\rho_1,\frac{\rho_1+\rho_2}2, \frac{\rho_1+\rho_2}2,\rho_2)$ do not vanish since $\rho_1\not=\pm\rho_2$. Following Rmk.~\ref{Rpol},  $f,\cP,\uk$ satisfy the assumptions of Thm.~\ref{DD}. Let $X=X_{2,4}(\log |\nu|,-\log |\nu|;2\pi\mu(\cP))$ be the corresponding Cahen-Wallach space. We claim that $X$ is not composed. Assume that $X$ is isometric to a composed example. Then $X\cong X_{2,4}(\log |\hat\nu|,-\log |\hat\nu|;2\pi(\hat \uk,\hat\mu))$, where $\hat\uk$ and $\hat \mu$ satisfy the assumptions in Thm.~\ref{CC}. In particular, we may assume  that $(\hat \uk,\hat\mu)$ and $\mu(\cP)$ are proportional vectors. Since $\rho_1\not=\pm\rho_2$,  $\hat\mu$ is the empty tuple and $\hat\uk=(k_1,\dots, k_4)$. Hence $k_4\rho_1=k_1\rho_2$. We want to show that $k_1=\pm k_4$, which will give the contradiction. To do so, we will argue as in Lemma~\ref{weswolf}. By construction, $\nu_1=e^{2\pi i\rho_1}$ and $\nu_2=e^{2\pi i\rho_2}$ are roots of $f$. We may assume that $k_1$ and $k_4$ are coprime, thus we find integers $m,n$ such that $nk_1+mk_4=1$.  Now we define $\nu:=\nu_1^n\nu_2^m=e^{2\pi i \rho}$, where $\rho:=n\rho_1+m\rho_2$. Let $g\in G(f)$ be such that $g(\nu_1)=\nu_2$ and let $l$ be its order. In the same way as in the proof of Lemma~\ref{weswolf} we see that $\nu^{{k_1}^l}=\nu^{{k_4}^l}$. This implies $\rho({k_1}^l-{k_4}^l)\in\ZZ$.
On the other hand, $\nu_1$ and $\nu_2$ are not roots of unity, see Rmk.~\ref{Reasy} (b). Hence $\rho_1$ and $\rho_2$ are irrational. Consequently, also $\rho$ is irrational and $k_1=\pm k_4$ follows.

}\end{ex}

\section{Properties of compact quotients and their fundamental groups}\label{last}
\subsection{Structural results for fundamental groups}

In Sections \ref{klumpfuss}--\ref{hans-werner} we were mainly concerned with necessary and sufficient conditions on a Cahen-Wallach space $X$ to admit a compact quotient $Y=\Gamma\backslash X$.
However, the proofs of our main results in this direction (Prop.~\ref{wulf}, Thm.~\ref{AA}, Thm.~\ref{BB}, Thm.~\ref{CC} and Thm.~\ref{DD}) also contain quite precise information
on the possible shape and structure of $\Gamma$, the fundamental group of $Y$. Indeed, we defined a group $S_\Gamma\subset G$ (Def.~\ref{vier.eins})
and have observed in Corollary~\ref{adam} that a conjugate of a finite index subgroup
of $\Gamma$ is a lattice in $S_\Gamma$. The above mentioned results then were obtained by finding conditions on and construction methods for
the possible groups $S_\Gamma$. The purpose of this subsection is to make a part of this information explicit.

Recall the decomposition $G= H\rtimes (\RR\times K)$. We start with the following observation  

\begin{pr}\label{zwitter}
Let $Y=\Gamma\backslash X$ be a compact quotient of a Cahen-Wallach space. If $\Gamma$ is not straight (see Def.~\ref{vier.sex}), then $X$ is of type $(0,2m)$ and
$$S_\Gamma= U\times D$$ 
for
a subgroup $U\subset H$ isomorphic to the Heisenberg group $H_m$ and a one parameter group $D\subset Z_G(U)$ projecting surjectively on the $\RR$-factor of $G$.
Now we assume that $\Gamma$ is straight and that $X$ is of type $(p,q)$. Let $\fa=\CC^p\oplus \CC^q$ be as in Example~\ref{Exosc2}. Let $\Gamma_0$ be
a conjugate of a finite index subgroup of $\Gamma$ such that $\Gamma_0\subset S_\Gamma$. Then
there exist a real $p$-dimensional subspace $V_{\Bbb R}\subset \CC^p$, a non-Lagrangian
real $q$-dimensional subspace $V_{I}\subset \CC^q$, and an element $\gamma_0\in \Gamma_0$ not belonging to $H\rtimes K$ 
and normalising $U:=\fz\oplus V_{\Bbb R}\oplus V_{I}$ such that
$$S_\Gamma= U\rtimes \langle \gamma_0\rangle\ .$$
There exists a power of $\gamma_0$ that acts unipotently on $U$ if and only if $p=0$.
In particular, $\Gamma$ is virtually nilpotent if and only if $p=0$.  
\end{pr}

\proof Concerning the non-straight case we have seen in the first part of the proof of Theorem~\ref{AA} that $S_\Gamma=U\cdot\psi(\Delta)$, $\Delta=\RR$, $\fu\cap\fa=\fa^{r(\psi(\Delta))}$
and that $\dim \fu\cap\fa=2m$. Since $\fa^{r(\psi(\Delta))}$ is symplectic, we conclude that $U\cong H_m$. It also follows that a generator $X$ of the Lie algebra
of $\psi(\Delta)$ acts on $\fu$ by an inner derivation $\ad X'$ for some $X'\in \fu$. We set $D:=\{\exp(t(X-X'))\mid t\in\RR\}$. Then $D\subset Z_G(U)$ and $S_\Gamma= U\times D$. 

For straight $\Gamma$ we have seen in the proof of Proposition~\ref{wulf} that $S_\Gamma=(\fz\oplus V)\rtimes \langle \gamma_0\rangle$ for a certain $(p+q)$-dimensional
subspace $V\subset\fa$. If $p=0$, it was observed that for an appropriate choice of $S_\Gamma$, i.e. of $\gamma_0$, the induced action of $\gamma_0$ on $V\cong(\fz\oplus V)/\fz$ is trivial. The decomposition $V=V_{\Bbb R}\oplus V_{I}$ was studied in the second part of the proof of Theorem~\ref{DD}.
In particular, all eigenvalues of the action of $\gamma_0$ on $V_{\Bbb R}$ have modulus different from one and there exists a special constellation $\cP$ and a $\cP$-admissible tuple $\uk$
such that $V_I\subset \CC^q$ is $(L(\cP),\phi(\cP))$-special and $\uk$-good.  
Now Proposition~\ref{lagrange} tells us that $V_I$ cannot be Lagrangian.
\qed

As a consequence of Proposition~\ref{zwitter} we state

\begin{pr}\label{eva}
Let $Y=\Gamma\backslash X$ be a compact quotient of a Cahen-Wallach space. Then $\Gamma$ is not abelian and 
$\Gamma\cap\fz$ is non-trivial.
The latter property is equivalent to: The null-leaves of the fibres (leaves) of the canonical fibration (foliation) of the quotient $Y$ are compact.
\end{pr}

\proof Let $\Gamma_0$ be
a conjugate of a finite index subgroup of $\Gamma$ such that $\Gamma_0\subset S_\Gamma$, and let $U$, $D$, $\gamma_0$ be as in Proposition~\ref{zwitter}. It suffices to prove the proposition for $\Gamma=\Gamma_0$. We set $\tilde U:=U$ if $\Gamma$ is straight and $\tilde U:=U\times D$ otherwise. We consider the lattice $\Lambda:=\tilde U\cap\Gamma_0\subset \tilde U$.
For $X$ of real type, it was shown in the second part of the proof of Theorem \ref{BB}
that $\Lambda\cap\fz\ne 0$. Moreover, the element $\gamma_0$ does not centralise $\Lambda$ by Proposition~\ref{zwitter}. If $X$ is not of real type, then we conclude from Proposition~\ref{zwitter} that $U$ is not abelian. Being a lattice in a non-abelian $1$-connected Lie group, the group $\Lambda\subset\Gamma_0$ is non-abelian, too, and the commutator group
$[\Lambda,\Lambda]\subset\fz\cap\Gamma_0$ is non-trivial.
\qed

Eventually, we use Proposition~\ref{zwitter} to obtain rather strong restrictions
for the structure of $\Gamma$ as an abstract group. By $H_n(\ZZ)$ we will denote the
following integral version of the Heisenberg group: Let $\omega$ be the standard symplectic form on $\RR^{2n}$ with standard basis $e_1,\dots,e_{2n}$, and let $H_n(\omega)=\RR\times \RR^{2n}$ be the corresponding Heisenberg group with multiplication given by (\ref{E*Hn}). Then $H_n(\ZZ)\subset H_n(\omega)$
is defined as the subgroup generated by the elements $(0,e_i)$, $i=1,\dots,2n$. It is a subgroup of index $2$ in $\frac{1}{2}\ZZ\times \ZZ^{2n}$.

\begin{pr}\label{fundi}
Let $Y=\Gamma\backslash X$ be a compact quotient of a Cahen-Wallach space of type $(p,q)$. If $p=0$, then there exists an integer $r$, $1\le r\le \frac{q+1}{2}$, such that $\Gamma$ has a finite index subgroup isomorphic to 
$$   H_r(\ZZ)\times \ZZ^{q+1-2r}\ .$$
If, in addition, $\Gamma$ is not straight, then $r=\frac{q}{2}$.

If $p\ne 0$, then either $q=0$ and a finite index subgroup of $\Gamma$ is isomorphic
to
\begin{equation}\label{lucky} \ZZ\times (\ZZ^p\rtimes_\alpha \ZZ)\ ,
\end{equation}
where $\alpha(1)\in \GL(p,\ZZ)$ is semisimple and has no eigenvalues on the unit circle,
or there exists an integer $r$, $1\le r\le \frac{p+q}{2}$, such that $\Gamma$ has a finite index subgroup isomorphic to 
\begin{equation}\label{pozzo} (H_r(\ZZ)\times \ZZ^{p+q-2r})\rtimes_\alpha\ZZ\ ,
\end{equation}
where the automorphism $\alpha(1)$ fixes the center of $H_r(\ZZ)$, and the induced element
$\tilde\alpha(1)\in \GL(p+q,\ZZ)$ is semisimple and has exactly $q$ eigenvalues on the unit circle.
\end{pr}

\proof We will use repeatedly the fact that every lattice in $H_m\times\RR^n$ has a finite
index subgroup isomorphic as an {\em abstract} group to $H_m(\ZZ)\times\ZZ^n$.
For information about lattices in Heisenberg groups the reader may consult \cite{T} or \cite{GW}.  

Let $p=0$. Then according to Proposition \ref{zwitter} there exists a finite
index subgroup of $\Gamma$ isomorphic to a lattice $\Gamma_0$ in $U\rtimes \langle \gamma_0\rangle$ ($\Gamma$ straight) or $U\times \RR$ (otherwise), and $U\cong H_s\times
\RR^{q-2s}$ for some $1\le s\le \frac{q}{2}$ ($s=\frac{q}{2}$ in the non-straight case). Moreover, we may assume that the induced action of $\gamma_0$ on $U/[U,U]\cong \RR^{q}$ is trivial. Therefore $\gamma_0$ acts on $U$ by conjugation by some $h_0\in\fa$. We define a homomorphism $\psi: \RR\rightarrow \Aut(U)$ that maps $t$ to conjugation by $th_0$. We obtain an embedding $U\rtimes \langle \gamma_0\rangle
\hookrightarrow U\rtimes_\psi\RR$, which sends $\Gamma_0$ to a lattice in $U\rtimes_\psi\RR$. We observe that $U\rtimes_\psi\RR\cong H_r\times
\RR^{q+1-2r}$, $r\in\{s,s+1\}$. Thus in both cases $\Gamma_0$ is a lattice in 
$H_r\times\RR^{q+1-2r}$, $r\in\{s,s+1\}$. Now we apply the above mentioned fact.

We turn to the case $p\ne 0$. Again, by Proposition \ref{zwitter} there exists a finite
index subgroup of $\Gamma$ isomorphic to a lattice $\Gamma_0$ in $U\rtimes \langle \gamma_0\rangle$. In particular, 
$\Gamma_0\cong \Lambda\rtimes_\alpha \ZZ$, where $\Lambda:=\Gamma_0\cap U$ is a lattice in $U$. We have $U=\fz\oplus V_{\Bbb R}\oplus V_{I}$. It is clear that $\gamma_0$ acts trivially on $\fz\subset Z(G)$. Moreover, its induced action
on $U/\fz\cong  V_{\Bbb R}\oplus V_{I}$ is semisimple, respects this decomposition, and its unimodular eigenvalues are exactly those on $V_I$.

Let us first discuss the case of abelian $U$.
By Proposition \ref{zwitter} this implies $V_I=0$, hence $q=0$. Let us denote the linear operator induced by the action of $\gamma_0$ on $U\cong \RR^{p+1}$ by $A$.
It was observed in the second part of the proof of Theorem \ref{BB} that $A$ is semisimple,
that $\ker (A-\id)=\fz$ and that $\fz\cap\Lambda\ne 0$ (see also Prop.~\ref{eva} for the latter property). We conclude that $\Lambda':=(\fz\cap \Lambda) \times (A-\id) \Lambda$ is a lattice of $U$ contained in $\Lambda$. Hence it has finite index in $\Lambda$. Moreover, it is $A$-stable.
We conclude that the finite index subgroup $\Lambda'\rtimes \langle \gamma_0\rangle\subset
\Gamma_0$ is of the form (\ref{lucky}).

It remains to discuss the case of non-abelian $U$. Thus there is an integer $1\le r\le \frac{p+q}{2}$ such that $U\cong H_r\times \RR^{p+q-2r}$. This isomorphism sends $\fz\subset U$ to the center of $H_r$. By our starting remark, the lattice $\Lambda$ has a finite index subgroup $\Lambda'$ isomorphic to $H_r(\ZZ)\times \ZZ^{p+q-2r}$. Since for any finitely generated group the set of its subgroups of fixed finite index is finite (see e.g. \cite{Hall}), we conclude that $\Lambda'$
is stable under conjugation by some power of $\gamma_0^k$ of $\gamma_0$.
It follows that $\Lambda'\rtimes \langle \gamma_0^k\rangle$ is of the form (\ref{pozzo}).     
\qed

Note that the eigenvalues of the operators $\alpha(1)$, $\tilde\alpha(1)$ appearing in the above proposition
are closely related to the parameters of the Cahen-Wallach space $X$. In fact,
let $f$ be the characteristic polynomial of that operator. Then there exist corresponding data $\cP$, $\uk$ as in Theorem~\ref{DD} such that the parameters of $X$ are given by (\ref{auf}).

The number $r=:r(Y)$ (in case (\ref{lucky}) we set $r(Y):=0$) is an interesting invariant for compact quotients that deserves further study. In particular, for quotients
of imaginary type it is the only homotopy invariant 
of $Y$ (besides dimension) invariant by going to finite covers. 
For instance, for
a fixed Cahen-Wallach space $X$ of imaginary type admitting a compact quotient one should try to determine
the minimum $r_X$ of the numbers $r(Y)$, where $Y$ runs over all compact quotients of $X$.
For example, if $X=X_{0,4}(1,1,1,1)$, then $r_X=1$.
Corollary \ref{open} implies that for all $r_X\le r\le  [(q+1)/2]$, there always exists a compact quotient $Y=\Gamma\backslash X$ such that $r(Y)=r$. 

\subsection{Quotients by groups of transvections}\label{trans}

As in Section 2, let $\hat G$ be the transvection group of a Cahen-Wallach space $X$. In this subsection we want to decide which of these spaces $X$ admit compact quotients by a group of transvections $\Gamma\subset \hat G$. 
Geometrically, these quotients are distinguished by their holonomy group, as the following proposition shows. 

\begin{pr} 
The holonomy group of a quotient $Y=\Gamma\setminus X$ of a Cahen-Wallach space is isomorphic to $\hat G_+\rtimes p_K(\Gamma)$, where $p_K:G=\hat G\rtimes K\rightarrow K$ denotes the projection to $K$. In particular, it is connected if and only if it is abelian, and this holds if and only if $\Gamma$ is contained in the group of transvections $\hat G$ of $X$. 
\end{pr}
\proof Let $\Hol_{y_0}(Y)$ denote the holonomy group of $Y$ at $y_0:=\Gamma x_0$, where $x_0=eG_+\in G/G_+$. It is well known that the identity component $\Hol^0_{y_0}(Y)$ of $\Hol_{y_0}(Y)$ is isomorphic to the holonomy group $\Hol_{x_0}(X)$ of $X$ at $x_0$ and that  $\hat G_+{\rightarrow}\,\Hol_{x_0}(X),\ g\mapsto dg_{x_0}$ is an isomorphism. Moreover, we have a surjective map
$$\pi_1(Y,y_0)\longrightarrow \Hol_{y_0}(Y)/\Hol^0_{y_0}(Y),\quad [\sigma] \longmapsto P(\sigma^{-1})\cdot \Hol^0_{y_0}(Y),$$
where $P(\sigma^{-1}):T_{x_0}X\rightarrow T_{x_0}X$ denotes the parallel translation along $\sigma^{-1}$. We compose this map with the isomorphism $\Gamma\cong \pi_1(Y,y_0)$, which sends $\gamma\in\Gamma$ to the following element $[\sigma_\gamma]$. We choose a curve $\tilde\sigma_\gamma:[0,1]\rightarrow X$ such that $\tilde\sigma_\gamma(0)=x_0$ and $\tilde\sigma_\gamma(1)=\gamma x_0$ and define $\sigma_\gamma:=\Gamma \cdot \tilde\sigma_\gamma$.

Now we will use a standard fact, which holds for arbitrary semi-Rie\-mannian symmetric spaces, see e.g. \cite{Nr} or \cite{baum}. If $\tau:[0,1]\rightarrow X$ is a curve with $\tau(0)=x_0$ and $\tau(1)=x_1$, then there exists a transvection $g\in  \hat G$ such that the parallel translation $P(\tau):T_{x_0}X\rightarrow T_{x_1}X$ along $\tau$ equals $dg_{x_0}$. Applying this to our curve $\tilde\sigma_\gamma$, we obtain a corresponding element $g\in \hat G$. In particular, $gG_+=\gamma G_+\in X=G/G_+$, thus  $g^{-1}\gamma\in G_+$.

Knowing that the parallel translation along  $\tilde\sigma_\gamma$ equals $dg_{x_0}$ we now determine the parallel translation along the inverse of the projection $\sigma_\gamma$ of $\tilde \sigma_\gamma$ to $Y$. Since $X$ is a covering of $Y$ we can identify $T_{y_0}Y\cong T_{x_0}X\cong T_{\gamma x_0}X$ in the usual way. With this identification $P(\sigma_\gamma^{-1}):T_{y_0}Y\rightarrow T_{y_0}Y$ equals
$$P(\sigma_\gamma^{-1})=(dg_{x_0})^{-1}\circ d\gamma_{x_0}=d(g^{-1}\gamma)_{x_0}.$$
Since $g^{-1}\gamma\in G_+$, we have $g^{-1}\gamma\in p_K (g^{-1}\gamma)\cdot\hat G_+=p_K\gamma\cdot\hat G_+$. Hence we obtain $\Hol_{y_0}(Y)/\Hol^0_{y_0}(Y)\cong \{p_K\gamma\mid \gamma\in\Gamma\}$. The second assertion now follows from the fact that $\Hol^0_{y_0}(Y)\cong\hat G_+$ is abelian and that $K\cap Z_G(\hat G_+)$ is trivial.\qed

We also observe
\begin{pr}\label{czecho}
A compact quotient of a Cahen-Wallach space by transvections is straight.
\end{pr}
\proof Assume that $Y=\Gamma\backslash X$ is a compact quotient that is not straight. Then we have seen in the first part of the proof of Thm.~\ref{AA}
that there exists a non-zero element $\phi\in\fk$ (in fact, $\phi$ has the same eigenvalues as $L$) such that $r(S_\Gamma)=r(\psi(\RR))=\{(t,\exp(t\phi))\mid t\in\RR\}$ . It follows that $p_K(\Gamma)=\{\exp(t\phi)\mid t\in p(\Gamma)\}$. Since $p(\Gamma)$ is dense in $\RR$, this group is non-trivial.
\qed

We now formulate a general criterion for the existence of such quotients in the language of Theorem~\ref{DD}. 

\begin{lm}\label{mickel}
Let $X$ be a Cahen-Wallach space of type $(p,q)$. Then $X$ admits a compact quotient by transvections
if and only if there exist data $(f,\cP,\uk)$ as in Theorem~\ref{DD} such that {\rm (\ref{auf})} holds, which satisfy the following
additional conditions:
\begin{enumerate}
\item[(i)] The roots of $f$ outside the unit circle are real.
\item[(ii)] The special constellation $\cP$ consists of blocks of Type I and Type II.b, only.
\end{enumerate}
As in Theorem~\ref{DD}, the assertion remains true if we require that $\nu_{c}(f)$ contains no roots of unity except $1$ and that $\cP$ is minimal.
\end{lm}

\proof By Prop.~\ref{czecho} quotients by transvections are straight. Therefore the proof of Prop.~\ref{wulf} shows that $X$ admits
a compact quotient by transvections if and only if there are
objects $V,\Lambda,t_0,\ph_0,h_0$ satisfying Conditions $(a)$ and $(b)$ of Prop.~\ref{wulf} with $\ph_0=1$.
Here we can replace the condition $\ph_0=1$ by `$\ph_0$ has finite order'.
We also observe that Condition $(ii)$ in the lemma is equivalent to $\phi(\cP)=0$.

In the proof of Thm.~\ref{DD} we constructed the objects $V,\Lambda,t_0,\ph_0,h_0$ required by Prop.~\ref{wulf} from given data $(f,\cP,\uk)$, and vice versa.
In particular, given $(f,\cP,\uk)$ we arrived at $\ph_0=\ph_{\Bbb R}\oplus\exp(2\pi\phi(\cP))$, see (\ref{didi}). If $(f,\cP)$ satisfies $(i)$
and $(ii)$, then
$\ph_0$ is of order $2$. The existence of a compact quotient by transvections follows. Vice versa, given $V,\Lambda,t_0,\ph_0,h_0$ with $\ph_0=1$
the resulting polynomial $f$ is the characteristic polynomial of a certain power of $e^{t_0L}$. Hence $f$ satisfies $(i)$. In addition, we
constructed a torus $\tilde T\subset O(\fa_I)$ together with a decomposition of its Lie algebra $\tilde\ft=\ft_+\oplus\ft_-$ such that $\phi(\cP)\in\ft_+$. The condition $\ph_0=1$ implies $\ft_+=0$, hence $\phi(\cP)=0$. Therefore $\cP$ satisfies Condition $(ii)$.
\qed

\begin{co}
The set of isometry classes of Cahen-Wallach spaces admitting a compact quotients by transvections is countable. 
\end{co}

\proof It is easily checked that the set of triples $(f,\cP,\uk)$ satisfying the conditions of Thm.~\ref{DD} and Lemma~\ref{mickel} is countable.
\qed

The following more explicit criteria for the existence of compact quotients by transvections of spaces of real and imaginary type
are direct consequences of Lemma \ref{mickel} (in the imaginary case one should take the second part of Remark~\ref{Rri} into account).
Alternatively, they can be derived by checking the proofs of Theorem~\ref{BB} and Theorem~\ref{CC} for consequences of the additional condition $\ph_0=1$.

\begin{pr}\label{bb}
Let $X$ be an $(n+2)$-dimensional Cahen-Wallach space of real type. Then $X$ admits a compact quotient by transvections if and only if there exists a polynomial
$f\in\ZZ[x]$ of degree $n$ of the form {\rm(\ref{otto})}
with all roots real and different from $\pm 1$ such that
$$X\cong X_{n,0}(\log |\nu_1|,\log |\nu_2|,\dots,\log |\nu_n|)\  ,$$
where $\nu_1,\nu_2,\dots,\nu_n$ are the roots of $f$.\qedohne
\end{pr}

\begin{pr}\label{cc}
Let $X$ be an $(n+2)$-dimensional Cahen-Wallach space of imaginary type. Then $X$ admits a compact quotient by transvections if and only if there exists an $\RR$-admissible $n$-tuple $\uk=(k_1,\dots,k_n)\in(\ZZ_{\not=0})^n$ such
that 
$X\cong X_{0,n}(k_1,\dots, k_n)$.\qedohne
\end{pr}

Also for spaces of mixed type the criterion provided by Lemma \ref{mickel} can be made much more explicit.
For this we need the following counterpart of the notion of $\RR$-admissibility.

\begin{de}
Let $n$ be even. An $n$-tuple $\uk=(k_1,\dots,k_n)\in\ZZ^n$ is called $\CC$-admissible if there exists a $\uk$-good complex vector subspace $V\subset \CC^n$.
\end{de}

For an $n$-tuple, the condition of $\CC$-admissibility is invariant under permutations of the coordinates
and under multiplication with a common factor $m\in \ZZ_{\not=0}$. It is also invariant under translation by tuples of the form $(k,k,\dots,k)$.
In contrast to $\RR$-admissibility, it is not invariant under independent sign changes of the coordinates.

For $\rho\in\RR^*$ and $d\in \NN$, let $\cP_d(\rho)$ be the special constellation consisting of $d$ copies of the block 
$\displaystyle {\rho\choose \rho}$ of Type II.b. 

\begin{lm}\label{hm}
\begin{enumerate}
\item[(a)] A $2d$-tuple $\uk=(k_1,\dots,k_{2d})\in\ZZ^{2d}$ is $\CC$-admissible if and only if it is $\cP_d(\rho)$-admissible for one (equivalently: for all) $\rho\in \RR\setminus\{0,-k_1,\dots,-k_{2d}\}$.
\item[(b)] Let $(\rho_1,\dots,\rho_r)\in (\RR\setminus\ZZ)^r$ be such that $|\rho_i|\ne |\rho_j|$ for $i\ne j$. We consider the special constellation $\cP=(0|\dots|0|\cP_{d_1}(\rho_1)|\dots|\cP_{d_r}(\rho_r))$ containing precisely 
$d_0$ blocks of Type I. Then a $d(\cP)$-tuple $\uk$ is $\cP$-admissible if and only if $\uk=(\uk^0,\uk^1,\dots,\uk^r)$, where $\uk^0$ is an $\RR$-admissible $d_0$-tuple,  and $\uk^i$, $i=1,\dots,r$, is a $\CC$-admissible $2d_i$-tuple.
\end{enumerate}
\end{lm}
\proof Assertion $(a)$ is an immediate consequence of the definitions. It is also clear that tuples of the form $\uk=(\uk^0,\uk^1,\dots,\uk^r)$ as in $(b)$
are $\cP$-admissible. Let now $\uk$ be $\cP$-admissible. We can write $\uk=(\uk^0,\uk^1,\dots,\uk^r)$ 
for some $\uk^0\in \ZZ^{d_0}$, $\uk^i\in \ZZ^{2d_i}$ for $i>0$. 
Set $n=d_0+2d_1+\dots+2d_r$, and let $V\in G_n(\CC^n)^{L(\cP)+\phi(\cP),i\rho(\cP)}$ be a $\uk$-good subspace. The splitting $\CC^n=\CC^{d_0}\oplus\CC^{2d_1}\oplus\dots\oplus\CC^{2d_r}$ coincides with the decomposition into eigenspaces of $(L(\cP)+\phi(\cP))^2$. The corresponding eigenvalues are $-\rho_0^2,-\rho_1^2,\dots,-\rho_n^2$ with $\rho_0=0$. Let 
$V= \bigoplus_{i=0}^r V(\rho_i)$ be the eigenspace decomposition of $V$. Then $V(0)\subset \CC^{d_0}$ is $\uk^0$-good, while for $i>0$ the spaces
$V(\rho_i)\subset \CC^{2d_i}$ are $(L(\cP_{d_i}(\rho_i)),\phi(\cP_{d_i}(\rho_i)))$-special and $\uk^i$-good. It follows that $\uk^0$ is $\RR$-admissible and,
taking Assertion $(a)$ into account, that $\uk^i$ is $\CC$-admissible for $i=1,\dots,r$.
\qed

We have the following analog of Proposition~\ref{charles}.

\begin{pr}
Let $\uk=(k_1,\dots,k_{2d})\in\ZZ^{2d}$ be $\CC$-admissible. Then there exists a decomposition $\{1,\dots,2d\}=I\cup J$ into two disjoint subsets of cardinality $d$ such that
\begin{equation}\label{dirichlet} 
\sum_{i\in I} k_i =\sum_{j\in J} k_j\ .
\end{equation}
\end{pr}

\proof By Lemma~\ref{hm}, $(a)$, the tuple $\uk=(k_1,\dots,k_{2d})\in\ZZ^{2d}$ is $\cP_d(\rho)$-admissible for some $\rho\ne 0$. By Proposition~\ref{fourier} there exists a sign vector $\kappa\in\{1,-1\}^{2d}$ such that
$\langle \kappa, \mu(\cP_d(\rho))\rangle=0$ and $\langle \kappa, \uk\rangle=0$. Set $I:=\{i\mid \kappa_i=1\}$, and let $J$ be its complement. Note that
$\mu(\cP_d(\rho))=(\rho,\rho,\dots,\rho)$. It follows that $I$ has cardinality $d$ and that (\ref{dirichlet}) holds.
\qed

It is clear that  $(k,k)\in\ZZ^2$ is $\CC$-admissible for all $k\in\ZZ$. In addition, Example~\ref{Exdim4} shows that all tuples $\uk\in\ZZ^4$
satisfying (\ref{dirichlet}) for some decomposition of equal cardinality are $\CC$-admissible. Therefore, Condition $(\ref{dirichlet})$ is equivalent to
$\CC$-admissiblity for dimensions $2d\le 4$. In analogy to the question at the end of Subsection~\ref{wagner} one may ask whether this is also true in higher dimensions.

We return to Cahen-Wallach spaces admitting compact quotients by transvections. We first consider some examples.

\begin{ex}\label{mufti}{\rm
Let $h\in \ZZ[x]$ be an irreducible polynomial of the form (\ref{otto}) such that all its roots outside the unit circle are real.
We assume that both types of roots (real and unimodular ones) occur. Then we can choose positive real numbers $\lambda_1,\dots,\lambda_r$, $\mu_1,\dots,\mu_s$ such that the absolute values of the real roots of $h$ are precisely $e^{\pm \lambda_1},\dots,e^{\pm \lambda_r}$, whereas the roots of $h$ on the unit circle are given by $e^{\pm i\mu_1}, \dots, e^{\pm i\mu_s}$ (see Remark~\ref{Reasy}, $(a)$). Then 
$$X:=X_{2r,2s}(\lambda_1,\lambda_1,\dots,\lambda_r,\lambda_r;\mu_1,\mu_1,\dots,\mu_s,\mu_s)$$
admits a compact quotient by transvections. Indeed, if $\cP:=(\cP_{1}(\frac{\mu_1}{2\pi})|\dots|\cP_{1}(\frac{\mu_s}{2\pi}))$ and $\uk:=0$, then $(h,\cP,\uk)$ satisfies the
conditions required by Thm.~\ref{DD} and Lemma~\ref{mickel} and $X$ is isometric to the space associated with $(h,\cP,\uk)$ by (\ref{auf}).
}\end{ex}

\begin{ex}\label{schuft}{\rm
Let $h,r,s,\lambda_i,\mu_j$ be as in the previous example. Fix $d\in\NN$, and let $\uk^1,\dots,\uk^s$ be a collection of $s$ $\CC$-admissible $2d$-tuples.
Let $X$ be the Cahen-Wallach space 
$$X_{2dr,2ds}(\underbrace{\lambda_1,..,\lambda_1}_{2d},\dots,\underbrace{\lambda_r,..,\lambda_r}_{2d};
\mu_1+2\pi k_1^1,
..,\mu_1+2\pi k^1_{2d},\dots,\mu_s+2\pi k_1^s,..,\mu_s+2\pi k^s_{2d})\,.$$
It admits a compact quotient by transvections. Again, this follows from Lemma~\ref{mickel} by observing that $X$ is isometric to a Cahen-Wallach
space associated with some data $(f,\cP,\uk)$. Here we take $f=h^d$, $\cP=(\cP_{d}(\frac{\mu_1}{2\pi})|\dots|\cP_{d}(\frac{\mu_s}{2\pi}))$, $\uk=(\uk^1,\dots,\uk^s)$.
Note that for $d=1$ we get exactly the same spaces as in Example~\ref{mufti}.
}\end{ex}


In Example~\ref{Excomp} we defined the notion of composition of two Cahen-Wallach spaces $X_1$ and $X_2$ admitting compact quotients resulting in a new Cahen-Wallach space with compact quotients. This composition depends on the chosen  data  $(f_1,\cP_1,\uk_1)$ and $(f_2,\cP_2,\uk_2)$ determining $X_1$ and $X_2$, respectively. Furthermore, we specialised compositions to the case, were $X_1$ and $X_2$ are of purely real or imaginary type only using the parametrisation of $X_1$ and $X_2$ according to Theorems~\ref{BB} and~\ref{CC}. Here we want to consider compositions based on the special data $(f,\cP,\uk)$ used in Example~\ref{schuft} and to spaces of real or imaginary type given as in Propositions~\ref{bb} and ~\ref{cc}. 
Let us describe the resulting compositions only using the parameters of $X_1$ and $X_2$. If $X_i$, $i=1,2$, is of real or of mixed type, then let $(\lambda^i;\mu^i)$ be the parameters of $X_i$ given by  Proposition~\ref{bb} (with $\mu^i=\emptyset$) and  Example~\ref{schuft}, respectively. If $X_i$ is of imaginary type given as in Proposition \ref{cc}, then put $(\lambda^i;\mu^i):= (\emptyset,2\pi\uk)$. Then the composition of $X_1$ and $X_2$ (with respect to these data) is given by the parameters $(\lambda^1,\lambda^2;\mu^1,\mu^2)$.

Compositions can obviously be defined also for a finite number of spaces $X_i$, $i=1,\dots,m$.

\begin{pr}\label{dd}
A Cahen-Wallach space admits a compact quotient by transvections if and only if  it is composed of one or  more of the following spaces:
\begin{itemize}
\item spaces of real type as in Proposition~\ref{bb},
\item spaces of imaginary type as in Proposition \ref{cc},
\item spaces of mixed type as in Example~\ref{schuft}.
\end{itemize} 
\end{pr}

\proof We already know that the spaces listed in the proposition admit compact quotients by transvections. Conditions $(i)$, $(ii)$ of Lemma~\ref{mickel} are compatible with composition. Thus also compositions of spaces in the above list have compact quotients by transvections. It remains to show that every Cahen-Wallach space $X$ defined by data $(f,\cP,\uk)$ satisfying the conditions of Thm.~\ref{DD} and Lemma~\ref{mickel} is composed of spaces in the above list.
We may assume that $\cP$ is minimal and that the only root of unity among the zeroes of $f$ is $1$ (if there is one).
We decompose $f=f_0 (x-1)^{d_0} f_1^{d_1}\dots f_l^{d_l}$, such that $f_i\in\ZZ[x]$, the roots of $f_0$ are in $\RR\setminus\{-1,0,1\}$, $f_1,\dots,f_l$ are irreducible, pairwise different, and have at least one pair of complex conjugate roots on the unit circle. 
Now Lemma~\ref{mickel}, $(ii)$, together with the minimality condition (\ref{stuber}) implies that $\cP$ is of the form
$$\cP=(0|\dots|0|\cP_{d_1}(\rho_{11})|\dots|\cP_{d_1}(\rho_{1s_1})|\dots|\cP_{d_l}(\rho_{l1})|\dots|\cP_{d_l}(\rho_{ls_l}))$$
with precisely 
$d_0$ blocks of Type I and such that the roots of $f_i$ on the unit circle are precisely $e^{\pm 2\pi i\rho_{i1}}, \dots, e^{\pm 2\pi i\rho_{i s_i}}$.
By Lemma~\ref{hm}, $(b)$, there exist an $\RR$-admissible $d_0$-tuple $\uk^0$ and $\CC$-admissible $2d_i$-tuples $\uk^{ij}$ such that 
$\uk=(\uk_0,\uk^{11},\dots,\uk^{ls_s})$. If follows that $(f,\cP,\uk)$ is composed of the following data
\begin{itemize}
\item $(f_0,\emptyset,\emptyset)$,
\item $((x-1)^{d_0},(0|\dots|0),\uk_0)$,
\item $\left( f_i^{d_i},(\cP_{d_i}(\rho_{i1})|\dots|\cP_{d_i}(\rho_{is_i})),(\uk^{i1},\dots,\uk^{is_i})\right)$, $i=1,\dots, k$.
\end{itemize}
This decomposition defines the desired decomposition of $X$.
\qed

\subsection{Solvmanifolds}\label{solv}
A Lorentzian solvmanifold (shortly solvmanifold in the following) is a quotient $\Gamma\backslash S$, where $S$ is
a $1$-connected solvable Lie group equipped with a left-invariant Lorentzian metric, and $\Gamma\subset S$ is a discrete subgroup. 
In this subsection we decide, using Theorem~\ref{DD}, which Cahen-Wallach spaces have compact quotients that are
solvmanifolds. 
Note that such a compact quotient of a Cahen-Wallach space $X$ is necessarily of the following form: $S$ is a connected solvable subgroup of the
group of $G\subset \Iso(X)$ acting simply transitively on $X$ and $\Gamma\subset S$ is a lattice. 
We also get information on the possible shapes of $S$ and $\Gamma$.

The proof of Proposition~\ref{dada} below combined with Proposition~\ref{gaga} shows even more: For a given compact quotient $\Gamma\backslash X$ we can decide, whether it is finitely covered by a solvmanifold. Equivalently, we can decide, whether a finite index subgroup $\Gamma_0\subset \Gamma$ has a syndetic hull in $G$, i.e. a connected
subgroup $\tilde S\subset G$ containing $\Gamma_0$ as a lattice (for syndetic hulls compare \cite{WM2, FG}).
In view of the following lemma, which could have been included in Section~\ref{klumpfuss}, it is not surprising that syndetic hulls and solvmanifolds are related. 


\begin{lm}\label{haha}
Let $S\subset G$ be a connected subgroup.
Then $S$ acts properly and cocompactly on $X$ if and only if there exists a connected solvable cocompact subgroup $S_1\subset S$ acting simply transitively
on $X$.
\end{lm}

We omit the proof. In fact, rather than Lemma~\ref{haha} we need an analogous statement for pairs $(S,\Gamma)$, where $S\subset G$ is connected and
$\Gamma\subset S$ is a lattice. Proposition \ref{gaga} below provides such a statement. 

A nilmanifold is a solvmanifold $\Gamma\backslash S$ with $S$ nilpotent.

\begin{pr}\label{hihi}
Every non-straight compact quotient $Y=\Gamma\backslash X$ of a Cahen-Wallach space is finitely covered by a nilmanifold.
\end{pr}
\proof From Proposition~\ref{zwitter} we see that $S_\Gamma$ is connected nilpotent and, using Lemma~\ref{kirsten}, that it acts simply transitively on $X$. There is a finite index subgroup $\Gamma_0$ of a conjugate of $\Gamma$ that is
a lattice in $S_\Gamma$. Thus $Y$ is covered by $\Gamma_0\backslash S_\Gamma$ (the conjugation is incorporated in the covering map).
\qed

We now give a construction of certain compact solvmanifolds covered by a given 
Cahen-Wallach space $X$ provided the following objects are given:
\begin{itemize}
\item Elements $\phi\in\fk$, $X\in\fa^{L+\phi}$. They define a one parameter group $\psi: \RR\rightarrow G$ by 
$\psi(t):=(tX,t,e^{t\phi})\in H\rtimes (\RR\times K)$.
\item An $(L+\phi)$-invariant subspace $V\subset\fa$ such that $\fa=V\oplus \fa_+$. It defines a subgroup $U:=\fz\oplus V\subset H$.
\item Elements $t_0\in\RR\setminus\{0\}$, $u_0\in U$ and a lattice $\Lambda\subset U$ stable under conjugation by $u_0\psi(t_0)$.
\end{itemize}
We set $S:=U\cdot\psi(\RR)$, $\Gamma:=\Lambda\cdot \langle u_0\psi(t_0)\rangle$.
Then $S$ is connected solvable, $\Gamma\subset S$ is a lattice and,
by Lemma~\ref{kirsten}, the action of $S$ on $X$ is simply transitive. Thus $\Gamma\backslash S$ is a compact solvmanifold covered by $X$. A {\bf standard solvmanifold} is a solvmanifold arising in this way.    

\begin{pr}\label{gaga}
Let $Y=\Gamma\backslash X$ be a compact quotient of a Cahen-Wallach space. Then the following
assertions are equivalent:
\begin{enumerate}
\item[(i)] A finite index subgroup of $\Gamma$ has a syndetic hull in $G$.
\item[(ii)] $Y$ is finitely covered by a solvmanifold.   
\end{enumerate}
If $Y$ is straight, then these assertions are equivalent to
\begin{enumerate}
\item[(iii)] $Y$ is finitely covered by a standard solvmanifold.    
\end{enumerate}
\end{pr}

\proof For non-straight $Y$ the proposition is an immediate consequence of Prop.~\ref{hihi}.
Let $Y$ be straight. The implications $(iii)\Rightarrow (ii)\Rightarrow (i)$ are obvious. We have to prove $(i)\Rightarrow (iii)$.
For that we may assume that $\Gamma\subset S_\Gamma$ and that $\Gamma$ has a syndetic hull $\tilde S$. We have $S_\Gamma = U\rtimes\langle \gamma_0\rangle$ as in Proposition~\ref{zwitter}.
In particular, $U=\fz\oplus V$ for some subspace $V\subset \fa$ invariant under conjugation by $\gamma_0$ satisfying $\fa=V\oplus \fa_+$,  $\Lambda:=
U\cap \Gamma$ is a lattice in $U$, and $\Gamma=\Lambda\cdot\langle \gamma_0\rangle$. 

The subgroup $r(\tilde S)\subset\RR\times K$ is connected, and $K$ is compact. It follows that there is a one parameter subgroup $C\subset r(\tilde S)$ containing $r(\gamma_0)$. Let $\phi\in\fk$ and $t_0\ne 0$ be such that $C=\{(t,e^{t\phi})\in\RR\times K\mid t\in\RR\}$ and $r(\gamma_0)=(t_0,e^{t_0\phi})$. The group $\tilde S$ acts properly and cocompactly on $X$. 
Therefore $\tilde S\cap H$ acts properly on the typical fibre $H/\fa_+$ of the canonical fibration. Let $\tilde U\subset H$ be the unique connected
subgroup of $H$ such that $(\tilde S\cap H)\backslash \tilde U$ is compact (see Lemma~\ref{drei.neun}). Then also $\tilde U$ acts properly on $H/\fa_+$.
Since $\Lambda\subset \tilde S\cap H$, we have $U\subset\tilde U$. Now Lemma~\ref{ulf} implies that $U=\tilde U$.
Since $\tilde S\cap H$ is a normal subgroup of $\tilde S$, the group $U=\tilde U$ is normalised by $\tilde S$.
We define $S':=U(\tilde S\cap r^{-1}(C))$. The  group $S'$ contains $U$ and $\gamma_0$, hence $\Gamma$. 

We finish the proof by showing that the pair $(S',\Gamma)$ is conjugate in $G$ to a pair $(S,\Gamma')$ such that $\Gamma'\backslash S$ is 
a standard solvmanifold. We choose $X'\in\fa$ such that $C':=\{\exp(t(X',1,\phi))\mid t\in\RR\}\subset S'$. 
Then $S'=U\cdot C'$ and $\gamma_0=u_0\exp(t_0(X',1,\phi))$ for some $u_0\in U$. 
We find elements $X\in \fa^{L+\phi}$, $Y\in\fa$ such that $X'= X+ (L+\phi)(Y)$. Now we view $h:=Y$ as an element of $H\subset G$
and conjugate by it. In particular, $hUh^{-1}=U$, $\Ad(h)(X',1,\phi)=e^{\ad(Y)} (X',1,\phi)=(Z+X,1,\phi)$ for some $Z\in\fz$.
We set $S:=hS'h^{-1}$, $\Gamma':=h\Gamma h^{-1}$, $\Lambda':=h\Lambda h^{-1}$, $\gamma_0':=h\gamma_0 h^{-1}$, $\psi(t):=\exp(t(X,1,\phi))=(tX,t,e^{t\phi})$.
It follows that $S=U\cdot \psi(\RR)$, $\gamma_0'=hu_0h^{-1}(t_0Z)\psi(t_0)$. We conclude that $\Gamma'\backslash S$ is 
the standard solvmanifold associated with the objects $\phi,X,V,t_0, u_0':=hu_0h^{-1}(t_0Z),\Lambda'$.
\qed

\begin{pr}\label{dudu}
Every compact quotient $Y=\Gamma\backslash X$ of a Cahen-Wallach space $X$ of real type is finitely covered by a solvmanifold.
\end{pr}

\proof We may assume that $\Gamma\subset S_\Gamma$, where $S_\Gamma=(\fz\oplus V)\cdot \langle \gamma_0\rangle$ as in Prop.~\ref{zwitter},
$r(\gamma_0)=:(t_0,\ph_0)\in \RR\times K$, $t_0\ne 0$ (compare Prop.~\ref{wulf}). For spaces of real type, we have observed in the second part of the proof of Thm.~\ref{BB} that $L(V)=V$. It follows that $V$ is also invariant under the closure in $K$ of the group generated by $\ph_0$. Let $T\subset K$
be the identity component of that closure. Then $\tilde S:=(\fz\oplus V)\rtimes (\RR\times T)\subset G$ is a connected subgroup
that contains $\gamma_0^k$ for some $k\in\NN$. It follows that $\tilde S$ is a syndetic hull for the finite index subgroup 
$(\Gamma\cap(\fz\oplus V))\cdot \langle \gamma_0^k\rangle\subset \Gamma$. Now we apply Proposition~\ref{gaga}.
\qed

\begin{pr}\label{dada}
Let $X$ be a Cahen-Wallach space of type $(p,q)$. Then $X$ covers a compact solvmanifold if and only if there exist
\begin{enumerate}
\item[(a)] a polynomial $f\in\ZZ[x]$ of degree $p+q$ of the form {\rm (\ref{otto})} having precisely $q$ roots on the unit circle (counted with multiplicity),
\item[(b)] a special constellation $\cP$ of dimension $q$ containing no blocks of Type I and Type II.c and satisfying $\nu(\cP)=\nu_c(f)$
\end{enumerate}
such that
\begin{equation}\label{auf1} X\cong X_{p,q}\Big(\log |\nu_1|,\dots,\log |\nu_p|;\,2\pi\mu(\cP)\,\Big)\ ,
\end{equation}
where $\nu_1,\dots,\nu_p$ are the roots of $f$ of modulus different from $1$.
\end{pr}

\proof Assume that a pair $(f,\cP)$ satisfying $(a)$ and $(b)$ is given. We consider the Cahen-Wallach space $X=X_{p,q}\Big(\log |\nu_1|,\dots,\log |\nu_p|;\,2\pi\mu(\cP)\,\Big)$ as in $(\ref{auf1})$. As in the proof of Theorem~\ref{DD} we work with the splitting
$\fa=\fa_{\Bbb{R}}\oplus\fa_I$, $L=L_{\Bbb R}\oplus L_I$. We have $L_I=2\pi L(\cP)$. The first part of that proof provides, starting from $(f,\cP)$ and $\uk=0$ (it is $\cP$-admissible),
a subspace $V= V_{\Bbb{R}}\oplus V_I\subset \fa$ transversal to $\fa_+$, an element $\ph_{\Bbb{R}}\in \grO(\fa_{\Bbb{R}})\cap K$ and a lattice
$\Lambda\subset \fz\oplus V$, where $V_{\Bbb{R}}$ is invariant under $L_{\Bbb R}$ and $\ph_{\Bbb{R}}$ and $V_I$ is invariant under $L(\cP)+\phi(\cP)$.
Moreover, if we set $\gamma_0:=(0,1,\ph_{\Bbb{R}}\oplus\exp(2\pi\phi(\cP)))\in H\rtimes(\RR\times K)=G$, then $\Lambda$ is stable under conjugation
by $\gamma_0$. The group $\Gamma:=\Lambda \cdot \langle \gamma_0\rangle$ gives rise to a compact quotient $Y=\Gamma\backslash X$.
Let $T\subset K$ be the identity component of the closure of the group generated by $\ph_{\Bbb{R}}$. We set $D:=\{(t,\exp(2\pi t \phi(\cP)))\mid t\in\RR\}\subset \RR\times K$. Then $D$ commutes with $T$, and as in the proof of Prop.~\ref{dudu} we see that $\tilde S:=(\fz\oplus V)\rtimes (D\times T)$
is a syndetic hull for a finite index subgroup of $\Gamma$. Proposition~\ref{gaga} implies that $X$ covers a compact solvmanifold.

Now we assume that a Cahen-Wallach space $X$ has a compact quotient $Y$ that is a solvmanifold. If $Y$ is not straight, then $X$ is a group manifold (see Theorem~\ref{AA}) which is of the form (\ref{auf1}) ($f$ is a power of $x-1$, $\cP$ consists of blocks of Type II.a, only). Thus we can assume that $Y$ is straight.
By Prop.~\ref{gaga} the solvmanifold $Y$ is finitely covered by a standard solvmanifold $\Gamma\backslash S$. Among its defining objects
we are particulary interested in $\phi,V,t_0,\Lambda=\Gamma\cap(\fz\oplus V)$. Proposition~\ref{eva} tells us that $\Lambda\cap\fz$ is non-trivial.
Thus the projection $\Lambda_0$ of $\Lambda$ to $V$ is a lattice in $V$. It is stabilised by the operator $\exp(t_0(L+\phi))$. Let $f$ be the
characteristic polynomial of the restriction of that operator to $V$. By Lemma~\ref{gigi} the polynomial $f$ is integral of the form (\ref{otto}).
Again we look at the splitting $\fa=\fa_{\Bbb{R}}\oplus\fa_I$, $L=L_{\Bbb R}\oplus L_I$. It induces splittings $\phi=\phi_{\Bbb R}\oplus\phi_I$, $V=V_{\Bbb R}\oplus V_I$. The subspace $V_I\subset\fa_I$ is $(L_I,\phi_I)$-special. By Proposition~\ref{veryspecial} there exists a special
constellation $\cP$ with $\nu(\cP)=\nu_c(f)$  such that $(\frac{t_0}{2\pi}L_I,\frac{t_0}{2\pi}\phi_I)\cong (L(\cP),\phi(\cP))$. Since $L_I$ is invertible, the special constellation $\cP$ does not contain
blocks of Type I and Type II.c. As in the proof of Theorem~\ref{DD} we eventually conclude that (\ref{auf1}) holds.
\qed

Note that in contrast to Theorem~\ref{DD} and Lemma~\ref{mickel} it is essential to allow non-minimal special constellations $\cP$ in Proposition~\ref{dada}. 

\begin{co}
Let $X\cong X_{0,n}(\alpha_1,\dots,\alpha_p,\beta_1,\beta_1,\dots,\beta_q,\beta_q)$, $n=p+2q$, $|\alpha_i|\ne|\alpha_j|$ for $i\ne j$, be a Cahen-Wallach
space of imaginary type. Then $X$ covers a compact solvmanifold if and only if 
\begin{itemize}
\item
$p=0$ or 
\item $p>0$, the quotients $\alpha_i/\alpha_1$, $i=1,\dots,p$, are rational, and Conditions $(i)$, $(ii)$, $(iii)$ in Proposition~\ref{special} are satisfied with the additional requirement $\beta_j/\alpha_1\in\QQ$ for
$j\in I_k$, $k=1,\dots,p_0,2p_0+1,\dots,p$.
\end{itemize}
\end{co}

\proof One can either specialise Proposition~\ref{dada} to 
spaces of imaginary type 
or combine Proposition~\ref{gaga} with Proposition~\ref{special} using the fact that the subspace $V\cap\fa_1\subset \fa_1$ constructed in the classification part of the proof of Theorem~\ref{CC} is $(L,\phi)$-special.
\qed


\subsection{Moduli spaces in small dimensions}\label{low}

Let $\cM_{p,q}$ be the space of isometry classes of Cahen-Wallach spaces of type $(p,q)$ as in Subsection~\ref{class}, and let $\cM^c_{p,q}$ be its subspace consisting
of classes of spaces having compact quotients. The main results of this paper, in particular Theorem~\ref{DD}, describe the space $\cM^c_{p,q}$.
Here we want to use the parametrisation of $\cM_{p,q}$ given by (\ref{humus}) to make this description completely explicit for $p,q\le 3$ and for $(p,q)=(0,4)$. Moreover, we also describe the subspaces $\cM^t_{p,q}\subset \cM^c_{p,q}$ and $\cM^s_{p,q}\subset\cM^c_{p,q}$ of classes of Cahen-Wallach spaces admitting quotients by transvections and of those covering compact solvmanifolds, respectively. The results are given in Table~\ref{obermufti}.

\begin{sidewaystable}[p]

\begin{center}
\renewcommand{\arraystretch}{1.2}
\begin{tabular}{|c|c|c|c|c|l|}
\hline
&&&&&\\
Type&$\cM^0_{p,q}$&$\cM^t_{p,q}$&$\cM^c_{p,q}\setminus\cM^t_{p,q}$&$\cM^s_{p,q}$&Further information\\
&&&&&\\
\hline
\hline
$(0,2)$&$(1,1)$&&&&\\
&&$\checkmark$&$\emptyset$&\checkmark&\\
\hline
$(0,3)$&$(1,\mu,\mu+1)$&&&&\\
&$\mu\in [1,\infty)$&$\mu\in\QQ$&$\emptyset$&$\mu=1$&Example~\ref{Exdim3}\\
\hline
$(0,4)$&$(1,\mu_1,\mu_2,\mu_1+\mu_2+\kappa)$&&&&\\
&$1\le\mu_1\le\mu_2$, $\kappa\in \{1,-1\}$&$\mu_1,\mu_2\in\QQ$&$\mu_1=1$, $\mu_2\notin\QQ,$ &$\mu_1=1$ or $\mu_1=\mu_2$&Example~\ref{Exdim4}\\
&&&$\kappa=-1$&&\\
\hline
\hline
$(2,0)$&$(1,1)$&&&&\\
&&$\checkmark$&$\emptyset$&\checkmark&\\
\hline
$(2,2)$&$(\lambda,\lambda;1,1)$&&&&\\
&$\lambda\in (0,\infty)$&$\lambda\in\cF_2\cup\cF_4$&$\lambda\notin\cF_2\cup\cF_4$&\checkmark&Example~\ref{Ex22}\\
\hline
$(2,3)$&$(\lambda,\lambda;1,\mu,\mu+1)$&&&&\\
&$\lambda\in (0,\infty)$, $\mu\in[1,\infty)$&$\lambda\in\cF_2$, $\mu\in\QQ$&$\lambda\in2\cF_4\setminus \cF_2$, $\mu=1$&$\mu=1$&Example~\ref{Ex23}\\
\hline
\hline
$(3,0)$&$(1,\lambda,\lambda+1)$&&&&\\
&$\lambda\in [1,\infty)$&$\lambda\in\cF_3$&$\lambda=1$&$\checkmark$&Proposition~\ref{duli}\\
\hline
$(3,2)$&$(\lambda_1,\lambda_2,\lambda_1+\lambda_2;1,1)$&&&&\\
&$0<\lambda_1\le\lambda_2$&$(\lambda_1,\lambda_2)\in\widetilde\cF_3$&$\lambda_2/\lambda_1\in\{1\}\cup\cF_3$,&$\checkmark$&Example~\ref{Ex32}\\
&&&$(\lambda_1,\lambda_2)\notin\widetilde\cF_3$&&\\
\hline
$(3,3)$&$(\lambda_1,\lambda_2,\lambda_1+\lambda_2;1,\mu,\mu+1)$&&&&\\
&$0<\lambda_1\le\lambda_2$, $\mu\in [1,\infty)$&$(\lambda_1,\lambda_2)\in\widetilde\cF_3$,&$\lambda_1=\lambda_2\in\cF_3'$,&$\mu=1$&\\
&&$\mu\in\QQ$&$\mu\in\QQ$&&\\
\hline
\end{tabular}
\caption{The parameters of low-dimensional Cahen-Wallach spaces with compact quotients}\label{obermufti}
\end{center}
\end{sidewaystable} 

Our starting point (the second column of the table) is to determine the parameters for the space $\cM^0_{p,q}\supset\cM^c_{p,q}$ of Cahen-Wallach spaces
satisfying the trace condition of Corollary~\ref{spurnull3}. It is a manifold with corners of dimension $p+q-3$ (mixed type), $p-2$ or $q-2$ (real
and imaginary type). In the next column we give the {\em additional} conditions  on the parameters that determine $\cM^t_{p,q}\subset \cM^0_{p,q}$.
These conditions are expressed in terms of certain sets $\cF_2$, $\cF_3$ \frenchspacing etc. that \nonfrenchspacing are defined below.
Note that $\cM^t_{p,q}$ is countable and dense in $\cM^0_{p,q}$. The next column describes the subset $(\cM^c_{p,q}\setminus\cM^t_{p,q})\subset \cM^0_{p,q}$ in the same way. Note that $\cM^c_{p,q}$ sometimes contains one-dimensional families of Cahen-Wallach spaces.
In the fifth column we list the conditions that determine $\cM^s_{p,q}$ as a subset of $\cM^c_{p,q}$.
Thus, in order to read off $\cM^s_{p,q}$ from that column, one first has to determine $\cM^c_{p,q}$ using the previous two columns. 
We see that $(p,q)=(3,3)$  is the only type appearing in the table such that  $\cM^c_{p,q}$ is strictly larger than $\cM^t_{p,q}\cup \cM^s_{p,q}$. 

The sign `$\emptyset$' says that the corresponding space is empty, whereas the empty condition is indicated by `$\checkmark$'.
Recall that for spaces of types $(1,q)$ and $(p,1)$ we have $\cM^c_{p,q}=\cM^0_{p,q}=\emptyset$. Therefore these types do not appear in the table.

We now define the relevant sets. They are related to certain integral polynomials of the form (\ref{otto}) (equivalently to units in certain number fields)
of degree 2, 3 and 4. Sometimes the description will be given in terms of the $\QQ$-vector space $\cH_K$ associated with a number field $K$ (see Definition~\ref{blasebalg}).

\begin{enumerate}
\item[$\cF_2:$] $=\left\{\frac{r}{\pi}\log(\frac{1}{2}(k+\sqrt{k^2-4}))\mid r\in\QQ^+,\, 3\le k\in\NN\right\}$. This set can be structured as follows (also avoiding repetitions in the list of elements). Let $d\ge 2$ be a square free integer, and let $\nu_d>1$ be a unit in $\QQ(\sqrt{d})$ (e.g. the fundamental one). If we define
$\Lambda_d:=\left\{\frac{r}{\pi}\log\nu_d\mid r\in\QQ^+\right\}=\left\{\lambda\in(0,\infty)\mid (\pi\lambda,-\pi\lambda)\in\cH_{{\Bbb Q}(\sqrt{d})}\right\}$, then 
$\cF_2=\bigcup_d \Lambda_d$.
The union is disjoint, and the sets $\Lambda_d\subset (0,\infty)$ are countable, dense and consist of transcendental numbers.\\
As it is well known, the fundamental unit $\nu_d$ can be determined as follows: If $d\equiv 1\  (4)$, then $\nu_d=\frac{1}{2}(k+l\sqrt{d})$, where
the pair $(l,k)\in\NN^2$ is the smallest solution in lexicographic ordering of one of the equations $l^2d\mp 4 = k^2$ (Pell's equation).
If $d\equiv 2,3\  (4)$, then $\nu_d=k+l\sqrt{d}$, where $(l,k)$ solves 
$l^2d\mp 1 = k^2$ instead.

\item[$\cF_3':$] A cubic field is called complex if it has a complex embedding. For a complex cubic field $K$ we have $\cH_K\subset\RR^3_1=\{(\lambda,\lambda,-2\lambda)\mid \lambda\in\RR\}$. We define\\
$\Lambda_K:=\{\lambda\in (0,\infty)\mid (\pi\lambda,\pi\lambda,-2\pi\lambda)\in \cH_K\}$, and $\cF_3':=\bigcup \Lambda_K$,
where the union is taken over all complex cubic fields (up to isomorphism). The union is disjoint, and the sets $\Lambda_K\subset (0,\infty)$ are countable, dense and consist of transcendental numbers.

\item[$\widetilde\cF_3:$] For a real cubic field $K$ we define\\ 
$\widetilde\Lambda_K:=
\{0<\lambda_1<\lambda_2\mid$ a permutation of $(\pi\lambda_1,\pi\lambda_2,-\pi(\lambda_1+\lambda_2))$
belongs to $\cH_K\}$ and set
$\widetilde\cF_3:=
\bigcup \widetilde\Lambda_K$,
where the union is taken over all real cubic fields (up to isomorphism). The union is disjoint, and the sets $\widetilde\Lambda_K$ are countable, dense in the corresponding region of $\RR^2$ and consist of elements with transcendental coordinates.

\item[$\cF_3:$] For a real cubic field $K$ let
$\Lambda_K$ be as in Proposition \ref{duli}. We have\\
$\Lambda_K=\{\lambda_2/\lambda_1\mid (\lambda_1,\lambda_2)\in\widetilde\Lambda_K\}$. 
We set
$\cF_3:=
\bigcup \Lambda_K$,
where the union is taken over all real cubic fields (up to isomorphism). The sets $\Lambda_K\subset (1,\infty)$ are countable, dense and consist of transcendental numbers. The union is disjoint provided the four exponentials conjecture is true, see Prop.~\ref{duli}.

\item[$\cF_4:$] Let $s$ be a Salem number of degree $4$, see the paragraph preceeding Example~\ref{Excomp}. We choose $\rho\in\RR$ such that 
$e^{2\pi i\rho}$ is a Galois conjugate of $s$. Then we define\\
$\Lambda_s:=\left\{\displaystyle\frac{\log s}{2\pi |\rho + r|}\mid r\in\QQ\right\}$. The set does not depend on the choice of $\rho$. Moreover,
$\Lambda_s=\Lambda_{s^k}$. We set $\cF_4:=
\bigcup \Lambda_s$, where the union is taken over all Salem numbers of degree of $4$ (up to taking powers).
The sets $\Lambda_s\subset (0,\infty)$ are countable, dense and consist of transcendental numbers. If the four exponentials conjecture is true, then the above union is disjoint and $\cF_2\cap 2\cF_4=\emptyset$, cf.~Example~\ref{Ex23}.

\end{enumerate} 



\vspace{1.6cm}
{\footnotesize 

Ines Kath\\ 
Institut f\"ur Mathematik und Informatik\\ 
der Ernst-Moritz-Arndt Universit\"at Greifswald\\
Walther-Rathenau-Str.\,47, D-17489 Greifswald, Germany\\ 
email: ines.kath@uni-greifswald.de\\[2ex] 
Martin Olbrich\\ 
Universit\'e du Luxembourg \\
6, rue Richard Coudenhove-Kalergi, 
L-1359 Luxembourg\\
Luxembourg\\ 
email: martin.olbrich@uni.lu} 


\begin{thebibliography}{MMMM}

\bibitem[A]{Aus} Auslander, L., {\sl Bieberbach's theorem on space groups and discrete uniform subgroups of Lie groups.} Amer. J. Math. {\bf 83} (1961), 276-280.
\bibitem[Ba]{baum} Baum, H., {\sl Eichfeldtheorie: Eine Einf\"uhrung in die Differentialgeometrie auf Faserb\"undeln.} Springer 2009.

\bibitem[Be]{Be} Benoist, Y., {\sl Actions propres sur les espaces homog\`enes r\'eductifs.} Ann. of Math. {\bf 144} (1996), 315-347.

\bibitem[Bo]{Bo} Borel, A., {\sl Compact Clifford-Klein forms of symmetric spaces.}
 Topology  {\bf 2}  (1963), 111--122.

\bibitem[B1]{B1} Boyd, D.\,W., {\sl Salem numbers of degree four have periodic expansions.} Th\'eorie des nombres (Quebec, PQ, 1987), de Gruyter, Berlin, 1989, 57-64. 

\bibitem[B2]{B2} Boyd, D.\,W., {\sl On the Beta Expansion for Salem Numbers of degree 6.} Math. Comp. {\bf 65} (1996), 861-875. 

\bibitem[CP]{CP80} Cahen, M., Parker, M. {\sl Pseudo-Riemannian symmetric 
spaces.} 
Mem. Amer.\ Math.\ Soc.  {\bf 24} (1980), no. 229.

\bibitem[CW]{CW} Cahen, M., Wallach, N., {\sl Lorentzian symmetric spaces.}
Bull. Amer. Math. Soc.  \textbf{76} (1970), 585-591.

\bibitem[CM] {CM} Calabi, E., Markus, L., {\sl Relativistic space forms.} Ann. Math. {\bf 75}  (1962), 63-76.
\bibitem[Ca]{Ca} Carri\`ere, Y. {\sl Autour de la conjecture de L. Markus sur les vari\'et\'es affines.} Invent. Math. {\bf 95} (1989),  no. 3, 615-628.

\bibitem[Co]{Co} Cohen, H., {\sl A course in computational algebraic number theory.} 
Springer, 1996.

\bibitem[ELPV]{ELMV} Einsiedler, M., Lindenstrauss, E., Phillippe, M., and Venkatesh, A., {\sl Distribution of periodic torus orbits and Duke's theorem for cubic fields.}
Ann. Math. \textbf{173} (2011), 815-885.

\bibitem[F]{F} Fischer, M., {\sl Lattices of oscillator groups.} arXiv:1303.1970 [math.GR].

\bibitem[FG]{FG} Fried, D., Goldman, W,  {\sl Three-dimensional affine crystallographic groups} Adv. in Math. {\bf 47}, (1983), no. 1, 1-49.

\bibitem[FT]{FT} Fr\"ohlich, A., Taylor, M.~J., {\sl Algebraic number theory.} 
Cambridge University Press, 1993.

\bibitem[Gh]{Gh} Ghys, \'E., {Flots d'Anosov dont les feuilletages stables sont diff\'erentiables.} Ann. Sci. \'Ecole Norm. Sup. (4), {\bf 20}, no. 2 (1987), 251-270.

\bibitem[G1]{Gads} Goldman, W., {\sl Nonstandard Lorentz space forms.} J. Differential Geom. {\bf 21} (1985), 301-308. 

\bibitem[G2]{Gicm} Goldman, W., {\sl Locally homogeneous geometric manifolds} Proceedings of the 2010 International Congress of Mathematicians, Hyderabad, India (2010), 717-744, Hindustan Book Agency, New Delhi.

\bibitem[GK]{GK} Goldman, W., Kamishima, Y., {\sl The fundamental group of a compact flat Lorentz space form is virtually policyclic.} J. Differential Geom. {\bf 19}  (1984), 233-240.

\bibitem[GW]{GW} Gordon, C. S., Wilson, E. N. {\sl The spectrum of the Laplacian on Riemannian Heisenberg manifolds.}  Mich. Math. J. {\bf 33} (1986),  253-271.

\bibitem[GS]{GS77} V. Guillemin and S. Sternberg,  {\em Geometric asymptotics.} AMS, 1977.
\bibitem[H]{Hall} Hall, Marshall, Jr, {\sl A topology for free groups and related groups.} Ann. of Math. {\bf 52} (1950), 127-139.
\bibitem[KO1]{KO1} Kath, I., Olbrich, M.,
{\sl On the structure of pseudo-Riemannian symmetric spaces.} Transform.\,Groups {\bf 14} (2009), no.~4, 847-885.
\bibitem[KO2]{KOesi} Kath, I., Olbrich, M., {\sl The classification problem for 
pseudo-Riemannian symmetric spaces} in Recent developments in 
pseudo-Riemannian Geometry,
ESI Lectures in Mathematics and Physics, EMS Publishing House, 2008, 1-52. 
\bibitem[Ki]{Ki} Kirby, D., {\sl Integer matrices of finite order.} Rend. Mat. {\bf 6} (1969), no. 2, 403-408.

\bibitem[Kl] {K} Klingler, B., {\sl Completude des vari\'et\'es lorentziennes \`a courbure constante.} Math. Ann. {\bf 306} (1996), 353-370.
\bibitem[Ko1]{Ko1} Kobayashi, T., {\sl Proper action on a homogeneous space of reductive type.} Math. Ann. {\bf 285} (1989) 249-263.
\bibitem[Ko2]{Ko2} Kobayashi, T., {\sl A necessary condition for the existence of compact Clifford-Klein forms of homogeneous spaces of reductive type.} Duke Math. J. {\bf 67} (1992) 653-664.

\bibitem[KY] {KY} Kobayashi, T., Yoshino, T., {\sl Compact Clifford-Klein forms of symmetric spaces -- revisited.} Pure and Appl. Math. Quaterly {\bf 1} (2005), 603-684. 

\bibitem[Ku]{Ku}Kulkarni, R., {\sl Proper actions and pseudo-Riemannian space forms.} Adv. Math. {\bf 40} (1981), 10-51.

\bibitem[LS]{LS} Leistner, T., Schliebner, D., {\sl Completeness of compact Lorentzian manifolds with special holonomy.} arXiv:1306.0120v2 [math.DG]
\bibitem[M]{M} Medina, A., {\sl Groupes de Lie munis de m\'etriques bi-invariantes.} Tohoku 
Math. J. (2) {\bf 37} (1985), no. 4, 405-421. 
\bibitem[MR1]{MR1} Medina, A., Revoy, P., {\sl Alg\`ebres de Lie et produit scalaire invariant.} 
Ann. Sci. \`Ecole Norm. Sup. (4) {\bf 18} (1985), no. 3, 553--561.
\bibitem[MR2]{MR2} Medina, A., Revoy, P., {\sl Les groups oscillateurs et leurs reseaux.} Manuscripta Mathematica {\bf 52} (1985), 81-95. 
\bibitem[N]{Neu} Neukirch, J., {\sl Algebraische Zahlentheorie.} 
Springer, 1992.
\bibitem[Nr]{Nr} Neukirchner, Th.,  {\sl Solvable Pseudo-Riemannian 
    Symmetric Spaces.} arXiv:math/0301326 [math.DG], 2003.
\bibitem[Ra]{Ragh} Raghunathan, M.~S., {\sl Discrete subgroups of Lie groups.} 
Springer, 1972.
\bibitem[Ri]{Ri} Rivin, I., {\sl Walks on groups, counting reducible matrices, polynomials, and surface and free group automorphisms.} Duke Math. J. {\bf 142} (2008), no. 2, 353-379. 
\bibitem[S]{Sch} Schneider, Th., {\sl Einf\"uhrung in die transzendenten Zahlen.} 
Springer, 1957.
\bibitem[T]{T} Tolimieri, R., {\sl Heisenberg manifolds and theta functions.} Trans. Amer. Math. {\bf 239} (1978), 293-319.
\bibitem[Wa]{Wa} Waldschmidt, M., {\sl Diophantine approximation on linear algebraic groups.} 
Springer, 2000.
\bibitem[Wi1]{WM1} Witte, D., {\sl Superrigidity of lattices in solvable Lie groups.}
 Invent. Math.  {\bf 122}  (1995),  no. 1, 147--193.

\bibitem[Wi2]{WM2} Witte, D., {\sl Superrigid subgroups and syndetic hulls in solvable Lie groups.}
 Rigidity in dynamics and geometry (Cambridge, 2000), 
 441--457, Springer, 2002. 
\bibitem[Z]{Zads} Zeghib, A., {\sl On closed anti de Sitter spacetimes}, Math. Ann. {\bf 310} (1998), 695-716. 

\end{thebibliography}
\end{document}